\documentclass[a4paper,oneside,12pt]{amsart}

\title{The index of mesoscopic Wilson--Dirac operators: I}
\author[M. Furuta]{Mikio Furuta}
\author[S. Matsuo]{Shinichiroh Matsuo \orcidlink{0000-0003-1240-6924}}
\address[MF]{%
  Graduate School of Mathematical Sciences \\
  The University of Tokyo \\
  Tokyo, Japan
}
\thanks{Mikio Furuta passed away on 28 June 2026.}
\address[SM]{%
  Department of Mathematical and Computing Science \\
  School of Computing \\
  Institute of Science Tokyo \\
  Tokyo, Japan
}
\email{matsuo.shinichiroh@comp.isct.ac.jp}

\subjclass[2020]{Primary 58J20; Secondary 19K56, 58J30, 81T25}
\keywords{Dirac operator, index, spectral flow, Wilson--Dirac operator, lattice gauge theory}

\usepackage{bxpapersize}
\calclayout

\usepackage[utf8]{inputenc}
\usepackage[T1]{fontenc}

\usepackage[p,osf]{cochineal}
\usepackage[scale=.95,type1]{cabin}
\usepackage[cochineal,vvarbb]{newtxmath}
\usepackage[cal=boondoxo]{mathalfa}

\usepackage{microtype}

\usepackage{mathtools}
\mathtoolsset{showonlyrefs=true}

\theoremstyle{plain}
\newtheorem{theorem}{Theorem}[section]
\newtheorem{proposition}[theorem]{Proposition}
\newtheorem{lemma}[theorem]{Lemma}
\newtheorem{corollary}[theorem]{Corollary}
\newtheorem*{maintheorem}{Main Theorem}

\theoremstyle{definition}
\newtheorem{definition}[theorem]{Definition}

\theoremstyle{remark}
\newtheorem{remark}[theorem]{Remark}

\usepackage{orcidlink}

\usepackage{tikz}
\usetikzlibrary{calc}

\usepackage{hyperref}
\hypersetup{
  bookmarksnumbered=true,
  hidelinks,
  pdftitle={The index of mesoscopic Wilson--Dirac operators: I},
  pdfauthor={Mikio Furuta and Shinichiroh Matsuo}
}

\usepackage[msc-links]{amsrefs}

\newcommand{\Z}{\mathbb{Z}}
\newcommand{\R}{\mathbb{R}}
\newcommand{\C}{\mathbb{C}}

\DeclarePairedDelimiter{\abs}{\lvert}{\rvert}
\DeclarePairedDelimiter{\norm}{\lVert}{\rVert}
\DeclarePairedDelimiter{\seminorm}{\lbrack}{\rbrack}
\DeclarePairedDelimiterX{\innerproduct}[2]{\langle}{\rangle}{#1,#2}
\DeclarePairedDelimiterX{\pointinnerproduct}[2]{\langle}{\rangle}{#1,#2}
\DeclarePairedDelimiterX{\commutator}[2]{\lbrack}{\rbrack}{#1,#2}

\DeclareMathOperator{\End}{End}
\DeclareMathOperator{\Hom}{Hom}
\DeclareMathOperator{\id}{id}

\DeclareMathOperator{\trace}{tr}

\DeclareMathOperator{\kernel}{Ker}
\DeclareMathOperator{\image}{Im}
\DeclareMathOperator{\cokernel}{Coker}

\DeclareMathOperator{\support}{supp}

\newcommand{\essentialspectrum}{\operatorname{spec}_{\mathrm{ess}}}

\DeclareMathOperator{\indexAtiyahSinger}{ind}
\DeclareMathOperator{\spectralflow}{sf}

\DeclareMathOperator{\injectivityradius}{inj}
\newcommand{\dvol}{d\mu}
\newcommand{\volume}{\mathrm{vol}}

\newcommand{\paralleltransport}[2]{\operatorname{P}_{#2 \gets #1}}
\newcommand{\paralleltransportLC}[2]{\paralleltransport{#1}{#2}^{\mathrm{LC}}}

\newcommand{\flattorus}{\mathbb{T}^n}
\newcommand{\Lattice}{L}
\newcommand{\duallattice}[1]{#1^{\vee}}

\newcommand{\cliffordmultiplication}{\rho}

\newcommand{\Ltwonorm}[2][L^2]{\norm{#2}_{#1}}
\DeclarePairedDelimiter{\operatornorm}{\lVert}{\rVert_{\mathrm{op}}}
\newcommand{\mesoscopicseminorm}[2][a]{\seminorm{#2}_{#1}}
\newcommand{\mesoscopicdifferenceseminorm}[2][a]{\seminorm{#2}_{\nabla,#1}}
\newcommand{\mesoscopicdifferenceseminormLocal}[3][a]{\seminorm{#3}_{\nabla,#1;#2}}
\newcommand{\flatmesoscopicdifferenceseminorm}[3][a]{\seminorm{#3}_{#1,#2,\mathrm{fl}}}

\newcommand{\mesoscopicDiracoperator}[1][a]{\mathcal{D}_{#1}}
\newcommand{\mesoscopicWilsonterm}[1][a]{\mathcal{W}_{#1}}
\newcommand{\mesoscopicWilsonDiracoperator}[1][a]{\mathcal{D}_{#1}^{\mathit{WD}}}

\newcommand{\cutoff}{\beta}
\newcommand{\cutoffDirac}{\cutoff^D}
\newcommand{\cutoffWilson}{\cutoff^W}
\newcommand{\cutoffnabla}{\cutoff^{\nabla}}
\newcommand{\cutoffupper}{\chi^+}
\newcommand{\cutofflower}{\chi^-}

\newcommand{\mesoscopicWilsontermcompactpert}[1][a]{\mathcal{B}_{#1}}

\newcommand{\integralkernel}{k}
\newcommand{\kernelDirac}{\integralkernel_a^D}
\newcommand{\kernelWilson}{\integralkernel_a^W}
\newcommand{\symbolDirac}{\sigma^D}
\newcommand{\symbolWilson}{\sigma^W}
\newcommand{\symbolWilsonDirac}{\sigma^{\mathit{WD}}}
\newcommand{\symbolnabla}{\sigma^{\nabla}}
\newcommand{\radialsymbolDirac}{\symbolDirac_{\mathrm{rad}}}
\newcommand{\radialsymbolWilson}{\symbolWilson_{\mathrm{rad}}}

\newcommand{\flatmesoscopicDiracoperator}{\mesoscopicDiracoperator^{\mathrm{fl}}}
\newcommand{\flatmesoscopicWilsonterm}{\mesoscopicWilsonterm^{\mathrm{fl}}}
\newcommand{\flatmesoscopicWilsonDiracoperator}[1][a]{\mathcal{D}_{#1}^{\mathit{WD},\mathrm{fl}}}

\newcommand{\flatmodel}[1]{#1^{\mathrm{fl}}}

\newcommand{\coupledoperator}[3][a]{\mathcal{C}_{#1}(#2,#3)}

\newcommand{\voronoicell}{\mathcal{V}}

\newcommand{\constantlocal}{C_{\mathrm{loc}}}
\newcommand{\constantflat}{C_{\mathrm{fl}}}
\newcommand{\constantgeometric}{C_{\mathrm{geom}}}
\newcommand{\constantLtwonorm}{C_{\mathrm{norm}}}
\newcommand{\constantnabla}{C_{\nabla}}
\newcommand{\constantWD}{C_{\mathit{WD}}}
\newcommand{\constantcommutator}{C_{\mathrm{com}}}
\newcommand{\radiusPoincare}{r_{\mathrm{Poin}}}
\newcommand{\constantPoincare}{C_{\mathrm{Poin}}}
\newcommand{\constantcomparison}{C_{\mathrm{F}}}
\newcommand{\constantmult}{C_{\mathrm{mult}}}

\newcommand{\constantAell}{a_{\mathrm{ell}}}
\newcommand{\constantAvan}{a_{\mathrm{van}}}
\newcommand{\constantAindex}{a_{\mathrm{ind}}}

\definecolor{ink}{HTML}{243347}
\definecolor{ballline}{HTML}{536C87}
\definecolor{ballouter}{HTML}{F7F9FC}
\definecolor{ballinner}{HTML}{E7EEF6}
\definecolor{shellblue}{HTML}{B9DCEE}
\definecolor{shellgold}{HTML}{EFC46D}
\definecolor{blueedge}{HTML}{2479B4}
\definecolor{goldedge}{HTML}{B77712}
\definecolor{bluewash}{HTML}{EDF6FB}
\definecolor{goldwash}{HTML}{FFF3D9}
\definecolor{divider}{HTML}{E3E8EE}

\begin{document}
\begin{abstract}
  This is the first of a series of papers which will be devoted to a study of the index of \emph{mesoscopic} Wilson--Dirac operators on closed Riemannian manifolds.
  The aim of the series is to discretise the index of the Dirac operator.
  On flat tori, this is achieved by the Wilson--Dirac operator of lattice gauge theory on a regular lattice.
  On a general closed Riemannian manifold, however, no regular lattice is available.
  As an intermediate step between the continuum Dirac operator and its discretisations, we introduce mesoscopic Wilson--Dirac operators on any closed Riemannian manifold, defined by integrating covariant differences over distances comparable to a fixed length scale.
  In this paper, we define the index of the mesoscopic Wilson--Dirac operators via spectral flow and prove that it agrees with that of the continuum Dirac operator at all sufficiently small scales.
  The proof rests on an elliptic estimate uniform in the scale and a Rellich-type compactness theorem.
  We also prove a mesoscopic version of the Nielsen--Ninomiya theorem.
  In the second paper, we shall pass from mesoscopic operators to discrete Wilson--Dirac operators.
\end{abstract}
\maketitle

\setcounter{tocdepth}{3}
\tableofcontents

\section{Introduction}

In this paper, we introduce mesoscopic Wilson--Dirac operators on closed Riemannian manifolds and study their index.
The purpose of introducing these operators is to represent a Dirac operator and its index at a finite length scale.
A Dirac operator is a differential operator, and a derivative is the limit of difference quotients as the distance tends to zero.
Instead of taking this limit, we integrate covariant difference quotients over distances comparable to a prescribed length scale.
This scale is small compared with the geometry of the manifold, which we regard as the macroscopic scale, but it is not infinitesimal, in contrast with the microscopic scale of differential operators.
We shall show that the index, a macroscopic invariant, is already captured at the mesoscopic scale.
The present paper is a first step towards our long-term goal of a discretisation of Seiberg--Witten theory that keeps the spatial locality of the equations.
We explain this motivation later in this introduction.

For flat tori, we showed in earlier joint work~\cite{MR5076866} that the index of the continuum Dirac operator can be recovered from a discrete operator on a sufficiently fine regular lattice.
The essential ingredient is the Wilson--Dirac operator from lattice gauge theory, a finite-difference operator whose index is defined through the spectral flow of a mass-dependent family.
Our next aim was to extend this construction and the index comparison to general closed Riemannian manifolds.
However, it was not clear how to construct a suitable discrete Wilson--Dirac operator without the regular lattice structure available on a flat torus.

We therefore introduce the mesoscopic Wilson--Dirac operators as an intermediate class of operators.
The mesoscopic operators are bounded on the same $L^2$ space as the unbounded continuum Dirac operator.
The index comparison follows the same strategy as in the flat torus case, but the mesoscopic formulation substantially simplifies the analytic argument.
In the second paper of this series, we shall pass from mesoscopic operators to discrete Wilson--Dirac operators by replacing the integrals with suitable weighted Riemann sums over finite subsets of the manifold.
This point-sampling approach to integration is reminiscent of Monte Carlo methods.
It will yield discrete Wilson--Dirac operators on general closed Riemannian manifolds, and we shall prove that the resulting discrete operators have the same index as the corresponding mesoscopic Wilson--Dirac operators when the approximation is sufficiently fine.
By the index theorem proved in the present paper, their index therefore also agrees with that of the continuum Dirac operator.
The mesoscopic formulation separates the index comparison from the discretisation problem and keeps the interaction scale independent of the sampling resolution.

We now describe the mesoscopic Wilson--Dirac operators and the main result in more detail.
Let $X$ be a closed Riemannian manifold, and let $E$ be a graded Hermitian Clifford module bundle over $X$ with grading operator $\gamma$.
Fix a compatible Clifford connection on $E$.
The associated Dirac operator $D$ is an unbounded self-adjoint Fredholm operator on $L^2(X,E)$.
For a sufficiently small length scale $a>0$, we define a mesoscopic Dirac operator $\mesoscopicDiracoperator$ in Section~\ref{subsection: the mesoscopic Dirac operator} and a mesoscopic Wilson term $\mesoscopicWilsonterm$ in Section~\ref{subsection: the mesoscopic Wilson term}.
They are given by
\begin{equation}
  (\mesoscopicDiracoperator \Phi)(x) = \int_{TX_x} \cliffordmultiplication_x\left(\frac{v}{\abs{v}}\right) \frac{\paralleltransport{\exp_x(v)}{x}(\Phi(\exp_x(v))) - \Phi(x)}{\abs{v}} \frac{1}{C_D} \cutoffDirac\left(\frac{\abs{v}}{a}\right)\,\frac{dv}{a^n}
\end{equation}
and
\begin{equation}
  (\mesoscopicWilsonterm \Phi)(x) = \int_{TX_x} \frac{\Phi(x) - \paralleltransport{\exp_x(v)}{x}(\Phi(\exp_x(v)))}{\abs{v}} \frac{1}{C_W} \cutoffWilson\left(\frac{\abs{v}}{a}\right)\,\frac{dv}{a^n}
\end{equation}
for $\Phi \in C^{\infty}(X,E)$ and $x \in X$.
Here, $n$ is the dimension of $X$, $\cliffordmultiplication$ is the Clifford multiplication, and $\paralleltransport{y}{x} \colon E_y \to E_x$ is the parallel transport along the minimising geodesic from $y$ to $x$.
The cut-off functions $\cutoffDirac$ and $\cutoffWilson$ are supported in compact subsets of $(0,1)$, and $C_D$ and $C_W$ are normalisation constants.
The combination
\begin{equation}
  \mesoscopicWilsonDiracoperator = \mesoscopicDiracoperator + \mesoscopicWilsonterm\gamma
\end{equation}
is a bounded self-adjoint Fredholm operator on $L^2(X,E)$, which we call the \emph{mesoscopic Wilson--Dirac operator}.
Its operator norm is of order $a^{-1}$, whereas
\begin{equation}
  \mesoscopicWilsonDiracoperator \Phi \longrightarrow D\Phi
\end{equation}
in $L^2(X,E)$ for every fixed smooth section $\Phi$ as $a \to 0$.
Thus, the unbounded Dirac operator $D$ is studied through a family $\{\mesoscopicWilsonDiracoperator\}_{a > 0}$ of bounded operators with the differential behaviour recovered as the scale $a$ tends to zero.
The main result of this paper is the following index theorem.

\begin{maintheorem}[Theorem~\ref{theorem: main theorem}]
  Fix the geometric data, the cut-off functions, and the normalisation constants in the construction of Section~\ref{section: mesoscopic WD operators}.
  There exists $\constantAindex > 0$ such that the family $\{\mesoscopicWilsonDiracoperator+m\gamma\}_{m\in[-1,1]}$ consists of self-adjoint Fredholm operators, is invertible at $m = \pm 1$, and has spectral flow equal to the index of $D$:
  \begin{equation}
    \spectralflow \bigl(\{\mesoscopicWilsonDiracoperator+m\gamma\}_{m\in[-1,1]}\bigr) = \indexAtiyahSinger D
  \end{equation}
  for every $0 < a < \constantAindex$.
\end{maintheorem}

The Wilson term has an essential role in this construction.
In lattice gauge theory, a symmetric-difference approximation of the Dirac operator produces additional low-energy modes at large momenta known as fermion doublers.
From a mathematical viewpoint, fermion doubling means that the naive symmetric-difference approximation of the Dirac operator does not satisfy an elliptic estimate uniform in the lattice spacing.
Adding the Wilson term restores such an estimate for the resulting Wilson--Dirac operator.
See \cite{MR5076866}*{Theorem~4.7}.
In physical terms, the Wilson term gives the doubler modes masses of order the inverse lattice spacing, while its effect on modes of fixed small momentum vanishes in the continuum limit.
In the mesoscopic setting, the Dirac part $\mesoscopicDiracoperator$ is compact, and hence is not Fredholm on $L^2(X,E)$.
In the flat model, its symbol tends to zero at high frequencies, so that all sufficiently high frequencies behave like fermion doublers.
Thus, fermion doubling already occurs at a finite length scale without a lattice.
The Wilson term provides control at high frequencies.
The resulting mesoscopic Wilson--Dirac operator is Fredholm and satisfies uniform elliptic estimates.
At low frequencies, it still approximates the Dirac operator.

The Dirac operator $D$ and the mesoscopic Dirac operator $\mesoscopicDiracoperator$ anti-commute with $\gamma$, whereas the mesoscopic Wilson--Dirac operator does not.
This is unavoidable.
We prove a mesoscopic version of the Nielsen--Ninomiya theorem~\cites{MR0621425, MR0639112, NNnogo} in lattice gauge theory (Theorem~\ref{theorem: no-go theorem}).
It says that, if a self-adjoint Fredholm operator is the sum of a multiplication operator and a compact operator and anti-commutes with $\gamma$, then the trace of $\gamma$ on its kernel is zero.
For $D$, this trace is the index.
Thus, we define the index of the mesoscopic Wilson--Dirac operator by spectral flow, as for the lattice Wilson--Dirac operator.

The proof of the main theorem follows the strategy used for lattice Wilson--Dirac operators on flat tori in~\cite{MR5076866}.
In that work, the continuum and lattice operators act on different Hilbert spaces, and their comparison requires interpolation maps together with estimates for these maps and their interaction with the operators.
The mesoscopic operators substantially simplify this part of the argument.
Both mesoscopic and continuum operators now act on the same $L^2$ space, so we need no interpolation maps to compare them.
The analytical basis of the comparison is an elliptic estimate uniform in $a$.
Writing $\mesoscopicdifferenceseminorm{\Phi}$ for the covariant mesoscopic difference seminorm introduced in Definition~\ref{definition: mesoscopic difference seminorm on manifolds}, we prove
\begin{equation}
  \mesoscopicdifferenceseminorm{\Phi} \le C \bigl(\Ltwonorm{\mesoscopicWilsonDiracoperator \Phi} + \Ltwonorm{\Phi}\bigr)
\end{equation}
for all sufficiently small $a$ with $C$ independent of $a$ and $\Phi$.
The mesoscopic Rellich lemma (Theorem~\ref{theorem: mesoscopic Rellich lemma}) then gives compactness as the scale tends to zero.
If $a_j \to 0$ and both $\Phi_j$ and $\mesoscopicWilsonDiracoperator[a_j]\Phi_j$ are bounded in $L^2$, a subsequence of $\Phi_j$ converges strongly in $L^2$ to a section $\Phi \in L^2_1$, and $\mesoscopicWilsonDiracoperator[a_j] \Phi_j$ converges weakly to $D \Phi$ along the same subsequence.
The proof of the main theorem involves sequences of sections that vary with the scale.
Convergence of the operators on each fixed smooth section alone is not sufficient for our argument.
Using this compactness, we show that, for small $a$, the families $\{D+m\gamma\}_{m \in [-1,1]}$ and $\{\mesoscopicWilsonDiracoperator+m\gamma\}_{m \in [-1,1]}$ define the same class in $K^1([-1,1],\{\pm 1\})$.
Since $\spectralflow (\{D+m\gamma\}) = \indexAtiyahSinger D$, the theorem follows.

We now discuss how our results relate to earlier works.
Index theorems for lattice Wilson--Dirac operators on flat tori have been proved by several methods.
Adams~\cite{MR1866670} computed the continuum limit of the local index density of the overlap Dirac operator.
Yamashita~\cite{MR4275791} used deformation quantisation and an algebraic index theorem, while Kubota~\cite{MR4407739} used higher index theory of almost flat bundles.
Aoki, Fukaya, Furuta, Matsuo, Onogi and Yamaguchi~\cite{MR5076866} gave a direct comparison of lattice and continuum mass families in relative $K$-theory.

Nonlocal difference energies and their relation to Sobolev regularity were studied by Bourgain, Brezis and Mironescu~\cite{MR3586796}.
Ponce~\cite{MR2041005} established Poincar\'e inequalities and compactness results for concentrating radial kernels.
After localisation and comparison in local trivialisations, his compactness results apply to the shell kernels used in this paper and yield the mesoscopic Rellich lemma.
We give two self-contained proofs adapted to covariant differences on the manifold: one uses local averages and a Poincar\'e inequality for differences over shells, and the other uses high-frequency estimates.

The use of a nonlocal intermediate operator also has precedents in geometric discretisation.
Burago, Ivanov and Kurylev developed graph approximations of the Laplace--Beltrami operator~\cite{MR3299811} and studied convolution-type Laplacians on metric-measure spaces~\cite{MR3990939}.
Burago, Ivanov, Kurylev and Lu~\cite{MR4437353} extended this approach to vector bundles with connection.
Their construction separates the interaction radius from the mesh size and uses parallel transport and volume-weighted sums to compare nonlocal and discrete connection Laplacians with the continuum operator.
The local averaging projection in our compactness proof is obtained by averaging a section over each cell and then extending the averages back to a section, using parallel transport in both steps.
For the same partition, these two steps are precisely their averaging and extension maps.
These works concern Laplacians and the convergence of their spectra.
Their nonlocal Laplacians and our mesoscopic Dirac operator both integrate the differences between the value of a section at a point and its values near the point.
However, there is an essential difference between them.
A second-order derivative does not change sign when the direction is reversed, whereas a first-order derivative does.
Thus, a nonlocal Laplacian adds the differences in opposite directions, whereas the mesoscopic Dirac operator subtracts them.
Hence, for a Laplacian, the term involving the value at the point remains, and it provides control at high frequencies.
For the Dirac operator, this term vanishes, as in the symmetric difference on a lattice, and $\mesoscopicDiracoperator$ is compact.
This is known as fermion doubling in lattice gauge theory.
It comes from the fact that the Dirac operator is of first order, and it does not occur for Laplacians.
The mesoscopic Wilson term is the scale $a$ times a nonlocal connection Laplacian, and it restores the term involving the value at the point.

The high-frequency obstruction to naive lattice approximations and the role of Wilson-type corrections have also been studied through continuum limits of discrete Dirac operators.
Fermion doubling prevents the naive symmetric-difference approximation of the Dirac operator from converging to the continuum operator in the norm-resolvent sense.
Cornean, Garde and Jensen~\cite{MR4590014} and Nakamura~\cite{MR4741062} proved that adding suitable Wilson-type corrections gives this convergence on Euclidean lattices.
To compare the discrete and continuum operators, they use suitable embeddings of the discrete spaces into the continuum $L^2$ spaces.
In the mesoscopic setting, no such embeddings are needed.
The mesoscopic Wilson--Dirac operators converge to $D$ in the norm-resolvent sense, whereas the mesoscopic Dirac operators do not.
See Corollary~\ref{corollary: norm-resolvent convergence} and Remark~\ref{remark: norm-resolvent nonconvergence}.

We now explain in more detail the motivation for introducing the mesoscopic Wilson--Dirac operators.
Donaldson and Seiberg--Witten invariants are defined using moduli spaces of solutions to gauge-theoretic partial differential equations on smooth four-manifolds.
In dimension four, a PL structure determines a compatible smooth structure up to diffeomorphism.
This correspondence motivates the search for descriptions of these invariants directly in terms of finite combinatorial data associated with the underlying PL structure.
One way to approach this problem is to discretise the defining equations.
Our long-term goal is a discretisation of Seiberg--Witten theory compatible with finite-dimensional approximation of the Seiberg--Witten map.
The discretisation may lead to a simpler description of gluing formulae~\cite{MR2312050} for Bauer and Furuta's stable cohomotopy refinement of the Seiberg--Witten invariants~\cite{MR2025298}.
We hope to construct finite-dimensional models in which gluing along a common boundary can be described by coupling the models on the two pieces through variables near that boundary.
To this end, the discretisation should keep the spatial locality of the original equations.
The treatment of the nonlinear Seiberg--Witten map, gauge symmetry, and boundary data remains part of this longer-term programme.
The present paper treats a linear problem towards this goal: how to represent a Dirac operator and its index at a finite length scale.

We seek a description that uses differences at a finite length scale while keeping the spatial locality of the original equations.
Accordingly, the covariant differences involve only pairs of points whose distances lie between two fixed positive multiples of that scale.
We refer to the upper bound on these distances as \emph{IR-locality} (infrared locality), and to the positive lower bound as \emph{UV-locality} (ultraviolet locality).
IR-locality confines interactions across a separating hypersurface to a neighbourhood of controlled width, as desired for applications to gluing.
UV-locality keeps the integral kernels used in the construction away from the diagonal at each fixed scale.
The resulting operators are defined directly on the manifold without choosing a lattice.
This leaves the interaction scale independent of the resolution of a subsequent discretisation.

The scale dependence of mesoscopic operators also suggests a connection with renormalisation.
In the flat model, the symbol of the mesoscopic Wilson--Dirac operators has the form $\sigma(a\xi)/a$.
For small dimensionless momentum $a\xi$, it agrees with the symbol of the Dirac operator to leading order, while at large dimensionless momentum $\sigma$ approaches a positive multiple of the grading operator.
The construction modifies the short-distance description while preserving the low-energy behaviour and the index.
We regard this as a linear model for a topological form of renormalisation, in which the information to be kept is specified by an index or, more generally, a homotopy-theoretic invariant.

The paper is organised as follows.
Section~\ref{section: mesoscopic WD operators} gives the construction of the mesoscopic operators and proves a no-go theorem of Nielsen--Ninomiya type.
Sections~\ref{section: Elliptic estimates on the flat torus} and~\ref{section: Elliptic estimates on closed manifolds} establish elliptic estimates, first on a flat torus and then on a closed Riemannian manifold.
Section~\ref{section: the mesoscopic rellich lemma} proves the mesoscopic Rellich lemma.
Section~\ref{section: the coupled operators} proves the technical core of the paper.
Section~\ref{section: Spectral flow and the index} recalls the relation between relative $K^1$ theory, spectral flow and the index of Dirac operators.
Section~\ref{section: index of mesoscopic Wilson Dirac operators} proves the main theorem.

\section{Mesoscopic Wilson--Dirac operators}\label{section: mesoscopic WD operators}

In this section we introduce the mesoscopic Wilson--Dirac operator associated with a graded Clifford module bundle over a closed Riemannian manifold.
We first fix notation.
We then define the mesoscopic Dirac operator, the Wilson term, and the mesoscopic Wilson--Dirac operator.
We finally prove a no-go theorem of Nielsen--Ninomiya type.

\subsection{Notation}

Let $X=(X,g)$ be a closed Riemannian manifold.
Set $n := \dim X$.
Let $\injectivityradius(X,g)$ be the injectivity radius of $(X,g)$.
We denote by $d_g$ the distance function on $X$ induced by $g$.
We denote by $\dvol_g$ the Riemannian measure, and by $\volume(X,g)$ the volume of $(X,g)$.
For $x,y \in X$ with $d_g(x,y) < \injectivityradius(X,g)$, we denote by
\begin{equation}
  \paralleltransportLC{y}{x} \colon TX_y \to TX_x
\end{equation}
the parallel transport with respect to the Levi-Civita connection along the unique minimising geodesic from $y$ to $x$.
We denote by $B(p,r)$ the open ball around $p \in X$ with radius $r > 0$.

A graded Clifford module bundle $E$ over $X$ is a Hermitian vector bundle with a grading operator $\gamma$ such that $\gamma^2 = \id$ and $\gamma^* = \gamma$, together with a homomorphism $\cliffordmultiplication \colon TX \to \End(E)$ such that
\begin{gather}
  \gamma\cliffordmultiplication_x(v)+\cliffordmultiplication_x(v)\gamma=0, \\
  \cliffordmultiplication_x(u)\cliffordmultiplication_x(v) +\cliffordmultiplication_x(v)\cliffordmultiplication_x(u) = -2g_x(u,v), \\
  \cliffordmultiplication_x(v)^*=-\cliffordmultiplication_x(v)
\end{gather}
for all $u,v \in TX_x$ and $x \in X$.
We formulate the Clifford multiplication as a homomorphism from $TX$, rather than from $T^*X$, for later convenience.

Fix a Clifford connection $\nabla$ on $E$.
Thus, $\nabla$ is a Hermitian connection on $E$ with $\nabla \cliffordmultiplication = 0$ and $\nabla \gamma = 0$.
The Dirac operator $D \colon C^{\infty}(X,E) \to C^{\infty}(X,E)$ is defined by
\begin{equation}
  D \Phi := \sum_{j=1}^n \cliffordmultiplication(e_j) \nabla_{e_j} \Phi
\end{equation}
for $\Phi \in C^{\infty}(X,E)$, where $\{e_1, \dots, e_n\}$ is a local orthonormal frame of $TX$.
The Dirac operator is self-adjoint and anti-commutes with $\gamma$.
We sometimes call $D$ the \emph{continuum} Dirac operator, in contrast with the \emph{mesoscopic} Dirac operator introduced below.

For $x,y \in X$ with $d_g(x,y) < \injectivityradius(X,g)$, we denote by
\begin{equation}
  \paralleltransport{y}{x} \colon E_y \to E_x
\end{equation}
the parallel transport with respect to $\nabla$ along the unique minimising geodesic from $y$ to $x$.
We have $\gamma_y \circ \paralleltransport{x}{y} = \paralleltransport{x}{y} \circ \gamma_x$.

For $x \in X$, we define the Jacobian density $J_x \colon TX_x \to [0,\infty)$ by
\begin{equation}
  (\exp_x)^* \dvol_g = J_x(v) dv
\end{equation}
where $dv$ is the Lebesgue measure on $TX_x$.
We have $J_x(v) = \abs{\det (d\exp_x)_v}$ under the canonical identification $T(TX_x)_v = TX_x$.

For $x,y \in X$ with $d_g(x,y) < \injectivityradius(X,g)$, we set $v_{xy} := \exp_x^{-1}(y) \in TX_x$.
We have $\abs{v_{xy}} = d_g(x,y)$.
Reversing the geodesic, we have
  \begin{equation}\label{equation: reversing the geodesic}
    v_{xy} =\paralleltransportLC{y}{x}\big(\paralleltransportLC{x}{y}(\exp_x^{-1}(y))\big) = \paralleltransportLC{y}{x}\big(-\exp_y^{-1}(x)\big) = \paralleltransportLC{y}{x}(-v_{yx}).
  \end{equation}
\begin{lemma}\label{lemma: symmetry of Jacobian density}
  For $x,y \in X$ with $d_g(x,y) < \injectivityradius(X,g)$, we have $J_x(v_{xy}) = J_y(v_{yx})$.
\end{lemma}
\begin{proof}
  Since $J_x(v_{xy}) = \abs{\det (d \exp_x^{-1})_y}^{-1}$ and $J_y(v_{yx}) = \abs{\det (d \exp_y^{-1})_x}^{-1}$, and since
  \begin{equation}
    g_x \big( (d \exp_x^{-1})_y (u_y),u_x \big) = g_y \big( u_y, (d \exp_y^{-1})_x (u_x) \big)
  \end{equation}
  for $u_x \in TX_x$ and $u_y \in TX_y$ by the symmetry of the mixed second derivatives of $d_g(x,y)^2/2$, we have $J_x(v_{xy}) = J_y(v_{yx})$.
\end{proof}

\subsection{The mesoscopic scale and cut-off functions}
We fix the geometric scale and cut-off functions that will be used in the definition of the mesoscopic operators.

Let $a\in\big(0,\injectivityradius(X,g)/10\big)$.
All constructions in this section are made at the scale $a$.
The parameter $a$ should be regarded as a small \emph{mesoscopic} length scale relative to the \emph{macroscopic} scale $\injectivityradius(X,g)$, while $a^{-1}$ should be regarded as the associated mesoscopic \emph{energy} scale.
The limit $a \to 0$ corresponds to the microscopic regime of the differential operator $D$.

Fix $0 < r_-^D < r_+^D < 1$ and $0 < r_-^W < r_+^W < 1$.
Let $\cutoffDirac, \cutoffWilson \colon [0,\infty) \to [0,\infty)$ be non-zero smooth cut-off functions such that $\big(\support \cutoffDirac\big) \subset [r_-^D,r_+^D]$ and $\big(\support \cutoffWilson\big) \subset [r_-^W,r_+^W]$.
We note that $ar_+^D, ar_+^W < \injectivityradius(X,g)$.

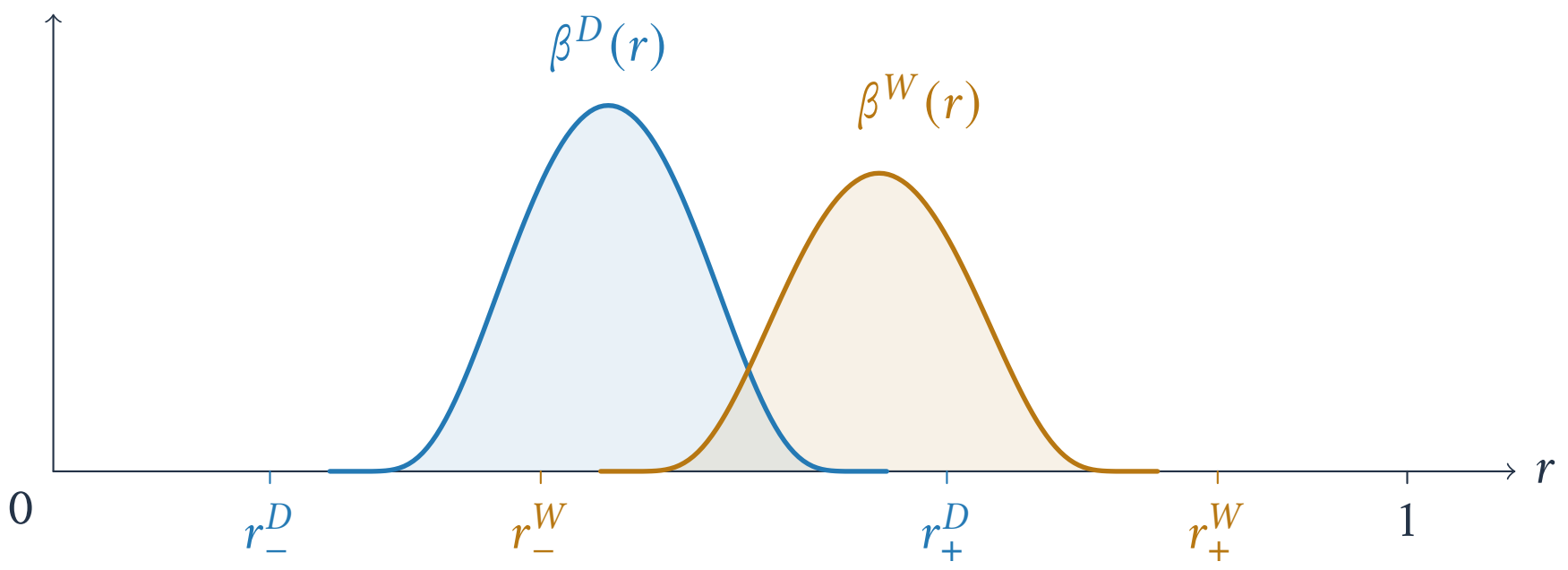
\begin{figure}[htbp]
  \centering
  \begin{tikzpicture}[
    x=10cm,y=2.5cm,
    font=\small,text=ink,
    line cap=round,line join=round
  ]
    \fill[blueedge,opacity=0.10]
      plot[domain=0.01:0.99,samples=121]
      ({0.20+0.42*\x},{1.08*exp(4-1/(\x*(1-\x)))})
      -- cycle;
    \fill[goldedge,opacity=0.10]
      plot[domain=0.01:0.99,samples=121]
      ({0.40+0.42*\x},{0.88*exp(4-1/(\x*(1-\x)))})
      -- cycle;

    \draw[->,thin,ink] (0,0) -- (1.08,0) node[right] {$r$};
    \draw[->,thin,ink] (0,0) -- (0,1.35);
    \node[below left] at (0,0) {$0$};
    \draw[thin,ink] (1,0) -- (1,-0.035) node[below] {$1$};

    \draw[thin,blueedge] (0.16,0) -- (0.16,-0.035)
      node[below,text=blueedge] {$r_-^D$};
    \draw[thin,blueedge] (0.66,0) -- (0.66,-0.035)
      node[below,text=blueedge] {$r_+^D$};
    \draw[thin,goldedge] (0.36,0) -- (0.36,-0.035)
      node[below,text=goldedge] {$r_-^W$};
    \draw[thin,goldedge] (0.86,0) -- (0.86,-0.035)
      node[below,text=goldedge] {$r_+^W$};

    \draw[blueedge,line width=1pt]
      plot[domain=0.01:0.99,samples=121]
      ({0.20+0.42*\x},{1.08*exp(4-1/(\x*(1-\x)))});
    \draw[goldedge,line width=1pt]
      plot[domain=0.01:0.99,samples=121]
      ({0.40+0.42*\x},{0.88*exp(4-1/(\x*(1-\x)))});

    \node[above,text=blueedge] at (0.41,1.12)
      {$\cutoffDirac(r)$};
    \node[above,text=goldedge] at (0.64,0.95)
      {$\cutoffWilson(r)$};
  \end{tikzpicture}
\end{figure}

\subsection{The mesoscopic Dirac operator}\label{subsection: the mesoscopic Dirac operator}

We define the \emph{mesoscopic Dirac operator}
\begin{equation}
  \mesoscopicDiracoperator \colon C^{\infty}(X,E) \to C^{\infty}(X,E)
\end{equation}
by
\begin{equation}
  (\mesoscopicDiracoperator \Phi)(x) := \int_{TX_x} \cliffordmultiplication_x\left(\frac{v}{\abs{v}}\right) \frac{\paralleltransport{\exp_x(v)}{x}(\Phi(\exp_x(v))) - \Phi(x)}{\abs{v}} \frac{1}{C_D} \cutoffDirac\left(\frac{\abs{v}}{a}\right)\,\frac{dv}{a^n},
\end{equation}
where
\begin{equation}
  C_D := \frac{1}{n} \int_{\R^n} \cutoffDirac(\abs{u}) \,du
\end{equation}
is a normalisation constant.
Here and below, $\R^n$ is equipped with its standard Euclidean structure and Lebesgue measure and is used only as a model space for radial integrals depending on $n$ and the cut-off functions.
Such an integral has the same value over any $n$-dimensional Euclidean vector space because a linear isometry preserves the Lebesgue measure.

Since $v \mapsto \cliffordmultiplication_x(v)$ is odd, we have
\begin{equation}
  \int_{TX_x} \cliffordmultiplication_x\big(v/\abs{v}\big) \frac{- \Phi(x)}{\abs{v}} \frac{1}{C_D} \cutoffDirac\big(\abs{v}/a\big)\,\frac{dv}{a^n} =  \frac{-1}{C_D a^n} \left( \int_{TX_x} \frac{\cutoffDirac\big(\abs{v}/a\big)}{\abs{v}^2} \cliffordmultiplication_x(v) \,dv \right) \Phi(x) = 0.
\end{equation}
Thus, we have
\begin{equation}
  (\mesoscopicDiracoperator \Phi)(x) = \int_{TX_x} \cliffordmultiplication_x\big(v/\abs{v}\big) \frac{\paralleltransport{\exp_x(v)}{x}(\Phi(\exp_x(v)))}{\abs{v}} \frac{1}{C_D} \cutoffDirac\big(\abs{v}/a\big)\,\frac{dv}{a^n}.
\end{equation}
Since the integrand is supported in the ball on which $\exp_x$ is injective, we have
\begin{equation}
  \begin{aligned}
& \int_{TX_x} \cliffordmultiplication_x\big(v/\abs{v}\big) \frac{\paralleltransport{\exp_x(v)}{x}\big(\Phi(\exp_x(v))\big)}{\abs{v}} \frac{1}{C_D} \cutoffDirac\big(\abs{v}/a\big)\,\frac{dv}{a^n} \\
  ={}& \int_X \cliffordmultiplication_x\big(v_{xy}/d_g(x,y)\big) \frac{\paralleltransport{\exp_x(v_{xy})}{x}\big(\Phi(\exp_x(v_{xy}))\big)}{d_g(x,y)} \frac{1}{C_D} \cutoffDirac\big(d_g(x,y)/a\big) \,\frac{\dvol_g(y)}{a^nJ_x(v_{xy})} \\
  ={}& \int_X \cliffordmultiplication_x\big(v_{xy}/d_g(x,y)\big) \frac{\paralleltransport{y}{x}(\Phi(y))}{d_g(x,y)} \frac{1}{C_D} \cutoffDirac\big(d_g(x,y)/a\big) \,\frac{\dvol_g(y)}{a^nJ_x(v_{xy})}.
  \end{aligned}
\end{equation}
Combining these identities, we obtain the following integral representation of $\mesoscopicDiracoperator$:
\begin{equation}
  (\mesoscopicDiracoperator \Phi)(x) = \int_X \kernelDirac(x,y) \Phi(y) \,\dvol_g(y),
\end{equation}
where
\begin{equation}
  \kernelDirac(x,y) := \frac{1}{C_D} \cutoffDirac\left(\frac{d_g(x,y)}{a}\right) \frac{1}{d_g(x,y)^2} \frac{1}{a^nJ_x(v_{xy})} \big( \cliffordmultiplication_x(v_{xy}) \circ \paralleltransport{y}{x} \big) \in \Hom(E_y,E_x).
\end{equation}
We set $\kernelDirac(x,y) := 0$ if $d_g(x,y)/a \notin \support(\cutoffDirac)$.

We state some immediate properties of $\kernelDirac(x,y)$ and $\mesoscopicDiracoperator$.

\begin{lemma}\label{lemma: smoothness of Dirac kernel}
  We have
  \begin{equation}
    \kernelDirac \in C^{\infty}(X \times X, \Hom(\pi_2^* E, \pi_1^* E)),
  \end{equation}
  where $\pi_1, \pi_2 \colon X \times X \to X$ are the first and second projections respectively.
\end{lemma}
\begin{proof}
  If $\kernelDirac(x,y) \ne 0$ for $x,y \in X$, we have
  \begin{equation}
    r_-^D a \le d_g(x,y) \le r_+^D a.
  \end{equation}
  Recall that $r_+^D a < \injectivityradius(X,g)$.
  Hence, on a neighbourhood of the support of $\kernelDirac$, the functions $d_g(x,y)$, $d_g(x,y)^{-1}$, $v_{xy} = \exp_x^{-1}(y)$, $\paralleltransport{y}{x}$, and $J_x(v_{xy})^{-1}$ depend smoothly on $(x,y)$.
  Thus, the result follows.
\end{proof}

\begin{corollary}\label{corollary: Dirac kernel is L2}
  $\kernelDirac \in L^2(X \times X, \Hom(\pi_2^* E, \pi_1^* E))$.
\end{corollary}

\begin{lemma}\label{lemma: symmetry of Dirac kernel}
  $\kernelDirac(y,x)^* = \kernelDirac(x,y)$.
\end{lemma}
\begin{proof}
  Since $\nabla$ is a Clifford connection, we have
  \begin{equation}
    \paralleltransport{y}{x} \circ \cliffordmultiplication_y(v_{yx}) = \cliffordmultiplication_x\big(\paralleltransportLC{y}{x}(v_{yx})\big) \circ \paralleltransport{y}{x} = - \cliffordmultiplication_x(v_{xy}) \circ \paralleltransport{y}{x}
  \end{equation}
  by \eqref{equation: reversing the geodesic}, and thus
  \begin{equation}
    \big( \cliffordmultiplication_y(v_{yx}) \circ \paralleltransport{x}{y} \big)^* = \paralleltransport{y}{x} \circ \big( -\cliffordmultiplication_y(v_{yx}) \big) = \cliffordmultiplication_x(v_{xy}) \circ \paralleltransport{y}{x}.
  \end{equation}
  By Lemma~\ref{lemma: symmetry of Jacobian density}, we have $J_x(v_{xy}) = J_y(v_{yx})$.
  Thus, we have
  \begin{equation}
    \begin{aligned}
    \kernelDirac(y,x)^* ={}& \left( \frac{1}{C_D} \cutoffDirac\big(d_g(y,x)/a\big) \frac{1}{d_g(y,x)^2} \frac{1}{J_y(v_{yx})a^n} \big( \cliffordmultiplication_y(v_{yx}) \circ \paralleltransport{x}{y} \big) \right)^* \\
    ={}& \frac{1}{C_D} \cutoffDirac\big(d_g(x,y)/a\big) \frac{1}{d_g(x,y)^2} \frac{1}{a^nJ_x(v_{xy})} \big( \cliffordmultiplication_x(v_{xy}) \circ \paralleltransport{y}{x} \big) = \kernelDirac(x,y),
    \end{aligned}
  \end{equation}
  which completes the proof.
\end{proof}

\begin{proposition}\label{proposition: mesoscopic Dirac is hilbert--schmidt and selfadjoint}
  The mesoscopic Dirac operator $\mesoscopicDiracoperator$ is a Hilbert--Schmidt operator on $L^2(X,E)$, and its operator norm satisfies $\operatornorm{\mesoscopicDiracoperator} = O(1/a)$ as $a \to 0$.
  Moreover, it is self-adjoint and anti-commutes with $\gamma$.
\end{proposition}
\begin{proof}
  Since $\kernelDirac \in L^2(X \times X, \Hom(\pi_2^* E, \pi_1^* E))$, it follows that $\mesoscopicDiracoperator$ is a Hilbert--Schmidt operator.
  The operator $\mesoscopicDiracoperator$ is self-adjoint because $\kernelDirac(y,x) = \kernelDirac(x,y)^*$, and it anti-commutes with $\gamma$ because the parallel transports commute with $\gamma$ and the Clifford multiplications anti-commute with $\gamma$.

  We will show that $\operatornorm{\mesoscopicDiracoperator} = O(1/a)$.
  Since $\norm{\cliffordmultiplication_x(v_{xy})} = \abs{v_{xy}} = d_g(x,y)$ and the parallel transport is unitary, we have
  \begin{equation}
    \norm*{\kernelDirac(x,y)} = \frac{1}{C_D} \cutoffDirac\big(d_g(x,y)/a\big) \frac{1}{d_g(x,y)} \frac{1}{a^nJ_x(v_{xy})}.
  \end{equation}
  Hence, fixing $x$ and changing the variable $y = \exp_x(v)$, we have
  \begin{equation}
    \begin{aligned}
    \int_X \norm*{\kernelDirac(x,y)} \,\dvol_g(y) ={}& \int_X \frac{1}{C_D} \cutoffDirac\big(d_g(x,y)/a\big) \frac{1}{d_g(x,y)} \frac{1}{a^nJ_x(v_{xy})} \,\dvol_g(y) \\
    ={}& \int_{TX_x} \frac{1}{C_D} \cutoffDirac\big(\abs{v}/a\big) \frac{1}{\abs{v}} \frac{1}{a^nJ_x(v)} J_x(v) \,dv \\
    ={}& \frac{1}{C_D} \int_{TX_x} \frac{\cutoffDirac\big(\abs{v}/a\big)}{\abs{v}} \,\frac{dv}{a^n} \\
    ={}& \frac{1}{C_D} \int_{\R^n} \frac{\cutoffDirac(\abs{u})}{a\abs{u}} \,du = \frac{M_D}{a},
    \end{aligned}
  \end{equation}
  where
  \begin{equation}
    M_D := \frac{1}{C_D} \int_{\R^n} \frac{\cutoffDirac(\abs{u})}{\abs{u}} \,du = \frac{n \int_0^{\infty} \cutoffDirac(r) r^{n-2} \,dr}{\int_0^{\infty} \cutoffDirac(r) r^{n-1} \,dr} < \infty.
  \end{equation}
  Since $\kernelDirac(y,x) = \kernelDirac(x,y)^*$, we have
  \begin{equation}
    \int_X \norm*{\kernelDirac(y,x)} \,\dvol_g(y) = \int_X \norm*{\kernelDirac(x,y)} \,\dvol_g(y) = \frac{M_D}{a}.
  \end{equation}
  Hence, we have, for any $x \in X$,
  \begin{equation}
    \begin{aligned}
    \abs{\mesoscopicDiracoperator \Phi(x)}^2 \le{}& \left( \int_X \norm{\kernelDirac(x,y)} \,\dvol_g(y) \right)\left( \int_X \norm{\kernelDirac(x,y)} \abs{\Phi(y)}^2 \,\dvol_g(y)\right) \\
    ={}& \frac{M_D}{a} \int_X \norm{\kernelDirac(x,y)} \abs{\Phi(y)}^2 \,\dvol_g(y)
    \end{aligned}
  \end{equation}
  for $\Phi \in C^{\infty}(X,E)$.
  Hence, we have
  \begin{equation}
    \begin{aligned}
    \Ltwonorm{\mesoscopicDiracoperator \Phi}^2 \le{}& \frac{M_D}{a} \int_X \left( \int_X \norm{\kernelDirac(x,y)} \abs{\Phi(y)}^2 \,\dvol_g(y) \right) \,\dvol_g(x) \\
    ={}& \frac{M_D}{a} \int_X \left( \int_X \norm{\kernelDirac(x,y)} \abs{\Phi(y)}^2 \,\dvol_g(x) \right) \,\dvol_g(y) \\
    ={}& \frac{M_D}{a} \int_X \left( \int_X \norm{\kernelDirac(x,y)} \,\dvol_g(x) \right) \abs{\Phi(y)}^2 \,\dvol_g(y) = \left(\frac{M_D}{a}\right)^2 \Ltwonorm{\Phi}^2.
    \end{aligned}
  \end{equation}
  Hence, we have $\Ltwonorm{\mesoscopicDiracoperator \Phi} \le (M_D/a) \Ltwonorm{\Phi}$ for $\Phi \in C^{\infty}(X,E)$.
  Thus, we have $\operatornorm{\mesoscopicDiracoperator} \le M_D/a$.
\end{proof}

\begin{proposition}\label{proposition: convergence of mesoscopic Dirac operators}
  We have
  \begin{equation}
    \lim_{a \to 0} \Ltwonorm{\mesoscopicDiracoperator \Phi - D\Phi} = 0
  \end{equation}
  for any $\Phi \in C^{\infty}(X,E)$.
\end{proposition}
\begin{proof}
  Fix $\Phi \in C^{\infty}(X,E)$.
  Fix $x \in X$.
  By the Taylor expansion along radial geodesics, we have
  \begin{equation}
    \paralleltransport{\exp_x(v)}{x} \bigl(\Phi(\exp_x(v))\bigr) = \Phi(x) + (\nabla_v \Phi)(x) + R_{\Phi}(x,v)
  \end{equation}
  with
  \begin{equation}
    \abs{R_{\Phi}(x,v)} \le \frac{1}{2}\norm{\nabla^2 \Phi}_{L^{\infty}} \abs{v}^2
  \end{equation}
  for $v \in TX_x$ with $\abs{v} < \injectivityradius(X,g)$.
  By substituting $v = au$ in the definition of $\mesoscopicDiracoperator$, we have
  \begin{equation}
    (\mesoscopicDiracoperator\Phi)(x) = \frac{1}{C_D} \int_{TX_x} \frac{\cliffordmultiplication_x(u) \bigl((\nabla_u\Phi)(x)\bigr)}{\abs{u}^2} \cutoffDirac(\abs{u}) \,du + E_a(x),
  \end{equation}
  where
  \begin{equation}
    \abs{E_a(x)} \le \left( \frac{\norm{\nabla^2 \Phi}_{L^{\infty}}}{2C_D} \int_{\R^n} \abs{u}\cutoffDirac(\abs{u}) \,du \right) a.
  \end{equation}
  For an orthonormal basis $\{e_1,\dots,e_n\}$ of $TX_x$, write $u=\sum_i u_i e_i$.
  By radial symmetry and the definition of $C_D$, we have
  \begin{equation}
    \int_{TX_x} \frac{u_i u_j}{\abs{u}^2} \cutoffDirac(\abs{u}) \,du = \frac{\delta_{ij}}{n} \int_{\R^n} \cutoffDirac(\abs{u}) \,du = C_D \delta_{ij}.
  \end{equation}
  Hence, we have
  \begin{multline}
    \frac{1}{C_D} \int_{TX_x} \frac{\cliffordmultiplication_x(u) \bigl((\nabla_u\Phi)(x)\bigr)}{\abs{u}^2} \cutoffDirac(\abs{u}) \,du = \frac{1}{C_D} \sum_{i,j} \cliffordmultiplication_x(e_i)(\nabla_{e_j}\Phi)(x) \int_{TX_x} \frac{u_i u_j}{\abs{u}^2} \cutoffDirac(\abs{u}) \,du \\
    = \sum_{i,j}\cliffordmultiplication_x(e_i)(\nabla_{e_j}\Phi)(x)\delta_{ij} = \sum_i \cliffordmultiplication_x(e_i)(\nabla_{e_i}\Phi)(x) = D\Phi(x).
  \end{multline}
  Thus, we have
  \begin{equation}
    \Ltwonorm{\mesoscopicDiracoperator\Phi - D\Phi} \le \left(\frac{\norm{\nabla^2 \Phi}_{L^{\infty}} \volume(X,g)^{1/2}}{2C_D} \int_{\R^n} \abs{u}\cutoffDirac(\abs{u}) \,du\right) a \longrightarrow 0
  \end{equation}
  as $a \to 0$.
\end{proof}

Since the mesoscopic Dirac operator $\mesoscopicDiracoperator$ is compact, it is not a Fredholm operator.
This motivates the introduction of the mesoscopic Wilson term.

\subsection{The mesoscopic Wilson term}\label{subsection: the mesoscopic Wilson term}

Fix a normalisation constant $C_W > 0$.
The choice of $C_W$ will not be important in this paper.
One may simply take $C_W=1$, but we keep the normalisation arbitrary.

We first define the integral kernel $\kernelWilson \in C^{\infty}(X \times X, \R)$ by
\begin{equation}
  \kernelWilson(x,y) := \frac{1}{C_W} \cutoffWilson(d_g(x,y)/a) \frac{1}{a^nJ_x(v_{xy})}
\end{equation}
for $x,y \in X$ with $d_g(x,y)\le ar_+^W$, where $v_{xy} = \exp_x^{-1}(y)$ as before.
We set $\kernelWilson(x,y) := 0$ if $d_g(x,y)> ar_+^W$.
By Lemma~\ref{lemma: symmetry of Jacobian density}, we have $\kernelWilson(x,y) = \kernelWilson(y,x)$, and $\kernelWilson(x,y) \ge 0$.

We define the \emph{mesoscopic Wilson term}
\begin{equation}
  \mesoscopicWilsonterm \colon C^{\infty}(X,E) \to C^{\infty}(X,E)
\end{equation}
by
\begin{equation}
  \begin{aligned}
  (\mesoscopicWilsonterm \Phi)(x) :={}& \int_X \frac{\Phi(x) - \paralleltransport{y}{x}(\Phi(y))}{d_g(x,y)} \kernelWilson(x,y) \,\dvol_g(y) \\
  ={}& \int_{TX_x} \frac{\Phi(x) - \paralleltransport{\exp_x(v)}{x}(\Phi(\exp_x(v)))}{\abs{v}} \frac{1}{C_W} \cutoffWilson(\abs{v}/a) \,\frac{dv}{a^n}.
  \end{aligned}
\end{equation}

Setting
\begin{equation}
  \begin{aligned}
  A_a :={}& \int_X \frac{1}{d_g(x,y)} \kernelWilson(x,y) \,\dvol_g(y) \\
  ={}& \int_{TX_x} \frac{1}{\abs{v}} \frac{1}{C_W} \cutoffWilson(\abs{v}/a) \,\frac{dv}{a^n} = \left( \int_{\R^n} \frac{1}{\abs{u}} \frac{1}{C_W} \cutoffWilson(\abs{u}) \,du \right) \frac{1}{a} =: \frac{m_W}{a}
  \end{aligned}
\end{equation}
and
\begin{equation}
  (\mesoscopicWilsontermcompactpert \Phi)(x) := \int_X \frac{-\paralleltransport{y}{x}(\Phi(y))}{d_g(x,y)} \kernelWilson(x,y) \,\dvol_g(y) = \int_X \left(-\frac{\cutoffWilson(d_g(x,y)/a)}{C_W a^n J_x(v_{xy}) d_g(x,y)} \paralleltransport{y}{x} \right) \Phi(y) \,\dvol_g(y),
\end{equation}
we have $\mesoscopicWilsonterm = A_a \id + \mesoscopicWilsontermcompactpert$.
Now $A_a = m_W/a > 0$.
Hence, $A_a \id$ is a bounded operator and $\operatornorm{A_a \id} = m_W/a$.
As was shown for the mesoscopic Dirac operator, $\mesoscopicWilsontermcompactpert$ is a Hilbert--Schmidt operator and $\operatornorm{\mesoscopicWilsontermcompactpert} = O(1/a)$.
Thus, the mesoscopic Wilson term $\mesoscopicWilsonterm$ is a bounded operator and $\operatornorm{\mesoscopicWilsonterm} = O(1/a)$.
Moreover, since $\mesoscopicWilsonterm = A_a \id + \mesoscopicWilsontermcompactpert$ is a compact perturbation by $\mesoscopicWilsontermcompactpert$ of the invertible bounded operator $A_a \id$, it is a Fredholm operator.

The mesoscopic Wilson term $\mesoscopicWilsonterm$ commutes with $\gamma$ because the parallel transports commute with $\gamma$.

Let $\Phi, \Psi \in C^{\infty}(X,E)$.
Since the parallel transports are unitary, we have
\begin{equation}
  \begin{aligned}
  & \pointinnerproduct*{\frac{\Phi(x) - \paralleltransport{y}{x}(\Phi(y))}{d_g(x,y)}}{\frac{\Psi(x) - \paralleltransport{y}{x}(\Psi(y))}{d_g(x,y)}} \\
  ={}& \pointinnerproduct*{\frac{\Phi(x) - \paralleltransport{y}{x}(\Phi(y))}{d_g(x,y)}}{\frac{\Psi(x)}{d_g(x,y)}} + \pointinnerproduct*{\frac{\Phi(x)-\paralleltransport{y}{x}(\Phi(y))}{d_g(x,y)}}{\frac{-\paralleltransport{y}{x}(\Psi(y))}{d_g(x,y)}} \\
  ={}& \pointinnerproduct*{\frac{\Phi(x) - \paralleltransport{y}{x}(\Phi(y))}{d_g(x,y)}}{\frac{\Psi(x)}{d_g(x,y)}} + \pointinnerproduct*{\frac{\paralleltransport{x}{y}\big(\Phi(x)-\paralleltransport{y}{x}(\Phi(y))\big)}{d_g(x,y)}}{\frac{-\Psi(y)}{d_g(x,y)}} \\
  ={}& \pointinnerproduct*{\frac{\Phi(x) - \paralleltransport{y}{x}(\Phi(y))}{d_g(x,y)}}{\frac{\Psi(x)}{d_g(x,y)}} + \pointinnerproduct*{\frac{\Phi(y)-\paralleltransport{x}{y}(\Phi(x))}{d_g(x,y)}}{\frac{\Psi(y)}{d_g(x,y)}}.
  \end{aligned}
\end{equation}
Hence, since $\kernelWilson(x,y) = \kernelWilson(y,x)$, we have
\begin{equation}
  \begin{aligned}
    \innerproduct*{\mesoscopicWilsonterm \Phi}{\Psi}_{L^2} ={}& \int_X \pointinnerproduct*{\int_X \kernelWilson(x,y) \frac{\Phi(x) - \paralleltransport{y}{x}(\Phi(y))}{d_g(x,y)} \,\dvol_g(y)}{\Psi(x)} \,\dvol_g(x) \\
    ={}& \int_X \int_X \kernelWilson(x,y) \pointinnerproduct*{\frac{\Phi(x) - \paralleltransport{y}{x}(\Phi(y))}{d_g(x,y)}}{\Psi(x)} \,\dvol_g(y) \,\dvol_g(x) \\
    ={}& \int_X \int_X \kernelWilson(x,y) d_g(x,y) \pointinnerproduct*{\frac{\Phi(x) - \paralleltransport{y}{x}(\Phi(y))}{d_g(x,y)}}{\frac{\Psi(x)}{d_g(x,y)}} \,\dvol_g(y) \,\dvol_g(x) \\
    ={}& \frac{1}{2} \int_X \int_X \kernelWilson(x,y) d_g(x,y) \pointinnerproduct*{\frac{\Phi(x) - \paralleltransport{y}{x}(\Phi(y))}{d_g(x,y)}}{\frac{\Psi(x) - \paralleltransport{y}{x}(\Psi(y))}{d_g(x,y)}} \,\dvol_g(y) \,\dvol_g(x).
  \end{aligned}
\end{equation}
From this representation, we immediately have $\innerproduct*{\mesoscopicWilsonterm \Phi}{\Psi}_{L^2} = \innerproduct*{\Phi}{\mesoscopicWilsonterm  \Psi}_{L^2}$.
Thus, $\mesoscopicWilsonterm$ is self-adjoint.
Since $\kernelWilson(x,y) \ge 0$, it also follows from this representation that $\mesoscopicWilsonterm \ge 0$.

In summary, we have proved the following proposition.
\begin{proposition}\label{proposition: mesoscopic Wilson term is Fredholm and self-adjoint}
  The mesoscopic Wilson term $\mesoscopicWilsonterm$ is a compact perturbation of an invertible operator on $L^2(X,E)$.
  In particular, it is a bounded Fredholm operator.
  Its operator norm satisfies $\operatornorm{\mesoscopicWilsonterm} = O(1/a)$ as $a \to 0$.
  Moreover, it is self-adjoint and commutes with $\gamma$.
\end{proposition}

\begin{proposition}\label{proposition: convergence of mesoscopic Wilson terms}
  We have
  \begin{equation}
    \lim_{a \to 0} \Ltwonorm{\mesoscopicWilsonterm\Phi} = 0
  \end{equation}
  for any $\Phi \in C^{\infty}(X,E)$.
  This convergence does not hold for every $\Phi \in L^2(X,E)$.
\end{proposition}
\begin{proof}
  Fix $\Phi \in C^{\infty}(X,E)$.
  By changing variables $y=\exp_x(au)$ and using the Taylor expansion as in the proof of Proposition~\ref{proposition: convergence of mesoscopic Dirac operators}, we have
  \begin{equation}
    (\mesoscopicWilsonterm \Phi)(x) = -\frac{1}{C_W} \int_{TX_x} \frac{(\nabla_u \Phi)(x)}{\abs{u}} \cutoffWilson(\abs{u}) \,du + O(a),
  \end{equation}
  uniformly in $x$.
  Since the integrand is odd in $u$, we have
  \begin{equation}
    -\frac{1}{C_W} \int_{TX_x} \frac{(\nabla_u \Phi)(x)}{\abs{u}} \cutoffWilson(\abs{u}) \,du = 0.
  \end{equation}
  Hence, $\mesoscopicWilsonterm \Phi(x) = O(a)$ uniformly in $x$.
  Thus, $\mesoscopicWilsonterm \Phi \to 0$ in $L^2(X,E)$ as $a \to 0$.

  Finally, recall that $\mesoscopicWilsonterm = (m_W/a)\id + \mesoscopicWilsontermcompactpert$, where $m_W > 0$ and $\mesoscopicWilsontermcompactpert$ is compact.
  If $\{e_{\ell}\}_{\ell=1}^{\infty}$ is an orthonormal basis in $L^2(X,E)$, then $\mesoscopicWilsontermcompactpert e_{\ell} \to 0$ for each fixed $a$.
  Thus, we have
  \begin{equation}
    \operatornorm{\mesoscopicWilsonterm} \ge \lim_{\ell \to \infty} \Ltwonorm{\mesoscopicWilsonterm e_\ell} = \frac{m_W}{a}.
  \end{equation}
  Choose a sequence $a_j \to 0$.
  If $\lim_{a \to 0} \Ltwonorm{\mesoscopicWilsonterm\Phi} = 0$ holds for every $\Phi \in L^2(X,E)$, the uniform boundedness principle implies
  \begin{equation}
    \sup_j \operatornorm{\mesoscopicWilsonterm[a_j]} < \infty.
  \end{equation}
  This contradicts the above lower bound.
  This completes the proof.
\end{proof}

\begin{remark}
  We have $\lim_{a \to 0} \Ltwonorm{\mesoscopicWilsonterm\Phi} = 0$ for $\Phi \in L^2_1(X,E)$.
  Indeed, $\Ltwonorm{\mesoscopicWilsonterm\Phi} \le C \norm{\Phi}_{L^2_1}$ with $C$ independent of $a$.
  We omit the details because this is not used below.
\end{remark}

\subsection{The mesoscopic Wilson--Dirac operator}

With these preparations, we are now ready to define the mesoscopic Wilson--Dirac operator.
\begin{definition}\label{definition: mesoscopic Wilson--Dirac operator}
  We define the \emph{mesoscopic Wilson--Dirac operator}
  \begin{equation}
    \mesoscopicWilsonDiracoperator \colon C^{\infty}(X,E) \to C^{\infty}(X,E)
  \end{equation}
  by
  \begin{equation}
    \mesoscopicWilsonDiracoperator := \mesoscopicDiracoperator + \mesoscopicWilsonterm \circ \gamma
  \end{equation}
  for $a \in \big(0,\injectivityradius(X,g)/10\big)$.
\end{definition}

Since $\mesoscopicDiracoperator$ is compact and $\mesoscopicWilsonterm$ is a Fredholm operator, $\mesoscopicWilsonDiracoperator$ is also a Fredholm operator.
Note that $\mesoscopicWilsonDiracoperator$ neither commutes nor anti-commutes with $\gamma$.
The following proposition is immediate from the discussion above and Weyl's theorem on the stability of the essential spectrum under compact perturbations.
Note that $\mesoscopicWilsonDiracoperator - (m_W/a)\gamma = \mesoscopicDiracoperator + \mesoscopicWilsontermcompactpert \gamma$ is compact.
\begin{proposition}\label{proposition: mesoscopic WD is bounded self-adjoint}
  The mesoscopic Wilson--Dirac operator $\mesoscopicWilsonDiracoperator$ is a bounded self-adjoint operator on $L^2(X,E)$, and its operator norm satisfies $\operatornorm{\mesoscopicWilsonDiracoperator} = O(1/a)$ as $a \to 0$.
  Moreover, it is a Fredholm operator, and its essential spectrum satisfies $\essentialspectrum(\mesoscopicWilsonDiracoperator) = \{ \pm m_W/a \}$.
\end{proposition}

By Propositions~\ref{proposition: convergence of mesoscopic Dirac operators} and~\ref{proposition: convergence of mesoscopic Wilson terms}, we have
\begin{equation}
  \lim_{a \to 0} \Ltwonorm{\mesoscopicWilsonDiracoperator \Phi - D \Phi} = 0
\end{equation}
for any $\Phi \in C^{\infty}(X,E)$.
Thus, in this sense, the mesoscopic Wilson--Dirac operator $\mesoscopicWilsonDiracoperator$ approximates the continuum Dirac operator $D$ as $a \to 0$.
In the following sections, we investigate $\mesoscopicWilsonDiracoperator$ in greater detail.
We retain the cut-off functions and normalisation constants fixed in this section.

\subsection{A no-go theorem of Nielsen--Ninomiya type}

The mesoscopic Dirac operator $\mesoscopicDiracoperator$ anti-commutes with $\gamma$, but it is compact and hence not Fredholm.
The mesoscopic Wilson--Dirac operator $\mesoscopicWilsonDiracoperator$ is Fredholm, but it does not anti-commute with $\gamma$.
Both operators are sums of a multiplication operator and a compact operator.
In this subsection, we prove that, if a self-adjoint Fredholm operator of this form anti-commutes with $\gamma$, then the trace of $\gamma$ on its kernel is zero.
The Nielsen--Ninomiya theorem~\cites{MR0621425, MR0639112, NNnogo} in lattice gauge theory says that, under natural assumptions such as locality and translation invariance, a lattice Dirac operator that anti-commutes with $\gamma$ has additional zero modes known as fermion doublers.
Theorem~\ref{theorem: no-go theorem} is a version of this theorem at the mesoscopic scale without translation invariance.

\begin{lemma}\label{lemma: multiplication compact invertible}
  Let $c \in C^0(X, \End(E))$.
  If the multiplication operator $c$ on $L^2(X,E)$ is Fredholm, then $c(x)$ is invertible for every $x \in X$.
  If $c$ is compact, then $c = 0$.
\end{lemma}
\begin{proof}
  We first construct a sequence of sections concentrating at a point.
  Fix $x_0 \in X$ and $v_0 \in E_{x_0}$ with $\norm{v_0} = 1$.
  Choose $j_0$ such that $1/j_0 < \injectivityradius(X,g)$.
  For $j \ge j_0$, let $B_j := B(x_0,1/j)$.
  We define $\Phi_j \in L^2(X,E)$ by
  \begin{equation}
    \Phi_j(x) := \frac{1}{\sqrt{\abs{B_j}}} \paralleltransport{x_0}{x}(v_0)
  \end{equation}
  for $x \in B_j$ and $0$ otherwise.
  We have $\Ltwonorm{\Phi_j} = 1$.
  By the Cauchy--Schwarz inequality, we have
  \begin{equation}
    \abs*{\innerproduct{\Phi_j}{\Theta}_{L^2}} \le \Ltwonorm[L^2(B_j)]{\Theta}
  \end{equation}
  for any $\Theta \in L^2(X,E)$.
  The right-hand side goes to zero as $j \to \infty$.
  Thus, $\Phi_j \rightharpoonup 0$ weakly in $L^2(X,E)$.
  Since $c$ is continuous, we have
  \begin{equation}\label{equation: cphij norm}
    \lim_{j \to \infty} \Ltwonorm{c\Phi_j}^2 = \lim_{j \to \infty} \frac{1}{\abs{B_j}} \int_{B_j} \norm{c(x) \paralleltransport{x_0}{x}(v_0)}^2 \,\dvol_g = \norm{c(x_0)v_0}^2.
  \end{equation}

  Suppose that $c$ is Fredholm.
  Let $\Pi$ be the orthogonal projection onto $\kernel c$.
  Since the range of $c$ is closed, there exists $C > 0$ such that
  \begin{equation}
    \Ltwonorm{c\Phi} \ge C \Ltwonorm{\Phi - \Pi \Phi}
  \end{equation}
  for any $\Phi \in L^2(X,E)$.
  Since $\Pi$ has finite rank, we have $\Pi \Phi_j \to 0$ strongly in $L^2(X,E)$.
  Thus, by \eqref{equation: cphij norm}, we have
  \begin{equation}
    \norm{c(x_0)v_0} = \lim_{j \to \infty} \Ltwonorm{c\Phi_j} \ge C \lim_{j \to \infty} \Ltwonorm{\Phi_j - \Pi \Phi_j} = C.
  \end{equation}
  Since $x_0$ and $v_0$ are arbitrary, $c(x)$ is invertible for every $x \in X$.

  Suppose next that $c$ is compact.
  Then, since $\Phi_j \rightharpoonup 0$ weakly, $c\Phi_j \to 0$ strongly in $L^2(X,E)$.
  Hence, by \eqref{equation: cphij norm}, we have $c(x_0)v_0 = 0$.
  Since $x_0$ and $v_0$ are arbitrary, we have $c = 0$.
\end{proof}

\begin{theorem}\label{theorem: no-go theorem}
  Let $c \in C^0(X,\End(E))$, and let $K$ be a compact operator on $L^2(X,E)$.
  Suppose that $T := c+K$ is self-adjoint and anti-commutes with $\gamma$.
  Then, $T$ is a Fredholm operator if and only if $c(x)$ is invertible for every $x \in X$.
  Moreover, if $T$ is Fredholm, then we have $\trace \bigl( \gamma |_{\kernel T} \bigr) = 0$.
\end{theorem}
\begin{proof}
  Since $K$ is compact, $T$ is Fredholm if and only if $c$ is Fredholm.
  By Lemma~\ref{lemma: multiplication compact invertible}, if $c$ is Fredholm, then $c(x)$ is invertible for every $x \in X$.
  Thus, the multiplication operator $c$ is invertible.
  The converse is clear.

  Assume that $T$ is Fredholm.
  Since $T$ anti-commutes with $\gamma$, we have
  \begin{equation}
    c \gamma + \gamma c = - (K \gamma + \gamma K).
  \end{equation}
  The left-hand side is a multiplication operator, and the right-hand side is compact.
  Hence, by Lemma~\ref{lemma: multiplication compact invertible}, we have $c \gamma + \gamma c = 0$.
  Thus, $T$ and $c$ have the forms
  \begin{equation}
    T = \begin{pmatrix}
      0 & T^- \\
      T^+ & 0
    \end{pmatrix}
  \end{equation}
  and
  \begin{equation}
    c = \begin{pmatrix}
      0 & c^- \\
      c^+ & 0
    \end{pmatrix}
  \end{equation}
  with respect to the eigenbundle decomposition of $E$ by $\gamma$.
  We have $(T^+)^* = T^-$.
  By the first assertion, the multiplication operator $c^+$ is invertible.
  Since $T-c = K$ is compact, the operator $T^+ - c^+$ is compact.
  Hence, $T^+$ is a compact perturbation of an invertible operator, and its Fredholm index is $0$.
  Thus, we have
  \begin{equation}
    \trace \bigl( \left. \gamma \right|_{\kernel T} \bigr) = \dim \kernel T^+ - \dim \kernel T^- = 0.
  \end{equation}
  This completes the proof.
\end{proof}

For the mesoscopic Dirac operator $\mesoscopicDiracoperator$, we have $c = 0$.
For the mesoscopic Wilson--Dirac operator $\mesoscopicWilsonDiracoperator$, we have $c = (m_W/a)\gamma$, which is invertible, but $\mesoscopicWilsonDiracoperator$ neither commutes nor anti-commutes with $\gamma$.

\section{Elliptic estimates on the flat torus}\label{section: Elliptic estimates on the flat torus}

In this section, we study the mesoscopic Wilson--Dirac operator in detail on a flat torus and, in particular, establish an elliptic estimate for it.
The analysis in this flat setting will serve as a basis for the corresponding analysis on general closed Riemannian manifolds in the next section.

\subsection{The flat model and Fourier expansion}

Let $V$ be an $n$-dimensional Euclidean vector space, and $\Lattice \subset V$ a lattice of full rank.
Set
\begin{equation}
  \flattorus := V / \Lattice,
\end{equation}
and equip $\flattorus$ with the induced flat Riemannian metric.
We use the same notation for a point in $\flattorus$ and any of its lifts to $V$.
We use the additive notation $x+v$ for the action of $V$ on $\flattorus$.
We have
\begin{equation}
  \injectivityradius(\flattorus) = \min_{\ell \in \Lattice \setminus \{0\}} \frac{\abs{\ell}}{2}.
\end{equation}
We fix the canonical trivialisation $T\flattorus = \flattorus \times V$.

Let $S$ be a graded Hermitian Clifford module over $V$ with a grading operator $\gamma$ and a Clifford multiplication $\cliffordmultiplication \colon V \to \End(S)$.
Fix an orthonormal basis $\{e_1, \dots, e_n\}$ of $V$, and let $\{e^1, \dots, e^n\}$ be the dual basis of $V^*$.
Set $c_j := \cliffordmultiplication(e_j)$ for $j = 1,\dots,n$.
Then, we have
\begin{equation}
  \gamma c_j + c_j \gamma = 0, \quad c_j c_k + c_k c_j = -2 \delta_{jk}, \quad c_j^* = -c_j
\end{equation}
for $j,k = 1,\dots,n$.
Set
\begin{equation}
  \flatmodel{E} := \flattorus \times S,
\end{equation}
and equip $\flatmodel{E}$ with the constant grading operator and Clifford multiplication induced by those on $S$ and with the product connection $\nabla$.
We use the same notation $\gamma$, $\cliffordmultiplication$, and $c_j$ for $\flatmodel{E}$.
Let
\begin{equation}
  \duallattice{\Lattice} := \{ \xi \in V^* \mid \xi(\ell) \in 2\pi\Z \text{ for every } \ell \in \Lattice \}
\end{equation}
be the dual lattice of $\Lattice$.
For $\xi \in \duallattice{\Lattice}$, we define
\begin{equation}
  e_{\xi}(x) := \frac{1}{\sqrt{\volume(\flattorus)}} e^{i\xi(x)}.
\end{equation}
Then, $\{e_{\xi}\}_{\xi \in \duallattice{\Lattice}}$ is an orthonormal basis of $L^2(\flattorus)$.
Every $\Phi \in L^2(\flattorus,\flatmodel{E})$ has a Fourier expansion
\begin{equation}
  \Phi = \sum_{\xi \in \duallattice{\Lattice}} e_{\xi} \Phi_{\xi}
\end{equation}
with $\Phi_{\xi} \in S$.
We have
\begin{equation}
  \Ltwonorm{\Phi}^2 = \sum_{\xi \in \duallattice{\Lattice}} \norm{\Phi_{\xi}}_S^2
\end{equation}
by the Plancherel formula.
We have
\begin{equation}
  D(e_{\xi} \Phi_{\xi}) = \sum_{j=1}^n c_j \partial_j (e_{\xi}\Phi_{\xi}) = \sum_{j=1}^n c_j (i\xi_j e_{\xi} \Phi_{\xi}) = e_{\xi} \cdot i \left(\sum_{j=1}^n \xi_j c_j\right) \Phi_{\xi}
\end{equation}
for $\xi = \sum_{j=1}^n \xi_j e^j \in \duallattice{\Lattice}$ and $\Phi_{\xi} \in S$.

\subsection{Fourier symbols of the mesoscopic operators}

Throughout this section, $a \in (0,\injectivityradius(\flattorus)/10)$ as in Section~\ref{section: mesoscopic WD operators}.
In the flat model $\flattorus$, the exponential map, parallel transport, and the Jacobian density are given by
\begin{equation}
  \exp_x(v) = x + v, \quad \paralleltransport{x+v}{x} = \id_S, \quad J_x(v) = 1
\end{equation}
for $x \in \flattorus$ and $v \in V$.
We denote the mesoscopic Dirac operator, the Wilson term, and the mesoscopic Wilson--Dirac operator associated with this flat model by $\flatmesoscopicDiracoperator$, $\flatmesoscopicWilsonterm$, and $\flatmesoscopicWilsonDiracoperator$.
Let $\Phi \in C^{\infty}(\flattorus,\flatmodel{E})$.
After changing variables $v = aw$ and using the symmetry $w \mapsto -w$, we can rewrite the mesoscopic Dirac operator as
\begin{equation}
  \begin{aligned}
  (\flatmesoscopicDiracoperator \Phi)(x) ={}& \int_{(T\flattorus)_x} \cliffordmultiplication\big(v/\abs{v}\big) \frac{\paralleltransport{\exp_x(v)}{x}(\Phi(\exp_x(v)))}{\abs{v}} \frac{1}{C_D} \cutoffDirac\big(\abs{v}/a\big)\,\frac{dv}{a^n} \\
   ={}& \int_V \cliffordmultiplication\left(\frac{aw}{\abs{aw}}\right) \frac{\Phi(x+aw)}{a\abs{w}} \frac{1}{C_D} \cutoffDirac(\abs{w}) \,dw \\
   ={}& \int_V \cliffordmultiplication\big(w/\abs{w}\big) \frac{\Phi(x+aw) - \Phi(x-aw)}{2a\abs{w}} \frac{1}{C_D} \cutoffDirac(\abs{w}) \,dw.
  \end{aligned}
\end{equation}
Similarly, for the mesoscopic Wilson term, we have
\begin{equation}
  \begin{aligned}
  (\flatmesoscopicWilsonterm \Phi)(x) ={}& \int_{\flattorus} \frac{\Phi(x) - \Phi(y)}{d_g(x,y)} \frac{1}{C_W} \cutoffWilson\big(d_g(x,y)/a\big) \frac{1}{a^n} \,\dvol_g(y) \\
  ={}& \int_V \frac{\Phi(x) - \Phi(x+v)}{\abs{v}} \frac{1}{C_W} \cutoffWilson\big(\abs{v}/a\big) \,\frac{dv}{a^n} \\
  ={}& \int_V \frac{\Phi(x) - \Phi(x+aw)}{a\abs{w}} \frac{1}{C_W} \cutoffWilson(\abs{w}) \,dw \\
  ={}& \int_V \frac{2\Phi(x) - \Phi(x+aw) - \Phi(x-aw)}{2a\abs{w}} \frac{1}{C_W} \cutoffWilson(\abs{w}) \,dw,
  \end{aligned}
\end{equation}
where the passage from the first line to the second uses the fact that the integrand is supported within the injectivity radius so that the map $v \mapsto x+v$ is injective on the relevant region in $V$.

We now compute the Fourier symbols of the mesoscopic operators.
Fix $\xi \in \duallattice{\Lattice}$ and $\phi \in S$.
We first consider the mesoscopic Dirac operator.
\begin{equation}
  \begin{aligned}
  \big(\flatmesoscopicDiracoperator (e_{\xi}\phi)\big)(x) ={}& \left(\int_V \cliffordmultiplication\big(w/\abs{w}\big) \frac{e_{\xi}(x+aw) - e_{\xi}(x-aw)}{2a\abs{w}} \frac{1}{C_D} \cutoffDirac(\abs{w}) \,dw\right) \phi \\
  ={}& \left(\int_V \cliffordmultiplication\big(w/\abs{w}\big) \frac{i \sin(a\xi(w))}{a\abs{w}} \frac{1}{C_D} \cutoffDirac(\abs{w}) \,dw\right) e_{\xi}(x)\phi.
  \end{aligned}
\end{equation}
Thus, defining $\symbolDirac \colon V^* \to \End(S)$ by
\begin{equation}
  \symbolDirac(\eta) := \int_V \cliffordmultiplication\big(w/\abs{w}\big) \frac{i \sin(\eta(w))}{\abs{w}} \frac{1}{C_D} \cutoffDirac(\abs{w}) \,dw,
\end{equation}
we have
\begin{equation}
  \flatmesoscopicDiracoperator (e_{\xi}\phi) = e_{\xi} \frac{\symbolDirac(a\xi)}{a} \phi
\end{equation}
and $\symbolDirac(\eta)$ is self-adjoint and anti-commutes with $\gamma$.
For the mesoscopic Wilson term, we define $\symbolWilson \colon V^* \to \R$ by
\begin{equation}
  \symbolWilson(\eta) := \int_V \frac{1 - \cos(\eta(w))}{\abs{w}} \frac{1}{C_W} \cutoffWilson(\abs{w}) \,dw.
\end{equation}
Then, we have
\begin{equation}
  \flatmesoscopicWilsonterm (e_{\xi}\phi) = e_{\xi} \frac{\symbolWilson(a\xi)}{a} \phi
\end{equation}
and $\symbolWilson(\eta) \in [0,\infty)$.
Thus, setting $\symbolWilsonDirac := \symbolDirac + \symbolWilson \gamma$, we have
\begin{equation}\label{equation: Fourier transform of mesoscopic WD operators}
  \flatmesoscopicWilsonDiracoperator (e_{\xi}\phi) = e_{\xi} \frac{\symbolWilsonDirac(a\xi)}{a} \phi.
\end{equation}
We call $\symbolDirac$, $\symbolWilson$, and $\symbolWilsonDirac$ the \emph{symbols} of the corresponding operators.
Note that the symbols are independent of a point in $\flattorus$ since the Clifford multiplication is constant under the trivialisation $\flatmodel{E}=\flattorus \times S$.

We emphasise that the dependence on the mesoscopic scale $a$ is completely separated in the form
\begin{equation}
  \frac{\symbolWilsonDirac(a\xi)}{a},
\end{equation}
where $\eta:=a\xi$ is the dimensionless frequency obtained by measuring the frequency $\xi$ at the mesoscopic length scale $a$.
The regime $\abs{\eta} \ll 1$ corresponds to frequencies much smaller than $1/a$, whereas $\abs{\eta} \gg 1$ corresponds to frequencies much larger than $1/a$.
Equivalently, the graph of the symbol
\begin{equation}
  \xi \mapsto \frac{\symbolWilsonDirac(a\xi)}{a}
\end{equation}
is obtained from the graph of the dimensionless symbol
\begin{equation}
  \eta \mapsto \symbolWilsonDirac(\eta)
\end{equation}
by a dilation by a factor $1/a$ in both the frequency and energy directions.
Thus, the shape of the symbol is independent of the mesoscopic scale $a$.
Changing $a$ only rescales the frequency and energy variables by the factor $1/a$.
In particular, any estimate for the dimensionless symbol $\symbolWilsonDirac(\eta)$ with constants independent of $\eta$ immediately yields an estimate uniform in $a$ for $\symbolWilsonDirac(a\xi)/a$.
Moreover, the transition scale $\abs{\eta} \sim 1$ corresponds precisely to the mesoscopic frequency scale $\abs{\xi} \sim 1/a$.

\subsection{Radial reduction of the symbols}

The symbols have a natural $O(V)$-symmetry, which reduces their analysis to the radial variable.
We denote by $\eta^{\sharp} \in V$ the metric dual of $\eta \in V^*$, and let $O(V)$ act on $V^*$ by $R \eta := \eta \circ R^{-1}$.
Note that $(R\eta)^{\sharp} = R \eta^{\sharp}$.

We have
\begin{equation}
  \symbolDirac(\eta) = \int_V \cliffordmultiplication\big(w/\abs{w}\big) \frac{i \sin(\eta(w))}{\abs{w}} \frac{1}{C_D} \cutoffDirac(\abs{w}) \,dw = i \cliffordmultiplication \left( \int_V \frac{w}{\abs{w}} \frac{\sin(\eta(w))}{\abs{w}} \frac{1}{C_D} \cutoffDirac(\abs{w}) \,dw \right)
\end{equation}
for any $\eta \in V^*$.
By the change of variables $w = Rw'$, we have
\begin{equation}
  \begin{aligned}
  \symbolDirac(R\eta) ={}& i \cliffordmultiplication \left( \int_V \frac{w}{\abs{w}} \frac{\sin((R\eta)(w))}{\abs{w}} \frac{1}{C_D} \cutoffDirac(\abs{w}) \,dw \right) \\
  ={}& i \cliffordmultiplication \left( \int_V \frac{Rw'}{\abs{Rw'}} \frac{\sin((R\eta)(Rw'))}{\abs{Rw'}}  \frac{1}{C_D} \cutoffDirac(\abs{Rw'}) \abs{\det R} \,dw' \right) \\
  ={}& i \cliffordmultiplication \left( R \left( \int_V \frac{w'}{\abs{w'}} \frac{\sin(\eta(w'))}{\abs{w'}} \frac{1}{C_D} \cutoffDirac(\abs{w'}) \,dw' \right) \right)
  \end{aligned}
\end{equation}
for any $R \in O(V)$.
Now fix any non-zero $\eta \in V^*$.
For any $R \in O(V)$ with $R \eta = \eta$, we have
\begin{equation}
  i \cliffordmultiplication \left( R \left( \int_V \frac{w}{\abs{w}} \frac{\sin(\eta(w))}{\abs{w}} \frac{1}{C_D} \cutoffDirac(\abs{w}) \,dw \right) \right) = i \cliffordmultiplication \left( \int_V \frac{w}{\abs{w}} \frac{\sin(\eta(w))}{\abs{w}} \frac{1}{C_D} \cutoffDirac(\abs{w}) \,dw \right).
\end{equation}
Since the Clifford multiplication $\cliffordmultiplication \colon V \to \End(S)$ is injective, we have
\begin{equation}
  R \left( \int_V \frac{w}{\abs{w}} \frac{\sin(\eta(w))}{\abs{w}} \frac{1}{C_D} \cutoffDirac(\abs{w}) \,dw \right) =\int_V \frac{w}{\abs{w}} \frac{\sin(\eta(w))}{\abs{w}} \frac{1}{C_D} \cutoffDirac(\abs{w}) \,dw.
\end{equation}
Since every vector invariant under the stabiliser of $\eta$ in $O(V)$ is parallel to $\eta^{\sharp}$, we have
\begin{equation}
   \int_V \frac{w}{\abs{w}} \frac{\sin(\eta(w))}{\abs{w}} \frac{1}{C_D} \cutoffDirac(\abs{w}) \,dw \in \R \frac{\eta^{\sharp}}{\abs{\eta}}.
\end{equation}
Thus, the $O(V)$-equivariance above and the transitivity of the $O(V)$-action on spheres in $V^*$ imply that there exists
\begin{equation}
  \radialsymbolDirac \colon [0,\infty) \to \R
\end{equation}
such that
\begin{equation}
  \symbolDirac(\eta) = i \radialsymbolDirac(\abs{\eta}) \cliffordmultiplication \big( \eta^{\sharp}/\abs{\eta} \big)
\end{equation}
for any non-zero $\eta \in V^*$.
Set $\radialsymbolDirac(0) = 0$.
We obtain an explicit expression for $\radialsymbolDirac$ as follows.
Fix $\zeta \in V^*$ with $\abs{\zeta} = 1$ and $s \in [0,\infty)$.
By definition, we have
\begin{equation}
  \int_V \frac{w}{\abs{w}} \frac{\sin\big((s\zeta)(w)\big)}{\abs{w}} \frac{1}{C_D} \cutoffDirac(\abs{w}) \,dw = \radialsymbolDirac(s) \zeta^{\sharp}.
\end{equation}
Thus, we have
\begin{equation}
  \radialsymbolDirac(s) = \int_V \frac{\zeta(w)}
  {\abs{w}} \frac{\sin(s\zeta(w))}{\abs{w}} \frac{1}{C_D} \cutoffDirac(\abs{w}) \,dw
\end{equation}
for any $\zeta \in V^*$ with $\abs{\zeta} = 1$ and $s \in [0,\infty)$.
Similarly, for the mesoscopic Wilson term, since we have
\begin{equation}
  \symbolWilson(R\eta) = \symbolWilson(\eta),
\end{equation}
there exists $\radialsymbolWilson \colon [0,\infty) \to \R$ such that
\begin{equation}
  \symbolWilson(\eta) = \radialsymbolWilson(\abs{\eta})
\end{equation}
for any $\eta \in V^*$.
We have
\begin{equation}
  \radialsymbolWilson(s) = \int_V \frac{1 - \cos(s\zeta(w))}{\abs{w}} \frac{1}{C_W} \cutoffWilson(\abs{w}) \,dw
\end{equation}
for any $\zeta \in V^*$ with $\abs{\zeta} = 1$ and $s \in [0,\infty)$.
In summary, we have
\begin{equation}
  \symbolWilsonDirac(\eta) = i \radialsymbolDirac(\abs{\eta}) \cliffordmultiplication\big(\eta^{\sharp}/\abs{\eta}\big) + \radialsymbolWilson(\abs{\eta}) \gamma
\end{equation}
for any non-zero $\eta \in V^*$.

We next examine the asymptotic behaviour of $\radialsymbolDirac$ and $\radialsymbolWilson$.
Set
\begin{equation}
  \kappa_W := \frac{1}{2nC_W}\int_{\R^n} \abs{u} \cutoffWilson(\abs{u}) \,du.
\end{equation}
Recall $C_D = (1/n)\int_{\R^n} \cutoffDirac(\abs{u}) \,du$ and $m_W = (1/C_W) \int_{\R^n} \cutoffWilson(\abs{u})/\abs{u}\,du$.

\begin{lemma}\label{lemma: asymptotics s to 0}
  We have
  \begin{gather}
    \radialsymbolDirac(s) = s + O(s^3), \\
    \radialsymbolWilson(s) = \kappa_W s^2 + O(s^4)
  \end{gather}
  as $s \to 0$.
\end{lemma}
\begin{proof}
  Fix $\zeta \in V^*$ with $\abs{\zeta}=1$.
  By the Taylor expansion, for $w \in V$, we have
  \begin{equation}
    \sin (s\zeta(w)) = s\zeta(w) + O\big(s^3 \abs{\zeta(w)}^3\big).
  \end{equation}
  The $O(V)$-symmetry yields
  \begin{equation}
    \int_V \frac{\zeta(w)}{\abs{w}} \frac{s\zeta(w)}{\abs{w}} \frac{1}{C_D} \cutoffDirac(\abs{w}) \,dw = s \frac{1}{C_D} \int_V \cutoffDirac(\abs{w}) \left(\frac{\zeta(w)}{\abs{w}}\right)^2 \,dw = s \frac{1}{C_D} \cdot \frac{1}{n} \int_V \cutoffDirac(\abs{w}) \,dw = s.
  \end{equation}
  Thus, since $\cutoffDirac$ is compactly supported, we have $\radialsymbolDirac(s) = s + O(s^3)$.
  Similarly, we have
  \begin{equation}
    1 - \cos(s\zeta(w)) = \frac{s^2 (\zeta(w))^2}{2} + O(s^4 \abs{\zeta(w)}^4),
  \end{equation}
  and radial symmetry gives
  \begin{equation}
    \int_V \frac{\zeta(w)^2}{\abs{w}} \cutoffWilson(\abs{w}) \,dw = \frac{1}{n} \int_V \abs{w} \cutoffWilson(\abs{w}) \,dw.
  \end{equation}
  Hence, we have $\radialsymbolWilson(s) = \kappa_W s^2 + O(s^4)$.
\end{proof}

\begin{lemma}\label{lemma: asymptotics s to infinity}
  We have
  \begin{gather}
    \radialsymbolDirac(s) \to 0, \\
    \radialsymbolWilson(s) \to m_W
  \end{gather}
  as $s \to \infty$.
\end{lemma}
\begin{proof}
  This follows from the Riemann--Lebesgue lemma, but it can also be directly proved as follows.
  Fix $\zeta \in V^*$ with $\abs{\zeta}=1$, and let $\zeta^{\sharp} \in V$ denote its metric dual.
  Since $\support \cutoffDirac \subset [r_-^D,r_+^D]$ and $r_-^D > 0$, we have $\cutoffDirac(\abs{w}) \zeta(w)/\abs{w}^2 \in C^{\infty}_c(V)$.
  Integration by parts gives
  \begin{equation}
    \begin{aligned}
    \radialsymbolDirac(s) ={}& \int_V \zeta(w) \sin(s\zeta(w)) \frac{1}{C_D} \frac{\cutoffDirac(\abs{w})}
  {\abs{w}^2} \,dw \\
  ={}& \int_V \zeta(w) \left( -\frac{1}{s} \partial_{\zeta^{\sharp}} \cos(s\zeta(w)) \right) \frac{1}{C_D} \frac{\cutoffDirac(\abs{w})}{\abs{w}^2} \,dw \\
  ={}& \frac{1}{C_D s} \int_V \left( \partial_{\zeta^{\sharp}} \left(\frac{\cutoffDirac(\abs{w})}{\abs{w}^2} \zeta(w)\right) \right) \cos(s\zeta(w)) \,dw.
    \end{aligned}
  \end{equation}
  Hence, we have
  \begin{equation}
    \abs{\radialsymbolDirac(s)} \le \left(\frac{1}{C_D} \int_V \abs*{\partial_{\zeta^{\sharp}} \left(\frac{\cutoffDirac(\abs{w})}{\abs{w}^2} \zeta(w)\right)} \,dw\right) \frac{1}{s}.
  \end{equation}
  Thus, $\radialsymbolDirac(s) \to 0$ as $s \to \infty$.
  The same argument applies to $\radialsymbolWilson - m_W$.
\end{proof}

\begin{lemma}\label{lemma: radialsymbolWilson positive}
  We have $\radialsymbolWilson(s) > 0$ for any $s \in (0,\infty)$.
\end{lemma}
\begin{proof}
  Fix $s > 0$ and $\zeta \in V^*$ with $\abs{\zeta}=1$.
  Since $\cutoffWilson$ is non-zero and non-negative, the function $w \mapsto \cutoffWilson(\abs{w})$ is positive on some non-empty open subset of $V$.
  The zero set of $w \mapsto 1-\cos(s\zeta(w))$ is a union of parallel affine hyperplanes
  and hence contains no non-empty open subset of $V$.
  Hence, the integrand is positive on some non-empty open subset of $V$ and non-negative everywhere.
  Thus, $\radialsymbolWilson(s) > 0$.
\end{proof}

These lemmata can be summarised as follows.
\begin{itemize}
  \item
    At low frequencies, the Dirac symbol $\symbolDirac$ is linear to leading order, whereas the Wilson symbol $\symbolWilson$ is quadratic to leading order.
    Thus, the Dirac term gives the leading contribution to the mesoscopic Wilson--Dirac operator $\flatmesoscopicWilsonDiracoperator$.
    This is consistent with the convergence $\flatmesoscopicWilsonDiracoperator \Phi \to D\Phi$ as $a \to 0$.
  \item
    At high frequencies, by contrast, the Dirac symbol $\symbolDirac$ tends to zero, whereas the Wilson symbol $\symbolWilson$ tends to the positive constant $m_W$.
    Thus, the high-frequency modes are controlled by the Wilson term.
    This is the mechanism underlying the elliptic estimates for $\flatmesoscopicWilsonDiracoperator$, which will be proved later.
    From the viewpoint of lattice gauge theory, it is also precisely the mechanism by which the doubler modes are removed.
\end{itemize}

\subsection{Examples}

We digress to one-dimensional examples in which the symbols can be computed explicitly.
Let $V := \R$ with its standard Euclidean structure and let $\Lattice := \Z$.
We identify both $V$ and $V^*$ with $\R$.
We set $S := \C^2$ and
\begin{equation}
  \gamma = \begin{pmatrix}1 & 0 \\ 0 & -1\end{pmatrix}, \quad c_1 = \begin{pmatrix}0 & -1 \\ 1 & 0\end{pmatrix}.
\end{equation}
Then, we have
\begin{equation}
  \symbolWilsonDirac(\eta) = \begin{pmatrix}
    \radialsymbolWilson(\abs{\eta}) & -i\radialsymbolDirac\big(\abs{\eta}\big) \big(\eta/\abs{\eta}\big) \\
    i\radialsymbolDirac\big(\abs{\eta}\big) \big(\eta/\abs{\eta}\big) & -\radialsymbolWilson(\abs{\eta})
  \end{pmatrix}
\end{equation}
for any non-zero $\eta \in V^*$.
Thus, the eigenvalues of $\symbolWilsonDirac(\eta)$ are
\begin{equation}
  \pm \sqrt{\radialsymbolDirac(\abs{\eta})^2 + \radialsymbolWilson(\abs{\eta})^2}
\end{equation}
for any $\eta \in V^*$.

\subsubsection{A sharp annular kernel}

Fix $0 < r_- < r_+ < 1$, and take
\begin{equation}
  \cutoffDirac(r) = \cutoffWilson(r) =
  \begin{cases}
    r & \text{for } r_- \le r \le r_+ \\
    0 & \text{otherwise}
  \end{cases}
\end{equation}
for $r \in [0,\infty)$.
These cut-off functions are not smooth and hence do not satisfy the standing smoothness assumptions, but this causes no difficulty in the present example.
Indeed, we can compute the symbols explicitly, and their low-frequency asymptotics and high-frequency asymptotics follow directly from the explicit formulae below.
In particular, no derivatives of the cut-off functions are used.
We have
\begin{equation}
  C_D = \int_{\R} \cutoffDirac(\abs{u}) \,du = r_+^2 - r_-^2.
\end{equation}
We choose $C_W := 2(r_+ - r_-)$.
Then,
\begin{equation}
  m_W = \frac{1}{C_W} \int_{\R} \frac{\cutoffWilson(\abs{u})}{\abs{u}} \,du = 1
\end{equation}
and
\begin{equation}
  \kappa_W = \frac{1}{2C_W} \int_{\R} \abs{u}\cutoffWilson(\abs{u}) \,du = \frac{r_+^2 + r_+r_- + r_-^2}{6}.
\end{equation}
We have
\begin{equation}
  \radialsymbolDirac(s) = \int_{\R} \frac{w}{\abs{w}} \frac{\sin(sw)}{\abs{w}} \frac{1}{C_D} \cutoffDirac(\abs{w}) \,dw = \frac{2}{r_+^2-r_-^2} \int_{r_-}^{r_+} \sin(sw) \,dw
\end{equation}
and
\begin{equation}
  \radialsymbolWilson(s) = \int_{\R} \frac{1-\cos(sw)}{\abs{w}} \frac{1}{C_W} \cutoffWilson(\abs{w}) \,dw = \frac{1}{r_+ - r_-} \int_{r_-}^{r_+} (1-\cos(sw)) \,dw
\end{equation}
for $s \ge 0$.

We compute the low-frequency asymptotics.
Since we have $\sin (sw) = sw - s^3w^3/6 + O(s^5)$ uniformly for $w \in [r_-,r_+]$, we have
\begin{equation}
  \radialsymbolDirac(s) = s - \frac{r_+^2+r_-^2}{12} s^3 + O(s^5)
\end{equation}
as $s \to 0$.
Similarly, since we have $1-\cos(sw) = s^2w^2/2 - s^4w^4/24 + O(s^6)$ uniformly for $w \in [r_-,r_+]$, we have
\begin{equation}
  \radialsymbolWilson(s) = \kappa_W s^2 - \frac{r_+^5-r_-^5}{120(r_+-r_-)}s^4 + O(s^6)
\end{equation}
as $s \to 0$.

At high frequencies, we have
\begin{equation}
  \radialsymbolDirac(s) = O(s^{-1})
\end{equation}
and
\begin{equation}
  \radialsymbolWilson(s) = 1 + O(s^{-1})
\end{equation}
as $s \to \infty$.

Moreover, $\radialsymbolWilson(s) > 0$ for every $s > 0$.
Thus, the full symbol $\symbolWilsonDirac$ has no zeros away from the origin.

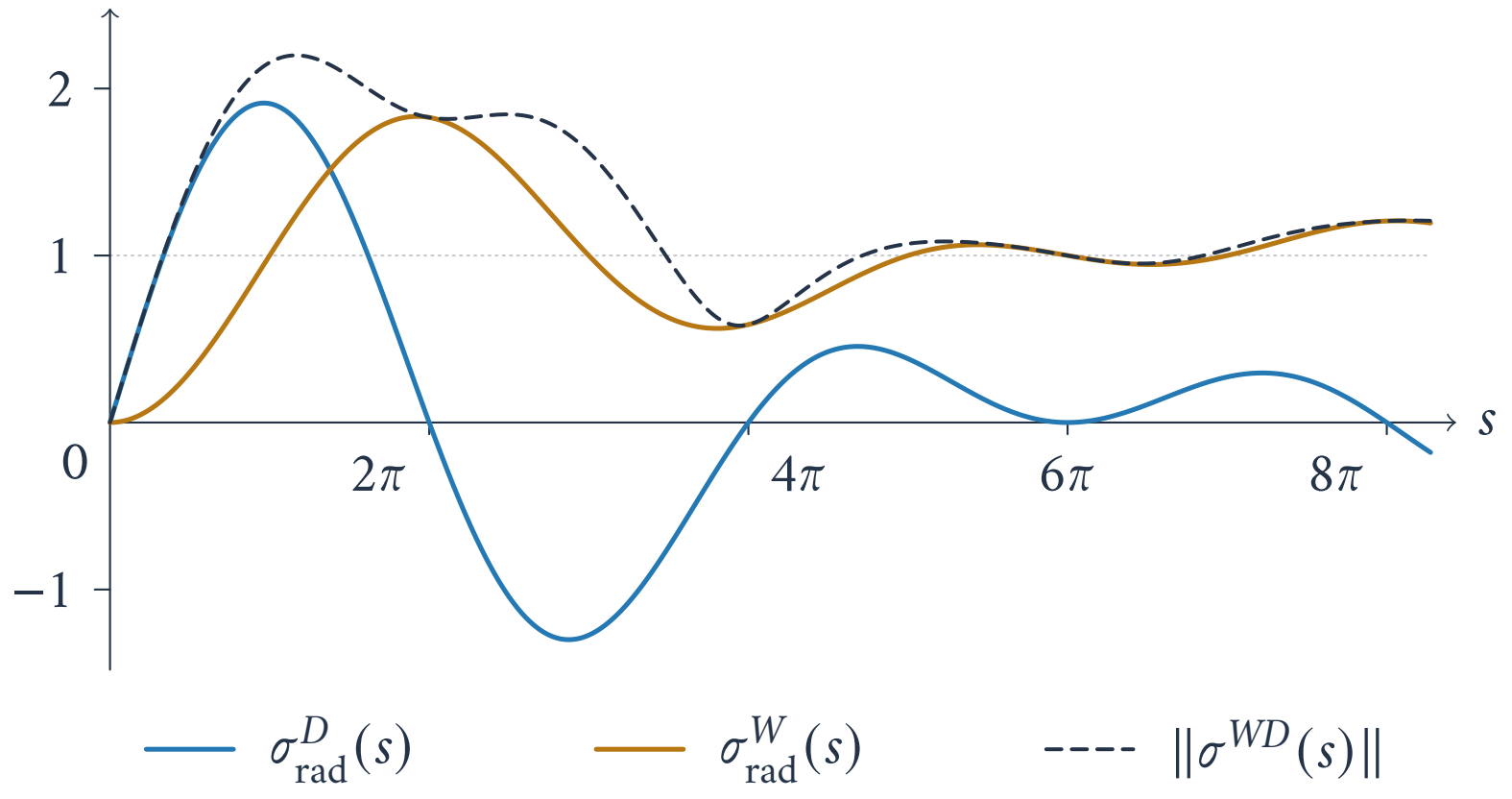
\begin{figure}[htbp]
  \centering
  \begin{tikzpicture}[
    x=1cm,y=1cm,font=\small,text=ink,
    line cap=round,line join=round,
    declare function={
      annulard(\s)=6*(cos(deg(\s/3))-cos(deg(2*\s/3)))
        /max(\s,0.001);
      annularw(\s)=ifthenelse(\s==0,0,
        1-3*(sin(deg(2*\s/3))-sin(deg(\s/3)))
        /max(\s,0.001));
    }
  ]
    \begin{scope}[x=0.35cm,y=1.15cm]
      \draw[ink!25,thin,densely dotted] (0,1) -- (26,1);
      \draw[->,thin,ink]
        (0,0) -- (26.5,0) node[right] {$s$};
      \draw[->,thin,ink] (0,-1.48) -- (0,2.48);
      \node[below left] at (0,0) {$0$};

      \draw[thin,ink] (6.283185,0) -- (6.283185,-0.07)
        node[below left] {$2\pi$};
      \draw[thin,ink] (12.566371,0) -- (12.566371,-0.07)
        node[below right] {$4\pi$};
      \draw[thin,ink] (18.849556,0) -- (18.849556,-0.07)
        node[below] {$6\pi$};
      \draw[thin,ink] (25.132741,0) -- (25.132741,-0.07)
        node[below left] {$8\pi$};

      \foreach \v in {-1,1,2}
        \draw[thin,ink] (0,\v) -- (-0.3,\v)
          node[left] {$\v$};

      \draw[blueedge,line width=0.95pt]
        plot[domain=0:26,samples=401,variable=\s]
        (\s,{annulard(\s)});

      \draw[goldedge,line width=0.95pt]
        plot[domain=0:26,samples=401,variable=\s]
        (\s,{annularw(\s)});

      \draw[ink,line width=0.75pt,densely dashed]
        plot[domain=0:26,samples=401,variable=\s]
        (\s,{sqrt(annulard(\s)^2+annularw(\s)^2)});
    \end{scope}

    \draw[blueedge,line width=0.95pt]
      (0.25,-2.25) -- (0.85,-2.25);
    \node[right] at (0.95,-2.25)
      {$\radialsymbolDirac(s)$};

    \draw[goldedge,line width=0.95pt]
      (3.35,-2.25) -- (3.95,-2.25);
    \node[right] at (4.05,-2.25)
      {$\radialsymbolWilson(s)$};

    \draw[ink,line width=0.75pt,densely dashed]
      (6.45,-2.25) -- (7.05,-2.25);
    \node[right] at (7.15,-2.25)
      {$\norm{\symbolWilsonDirac(s)}$};
  \end{tikzpicture}
\end{figure}

\subsubsection{A Gaussian kernel}

For a second explicit model, we take $\cutoffDirac(r) = r^2e^{-r^2}$ and $\cutoffWilson(r) = r e^{-r^2}$ for $r \ge 0$.
These are not compactly supported, so this example lies outside the finite-range assumptions that will be used on a general manifold; it is included only to show the symbols in closed form.
We define $\symbolDirac$ and $\symbolWilson$ directly by their integral formulae over $V$.

We have $C_D = \sqrt{\pi}/2$ and take $C_W := \sqrt{\pi}$.
Then, we have
\begin{equation}
  \radialsymbolDirac(s) = s e^{-s^2/4}
\end{equation}
and
\begin{equation}
  \radialsymbolWilson(s) = 1 - e^{-s^2/4}
\end{equation}
for $s \ge 0$.
Thus, the eigenvalues of $\symbolWilsonDirac(\eta)$ are
\begin{equation}
  \pm \sqrt{\eta^2 e^{-\eta^2/2} + (1-e^{-\eta^2/4})^2}
\end{equation}
for $\eta \in \R$.
The radial mesoscopic Dirac symbol is asymptotic to $s$ at the origin and tends to zero at high frequency, while the radial mesoscopic Wilson symbol vanishes quadratically at the origin and tends to $1$ at high frequency.
Thus, this example exhibits explicitly the two complementary mechanisms underlying the uniform estimate $\norm{\symbolWilsonDirac(\eta)} \asymp \min\{\abs{\eta},1\}$.

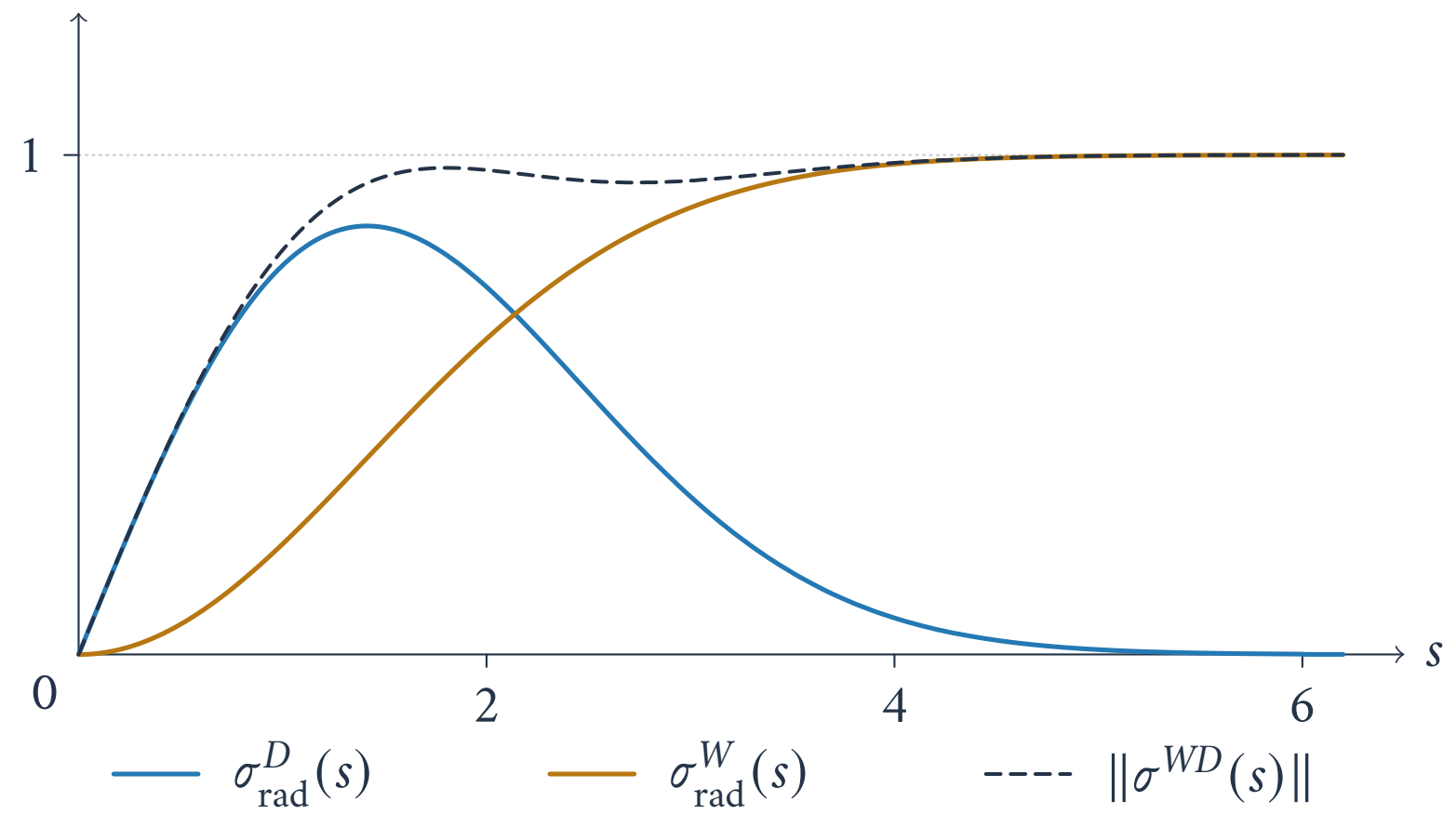
\begin{figure}[htbp]
  \centering
  \begin{tikzpicture}[
    x=1cm,y=1cm,font=\small,text=ink,
    line cap=round,line join=round
  ]
    \begin{scope}[x=1.45cm,y=3.55cm]
      \draw[ink!25,thin,densely dotted] (0,1) -- (6.2,1);
      \draw[->,thin,ink]
        (0,0) -- (6.5,0) node[right] {$s$};
      \draw[->,thin,ink] (0,0) -- (0,1.285);
      \node[below left] at (0,0) {$0$};

      \foreach \t in {2,4,6}
        \draw[thin,ink] (\t,0) -- (\t,-0.025)
          node[below] {$\t$};

      \draw[thin,ink] (0,1) -- (-0.07,1) node[left] {$1$};

      \draw[blueedge,line width=0.95pt]
        plot[domain=0:6.2,samples=241,variable=\s]
        (\s,{\s*exp(-\s^2/4)});

      \draw[goldedge,line width=0.95pt]
        plot[domain=0:6.2,samples=241,variable=\s]
        (\s,{1-exp(-\s^2/4)});

      \draw[ink,line width=0.75pt,densely dashed]
        plot[domain=0:6.2,samples=241,variable=\s]
        (\s,{sqrt(\s^2*exp(-\s^2/2)
          +(1-exp(-\s^2/4))^2)});
    \end{scope}

    \draw[blueedge,line width=0.95pt]
      (0.25,-0.85) -- (0.85,-0.85);
    \node[right] at (0.95,-0.85)
      {$\radialsymbolDirac(s)$};

    \draw[goldedge,line width=0.95pt]
      (3.35,-0.85) -- (3.95,-0.85);
    \node[right] at (4.05,-0.85)
      {$\radialsymbolWilson(s)$};

    \draw[ink,line width=0.75pt,densely dashed]
      (6.45,-0.85) -- (7.05,-0.85);
    \node[right] at (7.15,-0.85)
      {$\norm{\symbolWilsonDirac(s)}$};
  \end{tikzpicture}
\end{figure}

\subsection{The mesoscopic seminorm}

We introduce the mesoscopic seminorm naturally associated with the mesoscopic scale.

\begin{definition}\label{definition: mesoscopic seminorm}
  For $a > 0$ and $\Phi = \sum_{\xi \in \duallattice{\Lattice}} e_{\xi} \Phi_{\xi} \in L^2(\flattorus,\flatmodel{E})$, we define
  \begin{equation}
    \mesoscopicseminorm{\Phi}^2 := \sum_{\substack{\xi \in \duallattice{\Lattice} \\ \abs{\xi} \le 1/a}} \abs{\xi}^2 \norm{\Phi_{\xi}}_S^2 + \sum_{\substack{\xi \in \duallattice{\Lattice} \\ \abs{\xi} > 1/a}} \frac{1}{a^2} \norm{\Phi_{\xi}}_S^2 = \sum_{\xi \in \duallattice{\Lattice}} \min \left\{ \abs{\xi}^2, 1/a^2 \right\} \norm{\Phi_{\xi}}_S^2.
  \end{equation}
  We call $\mesoscopicseminorm{\Phi}$ the \emph{mesoscopic seminorm} of $\Phi$ at scale $a$.
\end{definition}
\noindent
In other words, $\mesoscopicseminorm{\cdot}$ is the Sobolev $L^2_1$-seminorm with the Fourier multiplier capped at scale $a^{-1}$.
This is consistent with the $O(a^{-1})$ operator norm bounds for the mesoscopic operators established in the preceding section.
In particular, this seminorm does not distinguish among frequencies above $a^{-1}$.
They are all measured uniformly at the scale $a^{-1}$.
This point of view is motivated by the idea of a low-energy effective theory, in which degrees of freedom above a prescribed energy scale are not resolved in detail.
Equivalently, we have
\begin{equation}
  \mesoscopicseminorm{\Phi}^2 = \sum_{\xi \in \duallattice{\Lattice}} \left(\frac{\theta(a\abs{\xi})}{a}\right)^2 \norm{\Phi_{\xi}}_S^2
\end{equation}
by setting $\theta(s) := \min \{s,1\}$ for $s \ge 0$.
By definition, we have $\mesoscopicseminorm{\Phi} \le \Ltwonorm{\Phi}/a$.
Also, $\mesoscopicseminorm{\Phi} = 0$ if and only if $\Phi$ is a constant section.
The monotone convergence theorem gives
\begin{equation}
  \lim_{a \to 0} \mesoscopicseminorm{\Phi} = \Ltwonorm{\nabla \Phi}
\end{equation}
for $\Phi \in L^2_1(\flattorus,\flatmodel{E})$.

\subsection{Uniform ellipticity of the symbols}

\begin{proposition}\label{proposition: uniform ellipticity of the symbols}
  There exist constants $0 < \lambda_- \le \lambda_+$ such that
  \begin{equation}
    \lambda_- \theta(\abs{\eta}) \norm{\phi}_S \le \norm{\symbolWilsonDirac(\eta)\phi}_S \le \lambda_+ \theta(\abs{\eta}) \norm{\phi}_S
  \end{equation}
  for any $\eta \in V^*$ and $\phi \in S$.
\end{proposition}
\begin{proof}
  Fix $\eta \in V^*$ and $\phi \in S$.
  The claim is immediate when $\eta = 0$.
  Assume that $\eta \ne 0$.
  By the Clifford relation, we have
  \begin{equation}
    \begin{aligned}
    \symbolWilsonDirac(\eta)^2 ={}& \left( i \radialsymbolDirac(\abs{\eta}) \cliffordmultiplication\big(\eta^{\sharp}/\abs{\eta}\big) + \radialsymbolWilson(\abs{\eta}) \gamma \right)^2 \\
    ={}& \radialsymbolDirac(\abs{\eta})^2 \left(i\cliffordmultiplication\big(\eta^{\sharp}/\abs{\eta}\big)\right)^2 + \radialsymbolDirac(\abs{\eta})\radialsymbolWilson(\abs{\eta}) \cdot i \left( \cliffordmultiplication\big(\eta^{\sharp}/\abs{\eta}\big) \gamma + \gamma \cliffordmultiplication\big(\eta^{\sharp}/\abs{\eta}\big) \right) + \radialsymbolWilson(\abs{\eta})^2 \gamma^2 \\
    ={}& \radialsymbolDirac(\abs{\eta})^2 + \radialsymbolWilson(\abs{\eta})^2.
    \end{aligned}
  \end{equation}
  Hence, we have
  \begin{equation}
    \norm{\symbolWilsonDirac(\eta)\phi}_S^2 = \big( \radialsymbolDirac(\abs{\eta})^2 + \radialsymbolWilson(\abs{\eta})^2 \big) \norm{\phi}_S^2.
  \end{equation}
  Since $\radialsymbolDirac(s) = s + O(s^3)$ and $\radialsymbolWilson(s) = O(s^2)$ as $s \to 0$ by Lemma~\ref{lemma: asymptotics s to 0}, there exist $s_1 > 0$ and $C_1 \ge 1$ such that
  \begin{equation}
    C_1^{-1} \abs{\eta}^2 \le \radialsymbolDirac(\abs{\eta})^2 + \radialsymbolWilson(\abs{\eta})^2 \le C_1 \abs{\eta}^2
  \end{equation}
  for $\abs{\eta} \in (0,s_1)$.
  Since $\radialsymbolDirac(s) \to 0$ and $\radialsymbolWilson(s) \to m_W > 0$ as $s \to \infty$ by Lemma~\ref{lemma: asymptotics s to infinity}, there exist $s_2 \ge s_1$ and $C_2  \ge 1$ such that
  \begin{equation}
    C_2^{-1} \le \radialsymbolDirac(\abs{\eta})^2 + \radialsymbolWilson(\abs{\eta})^2 \le C_2
  \end{equation}
  for $\abs{\eta} \in (s_2,\infty)$.
  On the compact interval $[s_1,s_2]$, the function $s \mapsto \radialsymbolDirac(s)^2 + \radialsymbolWilson(s)^2$ is continuous and strictly positive by Lemma~\ref{lemma: radialsymbolWilson positive}.
  Combining these three cases proves the assertion.
\end{proof}

\begin{theorem}\label{theorem: estimate of WD by seminorm}
  We have
  \begin{equation}
    \lambda_- \mesoscopicseminorm{\Phi} \le \Ltwonorm{\flatmesoscopicWilsonDiracoperator \Phi} \le \lambda_+ \mesoscopicseminorm{\Phi}
  \end{equation}
  for any $\Phi \in L^2(\flattorus,\flatmodel{E})$.
\end{theorem}
\begin{proof}
  Fix $\Phi = \sum_{\xi \in \duallattice{\Lattice}} e_{\xi} \Phi_{\xi} \in  L^2(\flattorus,\flatmodel{E})$.
  By \eqref{equation: Fourier transform of mesoscopic WD operators} and the Plancherel formula, we have
  \begin{equation}
    \Ltwonorm{\flatmesoscopicWilsonDiracoperator \Phi}^2 = \frac{1}{a^2} \sum_{\xi \in \duallattice{\Lattice}} \norm{\symbolWilsonDirac(a\xi) \Phi_{\xi}}_S^2.
  \end{equation}
  By Proposition~\ref{proposition: uniform ellipticity of the symbols}, we have
  \begin{equation}
    \lambda_-^2 \theta(a\abs{\xi})^2 \norm{\Phi_{\xi}}_S^2 \le \norm{\symbolWilsonDirac(a\xi) \Phi_{\xi}}_S^2 \le \lambda_+^2 \theta(a\abs{\xi})^2 \norm{\Phi_{\xi}}_S^2.
  \end{equation}
  Thus, we have
  \begin{equation}
     \lambda_-^2 \mesoscopicseminorm{\Phi}^2 = \lambda_-^2 \sum_{\xi \in \duallattice{\Lattice}} \left(\frac{\theta(a\abs{\xi})}{a}\right)^2 \norm{\Phi_{\xi}}_S^2 \le \Ltwonorm{\flatmesoscopicWilsonDiracoperator \Phi}^2 \le \lambda_+^2 \sum_{\xi \in \duallattice{\Lattice}} \left(\frac{\theta(a\abs{\xi})}{a}\right)^2 \norm{\Phi_{\xi}}_S^2 = \lambda_+^2 \mesoscopicseminorm{\Phi}^2,
  \end{equation}
  which completes the proof.
\end{proof}

In particular, the kernel of $\flatmesoscopicWilsonDiracoperator$ consists precisely of the constant sections.
Thus, the mesoscopic Wilson term removes the zero modes of $\flatmesoscopicDiracoperator$ at non-zero frequencies (the doublers) while preserving the zero-frequency mode.

\subsection{The mesoscopic difference seminorm on a flat torus}

We introduce yet another seminorm defined in terms of difference quotients on the flat torus.
In the next section, we extend this definition to closed Riemannian manifolds.
See Definition~\ref{definition: mesoscopic difference seminorm on manifolds}.
Fix $0 < r_-^{\nabla} < r_+^{\nabla} < 1$.
Let $\cutoffnabla \colon [0,\infty) \to [0,\infty)$ be a non-zero smooth cut-off function such that $\big(\support \cutoffnabla\big) \subset [r_-^{\nabla},r_+^{\nabla}]$.
\begin{definition}\label{definition: flat mesoscopic difference seminorm}
  For $a > 0$ and $\Phi \in L^2(\flattorus,\flatmodel{E})$, we define
  \begin{equation}
    \flatmesoscopicdifferenceseminorm{\cutoffnabla}{\Phi}^2 := \int_{\flattorus} \left( \int_V \norm*{\frac{\Phi(x+v)-\Phi(x)}{\abs{v}}}^2 \cutoffnabla\big(\abs{v}/a\big) \,\frac{dv}{a^n} \right) \,\dvol_g(x).
  \end{equation}
  We call $\flatmesoscopicdifferenceseminorm{\cutoffnabla}{\Phi}$ the \emph{mesoscopic difference seminorm} of $\Phi$ at scale $a$.
\end{definition}

By the change of variables $v=aw$, we have
\begin{equation}
  \flatmesoscopicdifferenceseminorm{\cutoffnabla}{\Phi}^2 = \frac{1}{a^2} \int_{\flattorus} \left( \int_V \norm*{\frac{\Phi(x+aw)-\Phi(x)}{\abs{w}}}^2 \cutoffnabla(\abs{w}) \,dw \right) \,\dvol_g(x).
\end{equation}
Define the symbol $\symbolnabla \colon V^* \to \R$ by
\begin{equation}
  \symbolnabla(\eta) := \int_V \frac{\abs{\exp (i\eta(w)) - 1}^2}{\abs{w}^2} \cutoffnabla(\abs{w}) \,dw.
\end{equation}
The symbol $\symbolnabla$ is non-negative and radial.
We have
\begin{equation}\label{equation: Fourier transform of difference seminorm}
  \flatmesoscopicdifferenceseminorm{\cutoffnabla}{\Phi}^2 = \frac{1}{a^2} \sum_{\xi \in \duallattice{\Lattice}} \symbolnabla(a\xi) \norm{\Phi_{\xi}}_S^2
\end{equation}
for  $\Phi = \sum_{\xi \in \duallattice{\Lattice}} e_{\xi} \Phi_{\xi} \in L^2(\flattorus,\flatmodel{E})$.

\begin{proposition}\label{proposition: uniform estimates for nabla symbol}
  There exist constants $0 < \nu_- \le \nu_+$ such that
  \begin{equation}
    \nu_-^2 \theta(\abs{\eta})^2 \le \symbolnabla(\eta) \le \nu_+^2 \theta(\abs{\eta})^2
  \end{equation}
  for any $\eta \in V^*$.
\end{proposition}
\begin{proof}
  Taylor expansion and radial symmetry give
  \begin{equation}
    \symbolnabla(\eta) = \frac{1}{n} \left( \int_V \cutoffnabla(\abs{w}) \,dw \right) \abs{\eta}^2 + O(\abs{\eta}^4)
  \end{equation}
  as $\eta \to 0$.
  Since $\abs{\exp(it)-1}^2 = 2(1-\cos t)$, we have
  \begin{equation}
    \symbolnabla(\eta) = 2\int_V \frac{\cutoffnabla(\abs{w})}{\abs{w}^2} \,dw - 2\int_V \cos(\eta(w)) \frac{\cutoffnabla(\abs{w})}{\abs{w}^2} \,dw.
  \end{equation}
  By the Riemann--Lebesgue lemma, or the same argument in Lemma~\ref{lemma: asymptotics s to infinity}, we have
  \begin{equation}
    \lim_{\abs{\eta} \to \infty} \symbolnabla(\eta) = 2\int_V \frac{\cutoffnabla(\abs{w})}{\abs{w}^2} \,dw > 0.
  \end{equation}
  Moreover, by the same argument in Lemma~\ref{lemma: radialsymbolWilson positive}, we have $\symbolnabla(\eta) > 0$ for $\eta \ne 0$.
  The rest of the proof is the same as that of Proposition~\ref{proposition: uniform ellipticity of the symbols}.
\end{proof}

\begin{proposition}\label{proposition: equivalence of seminorms}
  We have
  \begin{equation}
    \nu_- \mesoscopicseminorm{\Phi} \le \flatmesoscopicdifferenceseminorm{\cutoffnabla}{\Phi} \le \nu_+ \mesoscopicseminorm{\Phi}
  \end{equation}
  for any $\Phi \in L^2(\flattorus,\flatmodel{E})$.
\end{proposition}
\begin{proof}
  Fix $\Phi = \sum_{\xi \in \duallattice{\Lattice}} e_{\xi} \Phi_{\xi} \in  L^2(\flattorus,\flatmodel{E})$.
  By Proposition~\ref{proposition: uniform estimates for nabla symbol}, we have
  \begin{equation}
    \nu_-^2 \theta(a\abs{\xi})^2 \norm{\Phi_{\xi}}_S^2 \le \symbolnabla(a\xi) \norm{\Phi_{\xi}}_S^2 \le \nu_+^2 \theta(a\abs{\xi})^2 \norm{\Phi_{\xi}}_S^2.
  \end{equation}
  Thus, we have
  \begin{equation}
     \nu_-^2 \mesoscopicseminorm{\Phi}^2 = \nu_-^2 \sum_{\xi \in \duallattice{\Lattice}} \left(\frac{\theta(a\abs{\xi})}{a}\right)^2 \norm{\Phi_{\xi}}_S^2 \le \flatmesoscopicdifferenceseminorm{\cutoffnabla}{\Phi}^2 \le \nu_+^2 \sum_{\xi \in \duallattice{\Lattice}} \left(\frac{\theta(a\abs{\xi})}{a}\right)^2 \norm{\Phi_{\xi}}_S^2 = \nu_+^2 \mesoscopicseminorm{\Phi}^2
  \end{equation}
  by \eqref{equation: Fourier transform of difference seminorm}.
\end{proof}

Combining Theorem~\ref{theorem: estimate of WD by seminorm} and Proposition~\ref{proposition: equivalence of seminorms}, we obtain the elliptic estimate for the mesoscopic Wilson--Dirac operator on the flat torus.
\begin{theorem}\label{theorem: flat elliptic estimate}
  We have
  \begin{equation}
    \frac{\lambda_-}{\nu_+} \flatmesoscopicdifferenceseminorm{\cutoffnabla}{\Phi} \le \Ltwonorm{\flatmesoscopicWilsonDiracoperator \Phi} \le \frac{\lambda_+}{\nu_-} \flatmesoscopicdifferenceseminorm{\cutoffnabla}{\Phi}
  \end{equation}
  for any $\Phi \in L^2(\flattorus,\flatmodel{E})$.
\end{theorem}

\section{Elliptic estimates on closed manifolds}\label{section: Elliptic estimates on closed manifolds}

In this section, we pass from the translation-invariant estimate on a flat torus to an elliptic estimate on a closed Riemannian manifold.
We keep the geometric data $(X,g,E,\gamma,\cliffordmultiplication,\nabla)$, the cut-off functions $\cutoffDirac$, $\cutoffWilson$, and $\cutoffnabla$, and the normalisation constants $C_D$ and $C_W$.

We now extend the mesoscopic difference seminorm from the flat torus to a closed Riemannian manifold.
\begin{definition}\label{definition: mesoscopic difference seminorm on manifolds}
  For $a \in (0,\injectivityradius(X,g)/10)$ and $\Phi \in L^2(X,E)$, we define the covariant mesoscopic difference seminorm by
  \begin{equation}
    \mesoscopicdifferenceseminorm{\Phi}^2 := \int_X \left(\int_X \norm*{\frac{\Phi(x) - \paralleltransport{y}{x}(\Phi(y))}{d_g(x,y)}}^2 \cutoffnabla(d_g(x,y)/a) \,\frac{\dvol_g(y)}{a^n}\right) \,\dvol_g(x).
  \end{equation}
\end{definition}
By the support condition on $\cutoffnabla$, we have
\begin{equation}
  \mesoscopicdifferenceseminorm{\Phi} \le \frac{C}{a} \Ltwonorm{\Phi},
\end{equation}
where $C > 0$ is independent of $a$ and $\Phi$.

The main result of this section is the following elliptic estimate.
\begin{theorem}\label{theorem: global elliptic estimate}
  There exist constants $\constantAell > 0$ and $C > 0$ such that
  \begin{equation}
    \mesoscopicdifferenceseminorm{\Phi} \le C (\Ltwonorm{\mesoscopicWilsonDiracoperator \Phi} + \Ltwonorm{\Phi})
  \end{equation}
  for every $a \in (0,\constantAell)$ and $\Phi \in L^2(X,E)$.
\end{theorem}

We first prove the estimate for sections supported in a sufficiently small ball.
We then establish a uniform commutator estimate and use a partition of unity to prove the theorem.

\subsection{The local elliptic estimate}

\begin{proposition}\label{proposition: local elliptic estimate}
  There exist constants $r > 0$, $\constantAell > 0$, and $\constantlocal > 0$ such that we have
  \begin{equation}
    \mesoscopicdifferenceseminorm{\Phi} \le \constantlocal (\Ltwonorm{\mesoscopicWilsonDiracoperator \Phi} + \Ltwonorm{\Phi})
  \end{equation}
  for $p \in X$, $a \in (0,\constantAell)$, and $\Phi \in C^{\infty}(X,E)$ with $(\support \Phi) \subset B(p,r)$.
\end{proposition}

We compare the operators in normal coordinates with their flat models.
The radius $r$ will be chosen only at the end of the local argument.

\subsubsection{The flat model}\label{subsubsection: the flat model}

Choose an auxiliary cut-off function $\cutoffupper \colon [0,\infty) \to [0,\infty)$ such that $\big(\support \cutoffupper\big) \subset (0,1)$ and $\cutoffupper$ is identically $1$ on a neighbourhood of $\bigl( \support \cutoffDirac \cup \support \cutoffWilson \cup \support \cutoffnabla \bigr)$.

Choose $r_0 \in (0,\injectivityradius(X,g)/10)$ so small that every geodesic ball in $X$ of radius at most $10 r_0$ is strongly convex.
We shall impose further conditions on $r_0$ at the end of this subsection.
Fix $p \in X$.
Set $V := TX_p$ and $S := E_p$.
Choose a linear isometry $I \colon \R^n \to V$, and set $\Lattice := I\big((30r_0\Z)^n \big)$.
We consider $\flattorus := V/\Lattice$ and the product Clifford module bundle $\flatmodel{E} := \flattorus \times S$.
The lattices are so chosen that all these flat tori are isometric.
In particular, $\injectivityradius(\flattorus) = 15r_0 > 10r_0$.
We use the flat operators of the preceding section for this choice of $V$, $\Lattice$, $S$, and $\flatmodel{E}$.

By Theorem~\ref{theorem: flat elliptic estimate}, applied with $\cutoffupper$ in place of $\cutoffnabla$, there exists a constant $\constantflat > 0$ that depends only on the cut-off functions, the normalisation constants, and the dimension such that
\begin{equation}\label{inequality: flat elliptic estimate}
  \flatmesoscopicdifferenceseminorm{\cutoffupper}{\Phi} \le \constantflat \Ltwonorm{\flatmesoscopicWilsonDiracoperator \Phi}
\end{equation}
for any $\Phi \in C^{\infty}(\flattorus,\flatmodel{E})$.

The following estimates follow from standard arguments in normal coordinates and radial gauge with constants uniform in $p \in X$ by compactness of $X$.
\begin{lemma}\label{lemma: normal nbd comparison}
  There exists a constant $\constantgeometric \ge 1$ with the following properties.
  Fix $r \in (0,r_0)$ and $z \ne w \in B(0,3r) \subset V$.
  Set $x := \exp_p(z)$, $y := \exp_p(w)$, and $v_{xy} = \exp_x^{-1}(y)$.
  Then, we have
  \begin{gather}
    \constantgeometric^{-1} \abs{z-w} \le d_g(x,y) \le \constantgeometric \abs{z-w}, \label{inequality: uniform distance comparison} \\
    \constantgeometric^{-1}\le J_p(w)\le \constantgeometric, \label{inequality: Jpw bound} \\
    \constantgeometric^{-1}\le J_x(v_{xy})\le \constantgeometric \label{inequality: Jpv bound}
  \end{gather}
  and
  \begin{gather}
    \abs*{d_g(x,y) - \abs{z-w}} \le \constantgeometric r \abs{z-w}, \label{inequality: uniform distance ratio comparison} \\
    \abs{J_p(w)-1} + \abs*{\frac{J_p(w)}{J_x(v_{xy})}-1} \le \constantgeometric r, \label{inequality: jacobian ratio comparison} \\
    \abs*{\paralleltransportLC{x}{p}(v_{xy})-(w-z)} \le \constantgeometric r \abs{z-w},  \label{inequality: parallel transported vector comparison}\\
    \norm*{\big(\paralleltransport{x}{p} \circ \paralleltransport{y}{x} \circ \paralleltransport{p}{y}\big) - \id_{S}} \le \constantgeometric r d_g(x,y). \label{inequality: holonomy bound}
  \end{gather}
\end{lemma}
\begin{proof}
  All bounds below are uniform in the centre $p \in X$ by compactness of $X$.
  On the fixed normal coordinate balls of radius $3r_0$, the metric is uniformly equivalent to the Euclidean metric and satisfies $g_{ij}(v) = \delta_{ij} + O(\abs{v}^2)$.
  Let $c \colon [0,1] \to X$ be the unique minimising geodesic from $y$ to $x$.
  Since $x, y \in B(p,3r)$ and the ball $B(p,3r)$ is strongly convex, we have $c([0,1]) \subset B(p,3r)$.
  Hence, both $b := \exp_p^{-1} \circ c$ and the coordinate line segment from $z$ to $w$ lie in $B(0,3r)$.
  Comparing the length of $c$ with that of this line segment gives both the uniform distance comparison and
  \begin{equation}
    \frac{d_g(x,y)}{\abs{z-w}}=1+O(r^2).
  \end{equation}
  The Jacobian densities are uniformly bounded above and away from zero on the fixed coordinate balls.
  Moreover, the normal coordinate expansion gives $J_p(w)=1+O(r^2)$ and $J_x(v_{xy})=1+O(r^2)$.
  The asserted estimates for the densities and their ratio follow.

  Set $u(z,w) := \paralleltransportLC{x}{p}(v_{xy})$.
  This map is smooth across the diagonal, with $u(z,z)=0$ and $\partial_w u(0,0)=\id_V$.
  Its second derivatives are uniformly bounded, so $\partial_w u(z,w) = \id_V + O(r)$ for $z,w\in B(0,3r)$.
  Integrating along the segment from $z$ to $w$ yields
  \begin{equation}
    u(z,w)-(w-z) = \int_0^1 \big( \partial_w u(z,z+t(w-z)) - \id_V \big)(w-z) \,dt,
  \end{equation}
  which proves the displacement estimate.

  Finally, trivialise $E$ by radial parallel transport from $p$, and write $\nabla=d+A$.
  In the radial gauge, we have
  \begin{equation}
    A_v(\zeta) = \int_0^1 tF_A(tv)(v,\zeta) \,dt,
  \end{equation}
  where $F_A := dA + A \wedge A$ denotes the local expression of the curvature of $\nabla$, and
  \begin{equation}
    \norm{A_v(\zeta)} \le C\abs{v}\abs{\zeta}.
  \end{equation}
  Recall that $b:=\exp_p^{-1} \circ c$.
  The parallel transport matrix $U(t)$ along $b$ satisfies $U' = -A_b(\dot{b})U$ and $U(0)=\id_S$.
  Since $U(t)$ is unitary, we have
  \begin{equation}
    \norm{U(1)-\id_S} \le \int_0^1 \norm{A_{b(t)}(\dot b(t))} \,dt \le C\int_0^1 \abs{b(t)} \abs{\dot b(t)} \,dt .
  \end{equation}
  We have $\abs{b(t)} < 3r$ and $\int_0^1 \abs{\dot b(t)} \,dt \le C d_g(x,y)$ by the uniform equivalence of the metrics.
  Hence,
  \begin{equation}
    \norm{U(1)-\id_S} \le C r\, d_g(x,y).
  \end{equation}
  Since $U(1)=\paralleltransport{x}{p} \circ \paralleltransport{y}{x} \circ \paralleltransport{p}{y}$, this yields the last estimate.
\end{proof}

The neighbourhood on which $\cutoffupper = 1$ will be used to accommodate the change from $d_g(x,y)/a$ to $\abs{z-w}/a$.
To this end, we further decrease $r_0$ if necessary so that $\cutoffupper(\abs{z-w}/a) = 1$ whenever $d_g(x,y)/a \in \support \cutoffDirac \cup \support \cutoffWilson \cup \support \cutoffnabla$.
Indeed, \eqref{inequality: uniform distance ratio comparison} makes $d_g(x,y)/a$ and $\abs{z-w}/a$ uniformly close on this fixed compact set when $r_0$ is small.
In particular, we have
\begin{equation}\label{inequality: cutoff comparison}
  \cutoffnabla(d_g(x,y)/a) \le \norm{\cutoffnabla}_{L^{\infty}} \cutoffupper(\abs{z-w}/a).
\end{equation}
We now fix $r_0$ for the rest of the paper.

\subsubsection{Comparison of the difference seminorms}

Let $0 < r < r_0$ and assume from now on that $0 < a < r$.
Constants in the following local comparisons may depend on the fixed geometric data, cut-off functions, normalisation constants, and the initial choice of $r_0$, but not on $p$, $r$, $a$, or $\Phi$.

For $\Phi \in L^2(X,E)$ supported in $B(p,2r)$, we define $\flatmodel{\Phi} \in L^2(\flattorus, \flatmodel{E})$ by
\begin{equation}
  \flatmodel{\Phi}(z) := \paralleltransport{\exp_p(z)}{p}\big(\Phi(\exp_p(z))\big)
\end{equation}
for $z \in B(0,3r) \subset V = TX_p$, and extend it by zero to $\flattorus$.

\begin{lemma}\label{lemma: Ltwo norm comparison}
  There exists a constant $\constantLtwonorm \ge 1$ such that
  \begin{equation}\label{inequality: comparison of L2 norms}
    \constantLtwonorm^{-1} \Ltwonorm[L^2(X)]{\Phi} \le \Ltwonorm[L^2(\flattorus)]{\flatmodel{\Phi}} \le \constantLtwonorm \Ltwonorm[L^2(X)]{\Phi}
  \end{equation}
  for $\Phi \in L^2(X,E)$ supported in $B(p,2r)$.
\end{lemma}
\begin{proof}
  By unitarity of parallel transport and the change of variables $x=\exp_p(z)$, we have
  \begin{equation}
    \Ltwonorm[L^2(X)]{\Phi}^2 = \int_{B(0,3r)} \norm{\flatmodel{\Phi}(z)}^2 J_p(z)\,dz.
  \end{equation}
  By \eqref{inequality: Jpw bound}, the result follows with $\constantLtwonorm := \constantgeometric^{1/2}$.
\end{proof}

\begin{lemma}\label{lemma: comparison of difference seminorms}
  There exists a constant $\constantnabla > 0$ such that
  \begin{equation}
    \mesoscopicdifferenceseminorm{\Phi} \le \constantnabla \big( \flatmesoscopicdifferenceseminorm{\cutoffupper}{\flatmodel{\Phi}} + \Ltwonorm{\Phi} \big)
  \end{equation}
  for $\Phi \in C^{\infty}(X,E)$ supported in $B(p,2r)$.
\end{lemma}
\begin{proof}
  Let $\Phi \in C^{\infty}(X,E)$ be supported in $B(p,2r)$.
  Let $z \ne w \in B(0,3r) \subset V = TX_p$, and set $x := \exp_p(z)$ and $y := \exp_p(w)$.
  We have
  \begin{equation}
    \begin{aligned}
    \Phi(x) - \paralleltransport{y}{x}(\Phi(y)) ={}& \paralleltransport{p}{x}\big(\flatmodel{\Phi}(z) - \big(\paralleltransport{x}{p} \circ \paralleltransport{y}{x} \circ \paralleltransport{p}{y}\big)\big(\flatmodel{\Phi}(w)\big) \big) \\
    ={}& \paralleltransport{p}{x}\big(\big(\flatmodel{\Phi}(z) - \flatmodel{\Phi}(w)\big) - \big(\big(\paralleltransport{x}{p} \circ \paralleltransport{y}{x} \circ \paralleltransport{p}{y} - \id_S \big) \big(\flatmodel{\Phi}(w)\big)\big)\big).
    \end{aligned}
  \end{equation}
  Hence, we have
  \begin{equation}
    \begin{aligned}
      \norm*{\frac{\Phi(x) - \paralleltransport{y}{x}(\Phi(y))}{d_g(x,y)}}^2 \le{}& 2 \constantgeometric^2 \frac{\norm{\flatmodel{\Phi}(z)-\flatmodel{\Phi}(w)}^2}{\abs{z-w}^2} + 2 \constantgeometric^2 r^2 \norm*{\flatmodel{\Phi}(w)}^2
    \end{aligned}
  \end{equation}
  by Lemma~\ref{lemma: normal nbd comparison}.
  We also have $J_p(z) J_p(w) \le \constantgeometric^2$ by Lemma~\ref{lemma: normal nbd comparison}.
  We may replace $\cutoffnabla$ with $\cutoffupper$ by \eqref{inequality: cutoff comparison}.
  Since $\Phi$ is supported in $B(p,2r)$ and $a < r$, the integrand vanishes unless $x,y \in B(p,3r)$.
  Thus, we have
  \begin{equation}
    \begin{aligned}
    \mesoscopicdifferenceseminorm{\Phi}^2 ={}& \int_X \left(\int_X \norm*{\frac{\Phi(x) - \paralleltransport{y}{x}(\Phi(y))}{d_g(x,y)}}^2 \cutoffnabla(d_g(x,y)/a) \,\frac{\dvol_g(y)}{a^n}\right) \,\dvol_g(x) \\
    \le{}& \begin{multlined}[t]
      2 \constantgeometric^4 \norm{\cutoffnabla}_{L^{\infty}} \int_{B(0,3r)} \int_{B(0,3r)} \frac{\norm{\flatmodel{\Phi}(z)-\flatmodel{\Phi}(w)}^2}{\abs{z-w}^2} \cutoffupper(\abs{z-w}/a) \,\frac{dw}{a^n} \,dz \\
      + 2 \constantgeometric^4 \norm{\cutoffnabla}_{L^{\infty}} r^2 \int_{B(0,3r)} \int_{B(0,3r)} \norm*{\flatmodel{\Phi}(w)}^2 \cutoffupper(\abs{z-w}/a) \,\frac{dw}{a^n} \,dz.
    \end{multlined}
    \end{aligned}
  \end{equation}
  For the first term, we have
  \begin{equation}
    \begin{aligned}
    & \int_{B(0,3r)} \int_{B(0,3r)} \frac{\norm{\flatmodel{\Phi}(z)-\flatmodel{\Phi}(w)}^2}{\abs{z-w}^2} \cutoffupper(\abs{z-w}/a) \,\frac{dw}{a^n} \,dz \\
    ={}& \int_{B(0,3r)} \int_{B(0,3r)-z} \frac{\norm{\flatmodel{\Phi}(z+v)-\flatmodel{\Phi}(z)}^2}{\abs{v}^2} \cutoffupper(\abs{v}/a) \,\frac{dv}{a^n} \,dz \\
    \le{}& \int_{\flattorus} \int_{V} \frac{\norm{\flatmodel{\Phi}(z+v)-\flatmodel{\Phi}(z)}^2}{\abs{v}^2} \cutoffupper(\abs{v}/a) \,\frac{dv}{a^n} \,dz = \flatmesoscopicdifferenceseminorm{\cutoffupper}{\flatmodel{\Phi}}^2.
    \end{aligned}
  \end{equation}
  For the next term, the change of variables $z-w = av$ gives
  \begin{equation}
    \int_V \cutoffupper(\abs{z-w}/a) \,\frac{dz}{a^n} = \int_V \cutoffupper(\abs{v}) \,dv =: M_{\cutoffupper}
  \end{equation}
  independently of $a$ and $w$.
  Hence, we have
  \begin{equation}
    \begin{aligned}
    & \int_{B(0,3r)} \int_{B(0,3r)} \norm*{\flatmodel{\Phi}(w)}^2 \cutoffupper(\abs{z-w}/a) \,\frac{dw}{a^n} \,dz \\
    ={}& \int_{B(0,3r)} \norm*{\flatmodel{\Phi}(w)}^2 \left( \int_{B(0,3r)} \cutoffupper(\abs{z-w}/a) \,\frac{dz}{a^n} \right) \,dw \\
    \le{}& M_{\cutoffupper} \int_{B(0,3r)} \norm*{\flatmodel{\Phi}(w)}^2 \,dw \\
    \le{}& M_{\cutoffupper} \int_{\flattorus} \norm*{\flatmodel{\Phi}(w)}^2 \,dw = M_{\cutoffupper} \Ltwonorm[L^2(\flattorus)]{\flatmodel{\Phi}}^2 \\
    \le{}& \constantLtwonorm^2 M_{\cutoffupper} \Ltwonorm[L^2(X)]{\Phi}^2,
    \end{aligned}
  \end{equation}
  where we used \eqref{inequality: comparison of L2 norms} in the last line.
  Set $\constantnabla^2 := \max \{ 2 \constantgeometric^4 \norm{\cutoffnabla}_{L^{\infty}}, 2 \constantLtwonorm^2 \constantgeometric^4 \norm{\cutoffnabla}_{L^{\infty}} r_0^2 M_{\cutoffupper} \}$.
\end{proof}

\subsubsection{Freezing the coefficients}

We note that, if $(\support \Phi) \subset B(p,r)$, then we have
\begin{equation}
  \left( \support \big(\mesoscopicWilsonDiracoperator \Phi\big) \right) \subset B(p,r+a) \subset B(p,2r)
\end{equation}
and
\begin{equation}
  \left( \support \big(\flatmesoscopicWilsonDiracoperator \flatmodel{\Phi}\big) \right) \subset B(0,r+a) \subset B(0,2r).
\end{equation}
In particular, $\flatmodel{\left(\mesoscopicWilsonDiracoperator \Phi\right)}$ is defined by the same radial trivialisation.
The following lemma compares the mesoscopic Wilson--Dirac operators on $X$ and $\flattorus$ using radial trivialisation.
\begin{lemma}\label{lemma: comparison of two freezing operators}
  There exists a constant $\constantWD > 0$ such that
  \begin{equation}
    \norm*{\flatmesoscopicWilsonDiracoperator \flatmodel{\Phi} - \flatmodel{\left(\mesoscopicWilsonDiracoperator \Phi\right)}}_{L^2(\flattorus)} \le \constantWD \left( r\flatmesoscopicdifferenceseminorm{\cutoffupper}{\flatmodel{\Phi}} + \Ltwonorm{\Phi} \right)
  \end{equation}
  for $\Phi \in C^{\infty}(X,E)$ supported in $B(p,r)$.
\end{lemma}
\begin{proof}
  Let $\Phi \in C^{\infty}(X,E)$ be supported in $B(p,r)$.
  Both $\flatmesoscopicWilsonDiracoperator \flatmodel{\Phi}$ and $\flatmodel{\left(\mesoscopicWilsonDiracoperator \Phi\right)}$ are supported in $B(0,2r)$.
  Fix $z \in B(0,2r) \subset TX_p$ and set $x := \exp_p(z) \in X$.

  We first rewrite $(\mesoscopicDiracoperator \Phi)(x)$ as an integral in normal coordinates centred at $p$.
  For $v \in B(0,a) \subset TX_x$, we have
  \begin{equation}
    d_g(p,\exp_x(v)) \le d_g(p,x) + d_g(x,\exp_x(v)) = \abs{z} + \abs{v} < 2r + a < 3r.
  \end{equation}
  Thus, the change of variables $w := \exp_p^{-1}(\exp_x(v))$ takes values in $B(0,3r) \subset TX_p$.
  Since $J_x(v) \,dv = J_p(w) \,dw$, it gives
  \begin{equation}
    dv = \frac{J_p(w)}{J_x(v)} \,dw.
  \end{equation}
  For $w \in B(0,3r)$ with $w \ne z$, set $y := \exp_p(w)$ and $v_{xy} = \exp_x^{-1}(y)$.
  Then, $\abs{v_{xy}} = d_g(x,y)$.
  Thus, we have
  \begin{equation}
    \begin{aligned}
    (\mesoscopicDiracoperator \Phi)(x) ={}& \frac{1}{C_D} \int_{TX_x} \cliffordmultiplication_x\left(\frac{v}{\abs{v}}\right) \frac{\paralleltransport{\exp_x(v)}{x}(\Phi(\exp_x(v))) - \Phi(x)}{\abs{v}} \cutoffDirac(\abs{v}/a) \,\frac{dv}{a^n} \\
    ={}& \frac{1}{C_D} \int_{B(0,3r) \subset TX_p} \cliffordmultiplication_x\left(\frac{v_{xy}}{\abs{v_{xy}}}\right) \frac{\paralleltransport{y}{x}(\Phi(y)) - \Phi(x)}{\abs{v_{xy}}} \cutoffDirac(\abs{v_{xy}}/a) \frac{J_p(w)}{J_x(v_{xy})} \,\frac{dw}{a^n} \\
    ={}& \frac{1}{C_D a^n} \int_{B(0,3r) \subset TX_p} \frac{\cliffordmultiplication_x(v_{xy})\left(\paralleltransport{y}{x}(\Phi(y)) - \Phi(x)\right)}{d_g(x,y)^2} \cutoffDirac(d_g(x,y)/a) \frac{J_p(w)}{J_x(v_{xy})} \,dw.
    \end{aligned}
  \end{equation}
  The term involving $\Phi(x)$ vanishes by oddness, but keeping it is useful for comparison.

  We next apply the parallel transport $\paralleltransport{x}{p}$ to this expression.
  Since the connection is compatible with Clifford multiplication, we have
  \begin{equation}
    \paralleltransport{x}{p} \circ \cliffordmultiplication_x(v_{xy}) = \cliffordmultiplication_p \big(\paralleltransportLC{x}{p}(v_{xy})\big) \circ \paralleltransport{x}{p}.
  \end{equation}
  We also have $\Phi(y) = \paralleltransport{p}{y}\big(\flatmodel{\Phi}(w)\big)$ and $\paralleltransport{x}{p}(\Phi(x)) = \flatmodel{\Phi}(z)$.
  Hence, we have
  \begin{equation}
    \begin{aligned}
    & \flatmodel{\left(\mesoscopicDiracoperator \Phi\right)}(z) \\
    ={}& \paralleltransport{x}{p}\left(\mesoscopicDiracoperator \Phi(x)\right) \\
    ={}& \frac{1}{C_D a^n} \int_{B(0,3r)} \frac{\cliffordmultiplication_p \big(\paralleltransportLC{x}{p}(v_{xy})\big) \left( \big(\paralleltransport{x}{p} \circ \paralleltransport{y}{x} \circ \paralleltransport{p}{y}\big)\big(\flatmodel{\Phi}(w)\big) - \flatmodel{\Phi}(z) \right)}{d_g(x,y)^2} \cutoffDirac(d_g(x,y)/a) \frac{J_p(w)}{J_x(v_{xy})} \,dw \\
    ={}& \frac{1}{C_D a^n} \int_{B(0,3r)} \frac{\cliffordmultiplication_p \big(\paralleltransportLC{x}{p}(v_{xy})\big) }{d_g(x,y)^2} \cutoffDirac(d_g(x,y)/a) \frac{J_p(w)}{J_x(v_{xy})} \big( \big(\paralleltransport{x}{p} \circ \paralleltransport{y}{x} \circ \paralleltransport{p}{y}\big)\big(\flatmodel{\Phi}(w)\big) - \flatmodel{\Phi}(z) \big) \,dw. \label{equation: flat(D Phi)}
    \end{aligned}
  \end{equation}
  For the operator on the torus, by the substitution $v = w-z$, we have
  \begin{equation}\label{equation: flat D (flat Phi)}
    \flatmesoscopicDiracoperator \flatmodel{\Phi}(z) = \frac{1}{C_D a^n} \int_{B(0,3r)} \frac{\cliffordmultiplication_p(w-z)}{\abs{z-w}^2} \cutoffDirac(\abs{z-w}/a) \big(\flatmodel{\Phi}(w) - \flatmodel{\Phi}(z)\big) \,dw.
  \end{equation}
  Note that $\abs{v} < a$ implies $\abs{w} < 2r+a < 3r$.
  In the same way, we have
  \begin{multline}
    \flatmodel{\left(\mesoscopicWilsonterm \Phi\right)}(z) = \frac{1}{C_W a^n} \int_{B(0,3r)} \frac{1}{d_g(x,y)} \cutoffWilson(d_g(x,y)/a) \frac{J_p(w)}{J_x(v_{xy})} \\
    \big( \flatmodel{\Phi}(z) - \big(\paralleltransport{x}{p} \circ \paralleltransport{y}{x} \circ \paralleltransport{p}{y}\big)\big(\flatmodel{\Phi}(w)\big) \big) \,dw
  \end{multline}
  and
  \begin{equation}
    \flatmesoscopicWilsonterm \flatmodel{\Phi}(z) = \frac{1}{C_W a^n} \int_{B(0,3r)} \frac{1}{\abs{z-w}} \cutoffWilson(\abs{z-w}/a) \big(\flatmodel{\Phi}(z) - \flatmodel{\Phi}(w)\big)\,dw.
  \end{equation}

  We next subtract the two expressions for the Dirac term and separate the coefficient difference from the error in parallel transport.
  We use the identity
  \begin{equation}
    \big(\paralleltransport{x}{p} \circ \paralleltransport{y}{x} \circ \paralleltransport{p}{y}\big)\big(\flatmodel{\Phi}(w)\big) - \flatmodel{\Phi}(z) = \big(\flatmodel{\Phi}(w)-\flatmodel{\Phi}(z)\big) + \big(\paralleltransport{x}{p} \circ \paralleltransport{y}{x} \circ \paralleltransport{p}{y}-\id_S\big)\big(\flatmodel{\Phi}(w)\big).
  \end{equation}
  Substituting this into \eqref{equation: flat(D Phi)} and subtracting \eqref{equation: flat D (flat Phi)} yields
  \begin{equation}\label{equation: flat(D Phi) minus flat D (flat Phi)}
    \begin{aligned}
    & \flatmodel{\left(\mesoscopicDiracoperator \Phi\right)}(z) - \flatmesoscopicDiracoperator \flatmodel{\Phi}(z) \\
    ={}& \begin{multlined}[t]
      \frac{1}{C_D a^n} \int_{B(0,3r)} \left( \frac{\cliffordmultiplication_p \big(\paralleltransportLC{x}{p}(v_{xy})\big) }{d_g(x,y)^2} \cutoffDirac(d_g(x,y)/a) \frac{J_p(w)}{J_x(v_{xy})} - \frac{\cliffordmultiplication_p(w-z)}{\abs{z-w}^2} \cutoffDirac(\abs{z-w}/a) \right) \\ \big(\flatmodel{\Phi}(w) - \flatmodel{\Phi}(z)\big) \,dw
    \end{multlined} \\
    +{}& \frac{1}{C_D a^n} \int_{B(0,3r)} \left( \frac{\cliffordmultiplication_p \big(\paralleltransportLC{x}{p}(v_{xy})\big) }{d_g(x,y)^2} \cutoffDirac(d_g(x,y)/a) \frac{J_p(w)}{J_x(v_{xy})} \right) \big(\paralleltransport{x}{p} \circ \paralleltransport{y}{x} \circ \paralleltransport{p}{y} - \id_S\big) \big(\flatmodel{\Phi}(w)\big) \,dw.
    \end{aligned}
  \end{equation}
  The first integral contains the coefficient difference multiplying the ordinary difference of $\flatmodel{\Phi}$.
  The second integral contains the original coefficient multiplying the error in parallel transport.
  We estimate these two integrals separately.

  We first estimate the coefficient difference in the first integral.
  We claim that
  \begin{equation}\label{inequality: coefficients comparison}
    \norm*{\frac{\cliffordmultiplication_p \big(\paralleltransportLC{x}{p}(v_{xy})\big)}{d_g(x,y)^2} \cutoffDirac(d_g(x,y)/a) \frac{J_p(w)}{J_x(v_{xy})} - \frac{\cliffordmultiplication_p(w-z) }{\abs{z-w}^2} \cutoffDirac(\abs{z-w}/a)} \le C_1 r\frac{\cutoffupper(\abs{z-w}/a)}{\abs{z-w}}.
  \end{equation}
  To see this, we decompose the coefficient difference as
  \begin{multline}
      \frac{\cliffordmultiplication_p \big(\paralleltransportLC{x}{p}(v_{xy})\big)}{d_g(x,y)^2} \cutoffDirac(d_g(x,y)/a) \left(\frac{J_p(w)}{J_x(v_{xy})} - 1\right) \\
      + \cliffordmultiplication_p \big(\paralleltransportLC{x}{p}(v_{xy})\big) \left(\frac{1}{d_g(x,y)^2} - \frac{1}{\abs{z-w}^2}\right) \cutoffDirac(d_g(x,y)/a) \\
      + \frac{\cliffordmultiplication_p \big(\paralleltransportLC{x}{p}(v_{xy}) - (w-z)\big)}{\abs{z-w}^2} \cutoffDirac(d_g(x,y)/a) \\
      + \frac{\cliffordmultiplication_p(w-z)}{\abs{z-w}^2} \left(\cutoffDirac(d_g(x,y)/a) - \cutoffDirac(\abs{z-w}/a)\right).
  \end{multline}
  The four terms arise from successively replacing the Jacobian ratio, the distance in the denominator, the displacement vector, and the argument of the cut-off function by their counterparts in the flat model.
  We note that $\cutoffupper(\abs{z-w}/a) = 1$ if any of these four terms is not zero.
  We have $\norm{\cliffordmultiplication_p \big(\paralleltransportLC{x}{p}(v_{xy})\big)} = d_g(x,y)$.
  By \eqref{inequality: uniform distance comparison} and \eqref{inequality: jacobian ratio comparison}, we have
  \begin{equation}
    \norm*{\frac{\cliffordmultiplication_p \big(\paralleltransportLC{x}{p}(v_{xy})\big)}{d_g(x,y)^2} \cutoffDirac(d_g(x,y)/a) \left(\frac{J_p(w)}{J_x(v_{xy})} - 1\right)} \le \frac{d_g(x,y)}{d_g(x,y)^2} \cdot \norm{\cutoffDirac}_{L^{\infty}} \cdot \constantgeometric r \le \frac{\constantgeometric^2 \norm{\cutoffDirac}_{L^{\infty}} r}{\abs{z-w}}.
  \end{equation}
  By \eqref{inequality: uniform distance comparison} and \eqref{inequality: uniform distance ratio comparison}, we have
  \begin{equation}
    d_g(x,y) \abs*{\frac{1}{d_g(x,y)^2}-\frac{1}{\abs{z-w}^2}} = \frac{\abs*{d_g(x,y)-\abs{z-w}}}{\abs{z-w}^2} \left(1+\frac{\abs{z-w}}{d_g(x,y)}\right) \le\frac{\constantgeometric r}{\abs{z-w}} \cdot (1+\constantgeometric).
  \end{equation}
  Thus, we have
  \begin{equation}
    \norm*{\cliffordmultiplication_p \big(\paralleltransportLC{x}{p}(v_{xy})\big) \left(\frac{1}{d_g(x,y)^2} - \frac{1}{\abs{z-w}^2}\right) \cutoffDirac(d_g(x,y)/a)} \le \frac{(1+\constantgeometric)\constantgeometric \norm{\cutoffDirac}_{L^{\infty}} r}{\abs{z-w}}.
  \end{equation}
  By \eqref{inequality: parallel transported vector comparison}, we have
  \begin{equation}
    \norm*{\frac{\cliffordmultiplication_p \big(\paralleltransportLC{x}{p}(v_{xy}) - (w-z)\big)}{\abs{z-w}^2} \cutoffDirac(d_g(x,y)/a)} \le \frac{\constantgeometric \norm{\cutoffDirac}_{L^{\infty}} r}{\abs{z-w}}.
  \end{equation}
  For the fourth term, we note that, if either $\cutoffDirac(d_g(x,y)/a)$ or $\cutoffDirac(\abs{z-w}/a)$ is not zero, then $\cutoffupper(\abs{z-w}/a) = 1$ and $\abs{z-w}/a < 1$.
  Thus, by \eqref{inequality: uniform distance ratio comparison} and the mean value theorem, we have
  \begin{equation}
    \begin{aligned}
    \norm*{\frac{\cliffordmultiplication_p(w-z)}{\abs{z-w}^2} \left(\cutoffDirac(d_g(x,y)/a) - \cutoffDirac(\abs{z-w}/a)\right)} \le{}& \frac{1}{\abs{z-w}} \cdot \norm{(\cutoffDirac)'}_{L^{\infty}} \frac{\abs*{d_g(x,y)- \abs{z-w}}}{a} \\
    \le{}& \frac{\constantgeometric \norm{(\cutoffDirac)'}_{L^{\infty}} r}{\abs{z-w}}.
    \end{aligned}
  \end{equation}
  Thus, \eqref{inequality: coefficients comparison} holds with the constant
  \begin{equation}
    C_1 := \constantgeometric^2 \norm{\cutoffDirac}_{L^{\infty}} + (1+\constantgeometric)\constantgeometric \norm{\cutoffDirac}_{L^{\infty}} + \constantgeometric \norm{\cutoffDirac}_{L^{\infty}} + \constantgeometric \norm{(\cutoffDirac)'}_{L^{\infty}}.
  \end{equation}
  In the same way, we have
  \begin{equation}
    \abs*{\frac{\cutoffWilson(d_g(x,y)/a)}{d_g(x,y)} \frac{J_p(w)}{J_x(v_{xy})} - \frac{\cutoffWilson(\abs{z-w}/a)}{\abs{z-w}}} \le C_2 r \frac{\cutoffupper(\abs{z-w}/a)}{\abs{z-w}},
  \end{equation}
  where
  \begin{equation}
    C_2 := 2 \constantgeometric^2 \norm{\cutoffWilson}_{L^{\infty}} + \constantgeometric \norm{(\cutoffWilson)'}_{L^{\infty}}.
  \end{equation}

  We next estimate the second integral, which contains the error in parallel transport.
  By \eqref{inequality: Jpw bound} and \eqref{inequality: Jpv bound}, we have
  \begin{equation}
    \norm*{\frac{\cliffordmultiplication_p \big(\paralleltransportLC{x}{p}(v_{xy})\big) }{d_g(x,y)^2} \cutoffDirac(d_g(x,y)/a) \frac{J_p(w)}{J_x(v_{xy})}} \le \constantgeometric^2 \norm{\cutoffDirac}_{L^{\infty}} \frac{\cutoffupper(\abs{z-w}/a)}{d_g(x,y)}.
  \end{equation}
  Then, by \eqref{inequality: holonomy bound}, we have
  \begin{equation}
    \begin{aligned}
    & \norm*{\left( \frac{\cliffordmultiplication_p \big(\paralleltransportLC{x}{p}(v_{xy})\big) }{d_g(x,y)^2} \cutoffDirac(d_g(x,y)/a) \frac{J_p(w)}{J_x(v_{xy})} \right) \big(\paralleltransport{x}{p} \circ \paralleltransport{y}{x} \circ \paralleltransport{p}{y} - \id_S\big) \big(\flatmodel{\Phi}(w)\big)} \\
    \le{}& \constantgeometric^2 \norm{\cutoffDirac}_{L^{\infty}} \frac{\cutoffupper(\abs{z-w}/a)}{d_g(x,y)} \cdot \constantgeometric r d_g(x,y) \cdot \norm{\flatmodel{\Phi}(w)} \\
    ={}& \constantgeometric^3 \norm{\cutoffDirac}_{L^{\infty}} r \cutoffupper(\abs{z-w}/a) \norm{\flatmodel{\Phi}(w)}.
    \end{aligned}
  \end{equation}

  Combining the two estimates in \eqref{equation: flat(D Phi) minus flat D (flat Phi)} yields
  \begin{equation}
    \begin{aligned}
    & \norm*{\flatmodel{\left(\mesoscopicDiracoperator \Phi\right)}(z) - \flatmesoscopicDiracoperator \flatmodel{\Phi}(z)} \\
    \le{}&
    \begin{multlined}[t]
      \frac{C_1 r}{C_D a^n} \int_{B(0,3r)} \frac{\cutoffupper(\abs{z-w}/a)}{\abs{z-w}} \norm*{\flatmodel{\Phi}(w) - \flatmodel{\Phi}(z)} \,dw \\ + \frac{\constantgeometric^3 \norm{\cutoffDirac}_{L^{\infty}} r}{C_D a^n} \int_{B(0,3r)} \cutoffupper(\abs{z-w}/a) \norm{\flatmodel{\Phi}(w)} \,dw
    \end{multlined} \\
    \le{}& \frac{C_1 r}{C_D} \int_V \frac{\cutoffupper(\abs{v}/a)}{a^n} \norm*{\frac{\flatmodel{\Phi}(z+v) - \flatmodel{\Phi}(z)}{\abs{v}}} \,dv + \frac{\constantgeometric^3 \norm{\cutoffDirac}_{L^{\infty}} r}{C_D} \int_V \frac{\cutoffupper(\abs{v}/a)}{a^n} \norm*{\flatmodel{\Phi}(z+v)} \,dv.
    \end{aligned}
  \end{equation}
  The Wilson term satisfies
  \begin{multline}
    \norm*{\flatmodel{\left(\mesoscopicWilsonterm \Phi\right)}(z) - \flatmesoscopicWilsonterm \flatmodel{\Phi}(z)} \\
    \le \frac{C_2 r}{C_W} \int_V \frac{\cutoffupper(\abs{v}/a)}{a^n} \norm*{\frac{\flatmodel{\Phi}(z+v) - \flatmodel{\Phi}(z)}{\abs{v}}} \,dv + \frac{\constantgeometric^3 \norm{\cutoffWilson}_{L^{\infty}} r}{C_W} \int_V \frac{\cutoffupper(\abs{v}/a)}{a^n} \norm*{\flatmodel{\Phi}(z+v)} \,dv.
  \end{multline}
  So far we have assumed that $z\in B(0,2r)$.
  However, the support conditions on $\Phi$ and the cut-off functions imply that both sections on the left-hand side vanish outside $B(0,2r)$.
  Thus, the inequalities hold for every $z \in \flattorus$.

  Recall that $M_{\cutoffupper} = \int_V \cutoffupper(\abs{v}) \,dv$.
  By the Cauchy--Schwarz inequality, we have
  \begin{equation}
    \int_{\flattorus} \left( \int_V \frac{\cutoffupper(\abs{v}/a)}{a^n} \norm*{\frac{\flatmodel{\Phi}(z+v) - \flatmodel{\Phi}(z)}{\abs{v}}} \,dv \right)^2 \,dz \le M_{\cutoffupper} \flatmesoscopicdifferenceseminorm{\cutoffupper}{\flatmodel{\Phi}}^2.
  \end{equation}
  By the same inequality and \eqref{inequality: comparison of L2 norms}, we have
  \begin{equation}
    \int_{\flattorus} \left( \int_V \frac{\cutoffupper(\abs{v}/a)}{a^n} \norm{\flatmodel{\Phi}(z+v)} \,dv \right)^2 \,dz \le M_{\cutoffupper}^2 \Ltwonorm[L^2(\flattorus)]{\flatmodel{\Phi}}^2 \le \constantLtwonorm^2 M_{\cutoffupper}^2 \Ltwonorm[L^2(X)]{\Phi}^2.
  \end{equation}
  The same inequalities hold for the Wilson term.
  Since the connection preserves the grading, we have $\flatmodel{(\gamma \Phi)} = \gamma_p \flatmodel{\Phi}$.
  Recall that $r \le r_0$.
  Thus, we have
  \begin{equation}
    \begin{aligned}
    & \norm*{\flatmesoscopicWilsonDiracoperator \flatmodel{\Phi} - \flatmodel{\left(\mesoscopicWilsonDiracoperator \Phi\right)}}_{L^2(\flattorus)} \\
    \le{}& M_{\cutoffupper}^{1/2} \left(\frac{C_1}{C_D} + \frac{C_2}{C_W}\right) r \flatmesoscopicdifferenceseminorm{\cutoffupper}{\flatmodel{\Phi}} + \constantLtwonorm M_{\cutoffupper} \constantgeometric^3 \left(\frac{\norm{\cutoffDirac}_{L^{\infty}}}{C_D}+\frac{\norm{\cutoffWilson}_{L^{\infty}}}{C_W}\right) r_0 \Ltwonorm[L^2(X)]{\Phi}.
    \end{aligned}
  \end{equation}
  This completes the proof.
\end{proof}

\subsubsection{Completion of the local elliptic estimate}

\begin{proof}[Proof of Proposition~\ref{proposition: local elliptic estimate}]
  We choose $r \in (0,r_0)$ so that
  \begin{equation}
    r \le \frac{1}{2 \constantflat \constantWD}
  \end{equation}
  and set $\constantAell := r/4$.
  Fix $p \in X$ and $a \in (0,\constantAell)$.
  Let $\Phi \in C^{\infty}(X,E)$ be supported in $B(p,r)$.

  By the inequality~\eqref{inequality: flat elliptic estimate} and Lemma \ref{lemma: comparison of two freezing operators}, we have
  \begin{equation}
    \begin{aligned}
    \flatmesoscopicdifferenceseminorm{\cutoffupper}{\flatmodel{\Phi}} \le{}& \constantflat \Ltwonorm[L^2(\flattorus)]{\flatmesoscopicWilsonDiracoperator \flatmodel{\Phi}} \\
    \le{}& \constantflat \Ltwonorm[L^2(\flattorus)]{\flatmodel{\big(\mesoscopicWilsonDiracoperator \Phi\big)}} + \constantflat \constantWD \left( r \flatmesoscopicdifferenceseminorm{\cutoffupper}{\flatmodel{\Phi}} + \Ltwonorm[L^2(X)]{\Phi} \right).
    \end{aligned}
  \end{equation}
  Hence, we have
  \begin{equation}
    \left(1 - \constantflat \constantWD r\right) \flatmesoscopicdifferenceseminorm{\cutoffupper}{\flatmodel{\Phi}} \le \constantflat \Ltwonorm[L^2(\flattorus)]{\flatmodel{\big(\mesoscopicWilsonDiracoperator \Phi\big)}} + \constantflat \constantWD  \Ltwonorm[L^2(X)]{\Phi}.
  \end{equation}
  Since $r \le 1 / (2 \constantflat \constantWD)$, we have
  \begin{equation}
    \flatmesoscopicdifferenceseminorm{\cutoffupper}{\flatmodel{\Phi}} \le 2\constantflat \Ltwonorm[L^2(\flattorus)]{\flatmodel{\big(\mesoscopicWilsonDiracoperator \Phi\big)}} + 2\constantflat \constantWD  \Ltwonorm[L^2(X)]{\Phi}.
  \end{equation}
  Since $\mesoscopicWilsonDiracoperator \Phi$ is supported in $B(p,2r)$, by Lemma~\ref{lemma: Ltwo norm comparison}, we have
  \begin{equation}
    \flatmesoscopicdifferenceseminorm{\cutoffupper}{\flatmodel{\Phi}} \le 2 \constantflat \constantLtwonorm \Ltwonorm[L^2(X)]{\mesoscopicWilsonDiracoperator \Phi} + 2\constantflat \constantWD \Ltwonorm[L^2(X)]{\Phi}.
  \end{equation}
  By Lemma~\ref{lemma: comparison of difference seminorms}, we have
  \begin{equation}
    \begin{aligned}
    \mesoscopicdifferenceseminorm{\Phi} \le{}& \constantnabla \left( \flatmesoscopicdifferenceseminorm{\cutoffupper}{\flatmodel{\Phi}} + \Ltwonorm[L^2(X)]{\Phi} \right) \\
    \le{}& 2 \constantnabla \constantflat \constantLtwonorm \Ltwonorm[L^2(X)]{\mesoscopicWilsonDiracoperator \Phi} + \constantnabla(2\constantflat \constantWD + 1) \Ltwonorm[L^2(X)]{\Phi}.
    \end{aligned}
  \end{equation}
  This completes the proof.
\end{proof}

\subsection{The commutator estimate}

For $f \in C^{\infty}(X,\R)$, we also denote by $f$ the multiplication operator on sections of $E$.

\begin{proposition}\label{proposition: commutator estimate}
  There exists a constant $\constantcommutator > 0$ such that
  \begin{equation}
    \Ltwonorm{\commutator{\mesoscopicWilsonDiracoperator}{f}\Phi} \le \constantcommutator \norm{df}_{L^{\infty}}\Ltwonorm{\Phi}
  \end{equation}
  for $a \in (0,\injectivityradius(X,g)/10)$, $f \in C^{\infty}(X,\R)$, and $\Phi \in C^{\infty}(X,E)$.
\end{proposition}
\begin{proof}
  Fix $a \in (0,\injectivityradius(X,g)/10)$.
  Let $f \in C^{\infty}(X,\R)$ and $\Phi \in C^{\infty}(X,E)$.
  Fix $x \in X$.
  The integral representations of $\mesoscopicDiracoperator$ and $\mesoscopicWilsonterm$ give
  \begin{equation}
    (\commutator{\mesoscopicDiracoperator}{f}\Phi)(x) = \int_X (f(y)-f(x)) \kernelDirac(x,y) \Phi(y) \dvol_g(y)
  \end{equation}
  and
  \begin{equation}
    (\commutator{\mesoscopicWilsonterm}{f}\Phi)(x) = \int_X \frac{f(x)-f(y)}{d_g(x,y)} \kernelWilson(x,y) \paralleltransport{y}{x}(\Phi(y)) \dvol_g(y).
  \end{equation}
  We have
  \begin{equation}
    \commutator{\mesoscopicWilsonDiracoperator}{f} = \commutator{\mesoscopicDiracoperator}{f} + \commutator{\mesoscopicWilsonterm}{f} \circ \gamma.
  \end{equation}
  Thus, we have
  \begin{equation}
    \norm*{(\commutator{\mesoscopicWilsonDiracoperator}{f}\Phi)(x)} \le \norm*{df}_{L^{\infty}} \int_X \big( d_g(x,y) \norm{\kernelDirac(x,y)} + \kernelWilson(x,y) \big) \norm{\Phi(y)} \,\dvol_g(y).
  \end{equation}
  By changing variables $y = \exp_x(v)$ and $v=au$, we have
  \begin{equation}
    \begin{aligned}
    \int_X \big( d_g(x,y) \norm{\kernelDirac(x,y)} + \kernelWilson(x,y) \big) \,\dvol_g(y) ={}& \int_{TX_x} \left( \frac{\cutoffDirac(\abs{v}/a)}{C_D} + \frac{\cutoffWilson(\abs{v}/a)}{C_W} \right) \,\frac{dv}{a^n} \\
    ={}& \int_{\R^n} \left( \frac{\cutoffDirac(\abs{u})}{C_D} + \frac{\cutoffWilson(\abs{u})}{C_W} \right) \,du \\
    ={}& n + \frac{1}{C_W}\int_{\R^n} \cutoffWilson(\abs{u}) \,du =: \constantcommutator.
    \end{aligned}
  \end{equation}
  In particular, this integral is independent of $x$ and $a$.
  Symmetry gives
  \begin{equation}
    \int_X \big( d_g(x,y) \norm{\kernelDirac(x,y)} + \kernelWilson(x,y) \big) \,\dvol_g(x) = \constantcommutator
  \end{equation}
  for every $y \in X$.
  The Cauchy--Schwarz inequality yields
  \begin{equation}
    \begin{aligned}
      \Ltwonorm{\commutator{\mesoscopicWilsonDiracoperator}{f}\Phi}^2 ={}& \int_X \norm*{(\commutator{\mesoscopicWilsonDiracoperator}{f}\Phi)(x)}^2 \,\dvol_g(x) \\
      \le{}& \constantcommutator \norm*{df}_{L^{\infty}}^2 \int_X \int_X \big( d_g(x,y) \norm{\kernelDirac(x,y)} + \kernelWilson(x,y) \big) \norm{\Phi(y)}^2 \,\dvol_g(y) \,\dvol_g(x) \\
      ={}& \constantcommutator^2 \norm*{df}_{L^{\infty}}^2 \Ltwonorm{\Phi}^2.
    \end{aligned}
  \end{equation}
  This completes the proof.
\end{proof}

\subsection{The global estimate}

We now combine the local elliptic estimate with the commutator estimate.

\begin{proof}[Proof of Theorem~\ref{theorem: global elliptic estimate}]
  Fix the constants $r$, $\constantAell$, and $\constantlocal$ in Proposition~\ref{proposition: local elliptic estimate}.
  Choose finitely many points $p_1, \dots, p_N \in X$ such that the balls $B(p_i, r/2)$ cover $X$, and choose a partition of unity $f_1, \dots, f_N \colon X \to [0,1]$ with $\support f_j \subset B(p_j,r)$ and
  \begin{equation}
    \sum_{j=1}^N f_j = 1.
  \end{equation}
  Fix $a \in (0,\constantAell)$.
  We note that $N$ and $f_j$ are fixed independently of $a$.
  Let $\Phi \in C^{\infty}(X,E)$.

  Since each $f_j \Phi$ is supported in $B(p_j,r)$, by Proposition~\ref{proposition: local elliptic estimate}, we have
  \begin{equation}
    \begin{aligned}
    \mesoscopicdifferenceseminorm{f_j \Phi} \le{}& \constantlocal \left( \Ltwonorm{\mesoscopicWilsonDiracoperator(f_j \Phi)} + \Ltwonorm{f_j\Phi} \right) \\
    \le{}& \constantlocal \left( \Ltwonorm{\mesoscopicWilsonDiracoperator(f_j \Phi)} + \Ltwonorm{\Phi} \right).
    \end{aligned}
  \end{equation}
  Since
  \begin{equation}
    \mesoscopicWilsonDiracoperator(f_j \Phi) = f_j \big(\mesoscopicWilsonDiracoperator \Phi\big) + \commutator{\mesoscopicWilsonDiracoperator}{f_j} \Phi,
  \end{equation}
  by Proposition~\ref{proposition: commutator estimate}, we have
  \begin{equation}
    \begin{aligned}
    \Ltwonorm{\mesoscopicWilsonDiracoperator(f_j \Phi)} \le{}& \Ltwonorm{f_j \mesoscopicWilsonDiracoperator \Phi} + \constantcommutator \norm{df_j}_{L^{\infty}} \Ltwonorm{\Phi} \\
      \le{}& \Ltwonorm{\mesoscopicWilsonDiracoperator \Phi} + \constantcommutator \norm{df_j}_{L^{\infty}} \Ltwonorm{\Phi}.
    \end{aligned}
  \end{equation}
  Thus, we have
  \begin{equation}
    \begin{aligned}
    \mesoscopicdifferenceseminorm{\Phi} \le{}& \sum_{j=1}^N \mesoscopicdifferenceseminorm{f_j \Phi} \\
    \le{}& \sum_{j=1}^N \constantlocal \left( \Ltwonorm{\mesoscopicWilsonDiracoperator(f_j \Phi)} + \Ltwonorm{\Phi} \right) \\
    \le{}& \sum_{j=1}^N \constantlocal \left( \big( \Ltwonorm{\mesoscopicWilsonDiracoperator \Phi} + \constantcommutator \norm{df_j}_{L^{\infty}} \Ltwonorm{\Phi} \big) + \Ltwonorm{\Phi} \right) \\
    \le{}& \constantlocal N \Ltwonorm{\mesoscopicWilsonDiracoperator \Phi} +  \constantlocal N \big(1 + \constantcommutator (\max \norm{df_j}_{L^{\infty}})\big) \Ltwonorm{\Phi}.
    \end{aligned}
  \end{equation}
  Since both sides are continuous on $L^2(X,E)$ for fixed $a$, the estimate extends to all $\Phi \in L^2(X,E)$.
  This completes the proof.
\end{proof}

\section{The mesoscopic Rellich lemma}\label{section: the mesoscopic rellich lemma}

In this section, we prove a Rellich-type compactness theorem for sequences of sections whose covariant mesoscopic difference seminorms are bounded as the mesoscopic scale tends to zero.
We give two proofs of the compactness theorem.
The first proof is based on estimates at short distances, whereas the second proof is based on estimates at high frequencies.
Both proceed by finite-dimensional approximation, using local averages in the first case and Fourier truncation in the second.

\subsection{The lemma}

We keep the geometric data $(X,g,E,\nabla)$ and the cut-off function $\cutoffnabla$ fixed as in the preceding sections.
We also fix a non-zero smooth cut-off function $\cutofflower \colon [0,\infty) \to [0,\infty)$ whose support is contained in an open interval on which $\cutoffnabla$ is positive.
Then, there exists $C_{\cutofflower} > 0$ such that
\begin{equation}
  \cutofflower(t) + \abs{(\cutofflower)'(t)} \le C_{\cutofflower} \cutoffnabla(t)
\end{equation}
for all $t \ge 0$.

\begin{theorem}[The mesoscopic Rellich lemma]\label{theorem: mesoscopic Rellich lemma}
  Let $\{a_j\}_{j=1}^{\infty} \subset (0,\injectivityradius(X,g)/10)$ and $\{\Phi_j\}_{j=1}^{\infty} \subset L^2(X,E)$.
  We assume
  \begin{equation}
    \lim_{j \to \infty} a_j = 0
  \end{equation}
  and
  \begin{equation}
    \sup_{j \ge 1} \bigl( \Ltwonorm{\Phi_j} + \mesoscopicdifferenceseminorm[a_j]{\Phi_j} \bigr) < \infty.
  \end{equation}
  Then, there exists a subsequence $\{\Phi_{j_k}\}_{k=1}^{\infty}$ and $\Phi \in L^2(X,E)$ such that $\Phi_{j_k}$ converges to $\Phi$ strongly in $L^2(X,E)$.
\end{theorem}

\begin{corollary}\label{corollary: mesoscopic elliptic and Rellich}
  Let $\{a_j\}_{j=1}^{\infty} \subset (0,\injectivityradius(X,g)/10)$ and $\{\Phi_j\}_{j=1}^{\infty} \subset L^2(X,E)$.
  We assume
  \begin{equation}
    \lim_{j \to \infty} a_j = 0
  \end{equation}
  and
  \begin{equation}
    \sup_{j \ge 1} \bigl( \Ltwonorm{\Phi_j} + \Ltwonorm{\mesoscopicWilsonDiracoperator[a_j]\Phi_j} \bigr) < \infty.
  \end{equation}
  Then, there exists a subsequence $\{\Phi_{j_k}\}_{k=1}^{\infty}$ and $\Phi \in L^2_1(X,E)$ such that $\Phi_{j_k}$ converges to $\Phi$ strongly and $\mesoscopicWilsonDiracoperator[a_{j_k}] \Phi_{j_k}$ converges to $D\Phi$ weakly in $L^2(X,E)$.
\end{corollary}
\begin{proof}
  Since $a_j \to 0$, there exists $j_0$ such that $a_j < \constantAell$ for $j \ge j_0$.
  Then, by Theorem~\ref{theorem: global elliptic estimate}, we have
  \begin{equation}
    \sup_{j \ge j_0} \mesoscopicdifferenceseminorm[a_j]{\Phi_j} < \infty.
  \end{equation}
  By Theorem~\ref{theorem: mesoscopic Rellich lemma}, after passing to a subsequence, we have
  \begin{equation}
    \Phi_{j_k} \to \Phi
  \end{equation}
  strongly in $L^2(X,E)$ for some $\Phi \in L^2(X,E)$.
  Since $\{ \mesoscopicWilsonDiracoperator[a_{j_k}]\Phi_{j_k} \}_k$ is bounded in $L^2(X,E)$, after passing to a further subsequence, we have
  \begin{equation}
    \mesoscopicWilsonDiracoperator[a_{j_k}] \Phi_{j_k} \rightharpoonup \Psi
  \end{equation}
  weakly in $L^2(X,E)$ for some $\Psi \in L^2(X,E)$.
  Let $\Theta \in C^{\infty}(X,E)$.
  By Propositions~\ref{proposition: convergence of mesoscopic Dirac operators} and~\ref{proposition: convergence of mesoscopic Wilson terms}, we have
  \begin{equation}
    \mesoscopicWilsonDiracoperator[a_{j_k}]\Theta \to D\Theta
  \end{equation}
  strongly in $L^2(X,E)$.
  Since the mesoscopic Wilson--Dirac operators are self-adjoint, we have
  \begin{equation}
    \innerproduct*{\Psi}{\Theta}_{L^2} = \lim_{k \to \infty} \innerproduct*{\mesoscopicWilsonDiracoperator[a_{j_k}]\Phi_{j_k}}{\Theta}_{L^2} = \lim_{k \to \infty} \innerproduct*{\Phi_{j_k}}{\mesoscopicWilsonDiracoperator[a_{j_k}]\Theta}_{L^2} = \innerproduct*{\Phi}{D\Theta}_{L^2}.
  \end{equation}
  Thus, $D\Phi=\Psi$ weakly.
  By elliptic regularity for $D$, we have $\Phi \in L^2_1(X,E)$.
  Thus, $\Psi = D\Phi$ in $L^2(X,E)$.
  This completes the proof.
\end{proof}

\begin{corollary}\label{corollary: norm-resolvent convergence}
  The mesoscopic Wilson--Dirac operator $\mesoscopicWilsonDiracoperator$ converges to $D$ in the norm-resolvent sense as $a \to 0$.
  In other words, we have
  \begin{equation}
    \lim_{a \to 0} \operatornorm*{\bigl( \mesoscopicWilsonDiracoperator - z \bigr)^{-1} - (D-z)^{-1}} = 0
  \end{equation}
  for any $z \in \C \setminus \R$.
\end{corollary}
\begin{proof}
  Fix $z \in \C \setminus \R$.
  Suppose that there exist $\epsilon > 0$, $a_j \to 0$, and $\Psi_j \in L^2(X,E)$ with $\Ltwonorm{\Psi_j} = 1$ such that
  \begin{equation}
    \norm*{\bigl( \mesoscopicWilsonDiracoperator[a_j] - z \bigr)^{-1}\Psi_j - (D-z)^{-1} \Psi_j}_{L^2} \ge \epsilon
  \end{equation}
  for $j = 1,2,3,\dots$.
  Set $\Phi_j := \bigl( \mesoscopicWilsonDiracoperator[a_j] - z \bigr)^{-1}\Psi_j$.
  Then, both $\Phi_j$ and $\mesoscopicWilsonDiracoperator[a_j]\Phi_j = \Psi_j + z\Phi_j$ are bounded in $L^2(X,E)$.
  By Corollary~\ref{corollary: mesoscopic elliptic and Rellich}, after passing to a subsequence, $\Phi_j$ converges strongly to some $\Phi \in L^2_1(X,E)$, and $\Psi_j$ converges weakly to $(D-z)\Phi$.
  Since $(D-z)^{-1}$ is compact by the ordinary Rellich lemma, $(D-z)^{-1}\Psi_j$ converges strongly to $\Phi$.
  This contradicts the above lower bound.
\end{proof}

\begin{remark}\label{remark: norm-resolvent nonconvergence}
  The Wilson term has an essential role in this norm-resolvent convergence.
  Let $a \in (0,\injectivityradius(X,g)/10)$ and $z \in \C \setminus \R$.
  Since $\mesoscopicDiracoperator$ and $(D-z)^{-1}$ are compact, the operator
  \begin{equation}
    \bigl( \mesoscopicDiracoperator - z \bigr)^{-1} - (D-z)^{-1} + z^{-1}\id = z^{-1} \bigl(\mesoscopicDiracoperator - z\bigr)^{-1} \mesoscopicDiracoperator - (D-z)^{-1}
  \end{equation}
  is compact.
  Hence, as in the proof of Proposition~\ref{proposition: convergence of mesoscopic Wilson terms}, we have
  \begin{equation}
    \operatornorm*{\bigl(\mesoscopicDiracoperator - z\bigr)^{-1} - (D-z)^{-1}} \ge \frac{1}{\abs{z}}.
  \end{equation}
  Thus, $\mesoscopicDiracoperator$ does not converge to $D$ in the norm-resolvent sense.
  See \cite{MR4590014} and \cite{MR4741062} for the corresponding results on Euclidean lattices.
\end{remark}

Before proving Theorem~\ref{theorem: mesoscopic Rellich lemma}, we recall the corresponding proofs of the ordinary Rellich lemma.

\subsection{Two proofs of the ordinary Rellich lemma}

The ordinary Rellich lemma says that the inclusion map
\begin{equation}
  \iota \colon L^2_1(\flattorus) \to L^2(\flattorus)
\end{equation}
can be approximated in operator norm by finite-rank operators, and hence is compact.
We describe two ways to construct such finite-rank approximations:
one focuses on short distances, whereas the other focuses
on high frequencies.

\subsubsection{Approximation by block transformations}

For simplicity, take $\flattorus=\R^n/\Z^n$.
For each $N = 1,2,3,\dots$, we divide the torus $\flattorus$ into cubes $Q$ of side length $1/N$.
For $u \in L^2_1(\flattorus)$, we define $\iota_N(u)$ by
\begin{equation}
  \iota_N(u) := \sum_Q u_Q 1_Q,
\end{equation}
where
\begin{equation}
  u_Q := \frac{1}{\abs{Q}} \int_Q u(y) \,dy
\end{equation}
and $1_Q$ denotes the characteristic function of $Q$.
This averaging procedure is sometimes called a \emph{block transformation}.
The range of $\iota_N$ is finite dimensional.
By the Poincar\'e inequality, we have
\begin{equation}
  \int_Q \abs*{u(x) - u_Q}^2 \,dx \le \frac{1}{N^2} \int_Q \abs{\nabla u(x)}^2 \,dx.
\end{equation}
Summing over the cubes, we have
\begin{equation}
  \Ltwonorm{\iota(u) - \iota_N(u)} \le \frac{1}{N} \Ltwonorm{\nabla u}.
\end{equation}
Thus, $\iota$ is approximated by finite-rank operators $\iota_N$ as $N \to \infty$.

\subsubsection{Control of the high-frequency part}

For $u \in L^2_1(\flattorus)$, we write
\begin{equation}
  u = \sum_{\xi \in \duallattice{\Lattice}} e_{\xi} u_{\xi}
\end{equation}
in the normalised Fourier basis.
For each $R > 0$, we define
\begin{equation}
  \iota_R(u) := \sum_{\abs{\xi} \le R} e_{\xi} u_{\xi}.
\end{equation}
The range of $\iota_R$ is also finite dimensional.
We have
\begin{equation}
  \norm*{\iota(u) - \iota_R(u)}^2_{L^2} = \sum_{\abs{\xi} > R} \abs{u_{\xi}}^2 \le \frac{1}{R^2} \sum_{\xi} \abs{\xi}^2 \abs{u_{\xi}}^2 = \frac{1}{R^2}\Ltwonorm{\nabla u}^2.
\end{equation}
Thus, $\iota$ is approximated by finite-rank operators $\iota_R$ as $R \to \infty$.

\subsection{The first proof}

\subsubsection{A shell difference Poincar\'e inequality}

We first establish a variant of the Poincar\'e inequality.
Given $a \in (0,\injectivityradius(X,g)/10)$ and an open set $U \subset X$, we set
\begin{equation}
  \mesoscopicdifferenceseminormLocal{U}{\Phi}^2 := \int_U \int_U \norm*{\frac{\Phi(x) - \paralleltransport{y}{x}(\Phi(y))}{d_g(x,y)}}^2 \cutoffnabla(d_g(x,y)/a) \,\frac{\dvol_g(y)}{a^n} \,\dvol_g(x)
\end{equation}
for $\Phi \in L^2(U,E)$.
When $U = X$, this is the seminorm of Definition~\ref{definition: mesoscopic difference seminorm on manifolds}.
For a function $u$ on $U$ with values in a fixed finite-dimensional Hermitian vector space, we define $\mesoscopicdifferenceseminormLocal{U}{u}$ by the same formula with the parallel transport omitted.
By the support condition on $\cutoffnabla$, the right-hand side
involves only pairs of points $x,y$ with $r_-^{\nabla}a \le d_g(x,y) \le r_+^{\nabla}a$.
Thus, for each $x\in U$, the integration in $y$ is restricted
to a geodesic shell centred at $x$.
For an integrable function $u$ on a measurable set $U \subset X$ of positive volume, we write
\begin{equation}
  u_U := \frac{1}{\abs{U}} \int_U u \,\dvol_g.
\end{equation}

\begin{proposition}\label{proposition: shell poincare inequality}
  There exist $\radiusPoincare > 0$ and $\constantPoincare > 0$ such that
  \begin{equation}
    \int_{B(p,r)} \abs*{u(x) - u_{B(p,r)}}^2 \,\dvol_g(x) \le \constantPoincare r^2 \mesoscopicdifferenceseminormLocal{B(p,2r)}{u}^2
  \end{equation}
  for every $p \in X$, $r \in (0,\radiusPoincare)$, $a \in (0,r)$, and $u \in L^2(B(p,2r))$.
\end{proposition}

We first mollify $u$ by averaging over a geodesic shell at scale $a$.
The mesoscopic difference seminorm controls both the mollifying error and the gradient of the mollified function.
We then apply the ordinary Poincar\'e inequality to the mollified function.

\begin{proof}
  By uniform comparisons of the metric and volume element in normal coordinates, we can choose $\radiusPoincare > 0$ and $C_1 > 0$ such that the ordinary Poincar\'e inequality
  \begin{equation}
    \int_{B(p,r)} \abs*{f(x) - f_{B(p,r)}}^2 \,\dvol_g(x) \le C_1^2 r^2 \int_{B(p,r)} \abs{\nabla f(x)}^2 \,\dvol_g(x)
  \end{equation}
  holds for every $p \in X$, $r \in (0,\radiusPoincare)$ and $f \in L^2_1(B(p,r))$.

  Fix $p \in X$, $r \in (0,\radiusPoincare)$, and $a \in (0,r)$.
  Set $B_r := B(p,r)$ and $B_{2r} := B(p,2r)$.
  Let $u \in L^2(B_{2r})$.
  We will construct a smooth function $u_a$ on $B_r$ such that
  \begin{equation}\label{inequality: shell mollification error}
    \Ltwonorm[L^2(B_r)]{u-u_a}^2 \le C_2^2 a^2 \mesoscopicdifferenceseminormLocal{B_{2r}}{u}^2
  \end{equation}
  and
  \begin{equation}\label{inequality: shell gradient estimate}
    \Ltwonorm[L^2(B_r)]{\nabla u_a}^2 \le C_2^2 \mesoscopicdifferenceseminormLocal{B_{2r}}{u}^2,
  \end{equation}
  where $C_2 > 0$ is independent of $p$, $r$, $a$ and $u$.
  These estimates suffice to prove the inequality.
  Indeed, since $u_{B_r}$ is the best constant approximation to $u$ in $L^2(B_r)$, the ordinary Poincar\'e inequality gives
  \begin{equation}
    \begin{aligned}
      \Ltwonorm[L^2(B_r)]{u-u_{B_r}} \le{}& \Ltwonorm[L^2(B_r)]{u-(u_a)_{B_r}}\\
      \le{}& \Ltwonorm[L^2(B_r)]{u-u_a} + \Ltwonorm[L^2(B_r)]{u_a - (u_a)_{B_r}} \\
      \le{}& \Ltwonorm[L^2(B_r)]{u-u_a} + C_1 r \Ltwonorm[L^2(B_r)]{\nabla u_a} \\
      \le{}& C_2 a \mesoscopicdifferenceseminormLocal{B_{2r}}{u} + C_1 r C_2 \mesoscopicdifferenceseminormLocal{B_{2r}}{u} \\
      \le{}& (C_1+1)C_2 r \mesoscopicdifferenceseminormLocal{B_{2r}}{u}.
    \end{aligned}
  \end{equation}

  \begin{figure}[htbp]
    \centering
    \begin{tikzpicture}[
      x=0.012cm,y=-0.012cm,
      font=\small\sffamily,text=ink,
      every node/.style={inner sep=0pt,outer sep=0pt},
      fine dots/.style={dash pattern=on 0.7pt off 1.35pt,line cap=round},
      ball boundary/.style={draw=ballline,line width=0.55pt,fine dots},
      blue boundary/.style={draw=blueedge,line width=0.45pt},
      gold boundary/.style={draw=goldedge,line width=0.45pt},
      guide/.style={draw=blueedge,line width=0.4pt}
    ]
    \def\BallRadius{220}
    \def\ShellScale{220}
    \pgfmathsetmacro{\InnerBlue}{0.28*\ShellScale}
    \pgfmathsetmacro{\OuterBlue}{0.90*\ShellScale}
    \pgfmathsetmacro{\InnerGold}{0.46*\ShellScale}
    \pgfmathsetmacro{\OuterGold}{0.70*\ShellScale}
    \coordinate (p) at (560,650);
    \coordinate (x) at (766.8,650);
    \coordinate (y) at ($(x)+({0.58*\ShellScale*cos(25)},
                              {-0.58*\ShellScale*sin(25)})$);

    \fill[ballouter] (p) circle[radius={2*\BallRadius}];
    \fill[ballinner] (p) circle[radius=\BallRadius];
    \draw[ball boundary] (p) circle[radius={2*\BallRadius}];
    \draw[draw=ballline!65,line width=0.4pt] (p) -- ++(-\BallRadius,0);
    \node[anchor=base west] at (437,637) {$r$};
    \fill[ink] (p) circle[radius=4];
    \node[anchor=base west] at (543,682) {$p$};

    \fill[shellblue,even odd rule]
      (x) circle[radius=\OuterBlue] (x) circle[radius=\InnerBlue];
    \fill[shellgold,even odd rule]
      (x) circle[radius=\OuterGold] (x) circle[radius=\InnerGold];
    \draw[blue boundary] (x) circle[radius=\OuterBlue];
    \draw[blue boundary] (x) circle[radius=\InnerBlue];
    \draw[gold boundary] (x) circle[radius=\OuterGold];
    \draw[gold boundary] (x) circle[radius=\InnerGold];
    \draw[ball boundary] (p) circle[radius=\BallRadius];

    \node[anchor=base west] at (448,500) {$B_r=B(p,r)$};
    \node[anchor=base west] at (430,1020) {$B_{2r}=B(p,2r)$};
    \fill[ink] (x) circle[radius=3.7];
    \node[anchor=base west] at ($(x)+(-12,25)$) {$x$};
    \fill[ink] (y) circle[radius=3.7];
    \node[anchor=base west] at ($(y)+(-12,25)$) {$y$};

    \coordinate (outer-radius-end) at
      ($(x)+({\OuterBlue*cos(45)},{\OuterBlue*sin(45)})$);
    \draw[guide] (x) --
      node[midway,above=2pt,sloped,text=blueedge] {$ar_+^{\nabla}$}
      (outer-radius-end);
    \draw[guide] (x) --
      node[midway,above=2pt,text=blueedge]
        {$ar_-^{\nabla}$}
      ++(-\InnerBlue,0);
    \end{tikzpicture}
  \end{figure}

  We now construct $u_a$ by averaging over a geodesic shell.
  Recall the cut-off function $\cutofflower$ fixed at the beginning of this section.
  We define, for $x,y \in X$,
  \begin{equation}
    q_a(x,y) := \frac{1}{c_a(x)} \frac{\cutofflower(d_g(x,y)/a)}{a^n},
  \end{equation}
  where
  \begin{equation}
    c_a(x) := \int_X \frac{\cutofflower(d_g(x,y)/a)}{a^n} \,\dvol_g(y).
  \end{equation}
  The kernel $q_a$ is smooth.
  For $x \in B_r$, the support of $q_a(x,\cdot)$ is contained in $B_{2r}$.
  Hence, the normalisation of $q_a$ yields
  \begin{equation}
    \int_{B_{2r}} q_a(x,y) \,\dvol_g(y) = 1
  \end{equation}
  and
  \begin{equation}
    \int_{B_{2r}} \nabla_x q_a(x,y) \,\dvol_g(y) = 0.
  \end{equation}
  We have $C_c^{-1} \le c_a(x) \le C_c$ by the uniform bounds on $J_x$, and $\abs{\nabla c_a(x)} \le C_c/a$ since $\abs{\nabla_x d_g(x,y)} = 1$.
  Thus, we have
  \begin{equation}\label{inequality: mollified qa estimate}
    q_a(x,y) + a\abs{\nabla_x q_a(x,y)} \le \frac{C_3}{a^n}\cutoffnabla(d_g(x,y)/a)
  \end{equation}
  and
  \begin{equation}
    \int_X \abs{\nabla_x q_a(x,y)} \,\dvol_g(y) \le \frac{C_3}{a},
  \end{equation}
  where $C_3$ is independent of $x$ and $a$.

  The mollified function
  \begin{equation}
    u_a(x) := \int_{B_{2r}} q_a(x,y)u(y) \,\dvol_g(y)
  \end{equation}
  is smooth on $B_r$.
  Thus, we have
  \begin{equation}
      u(x)-u_a(x) = \int_{B_{2r}} q_a(x,y) (u(x)-u(y)) \,\dvol_g(y)
  \end{equation}
  and
  \begin{equation}
      \nabla u_a(x) = \int_{B_{2r}} \nabla_x q_a(x,y) (u(y)-u(x)) \,\dvol_g(y)
  \end{equation}
  for $x \in B_r$.
  These identities express both the mollifying error and the gradient in terms of differences on the shell.

  Note that $d_g(x,y)\le a$ whenever $\cutoffnabla(d_g(x,y)/a)\ne 0$.
  By Jensen's inequality and \eqref{inequality: mollified qa estimate}, we have
  \begin{equation}
    \begin{aligned}
      \abs{u(x)-u_a(x)}^2 \le{}& \int_{B_{2r}} q_a(x,y)\abs{u(x)-u(y)}^2 \,\dvol_g(y) \\
      \le{}& C_3 \int_{B_{2r}} \frac{\cutoffnabla(d_g(x,y)/a)}{a^n} \abs{u(x)-u(y)}^2 \,\dvol_g(y) \\
      \le{}&  C_3 a^2 \int_{B_{2r}} \frac{\cutoffnabla(d_g(x,y)/a)}{a^n}
        \frac{\abs{u(x)-u(y)}^2}{d_g(x,y)^2}\,\dvol_g(y).
    \end{aligned}
  \end{equation}
  Integrating over $x \in B_r$ proves \eqref{inequality: shell mollification error}.
  Similarly, the Cauchy--Schwarz inequality gives
  \begin{equation}
    \begin{aligned}
      \abs{\nabla u_a(x)}^2 \le{}& \left(\int_{B_{2r}} \abs{\nabla_x q_a(x,y)} \,\dvol_g(y)\right) \left(\int_{B_{2r}} \abs{\nabla_x q_a(x,y)}\abs{u(y)-u(x)}^2 \,\dvol_g(y)\right) \\
      \le{}& C_3^2 \int_{B_{2r}} \frac{\cutoffnabla(d_g(x,y)/a)}{a^n} \frac{\abs{u(y)-u(x)}^2}{d_g(x,y)^2} \,\dvol_g(y),
    \end{aligned}
  \end{equation}
  which proves \eqref{inequality: shell gradient estimate}.
  The proof is complete.
\end{proof}

We remark that the same inequality holds, with the same constant, for functions taking values in a finite-dimensional Hermitian vector space.

\subsubsection{Voronoi cells and finite-dimensional approximation}

Fix $\ell>0$.
A subset of $X$ is called $\ell$-separated if the distance between any two distinct points in it is at least $\ell$.
It is maximal if it is not properly contained in any other $\ell$-separated subset of $X$.
By compactness of $X$, we can choose a finite maximal $\ell$-separated set $\{p_1,\dots,p_N\}$.

For each $x\in X$, there is a unique index $\alpha \in \{1,\ldots,N\}$ such that
\begin{equation}
  \begin{cases}
    1 \le j < \alpha \Longrightarrow d_g(x,p_{\alpha}) < d_g(x,p_j) \\
    \alpha < j \le N \Longrightarrow d_g(x,p_{\alpha}) \le d_g(x,p_j)
  \end{cases}
\end{equation}
for any $j = 1,\dots,N$.
We denote this index by $\alpha(x)$ and define the Voronoi cells $\{\voronoicell_{\alpha}\}$ by
\begin{equation}
  \voronoicell_{\alpha} := \{ x\in X \mid \alpha(x)=\alpha\}
\end{equation}
for each $\alpha = 1,\dots,N$.
These cells form a measurable partition of $X$ and satisfy
\begin{equation}
  B(p_\alpha,\ell/2) \subset \voronoicell_\alpha \subset B(p_\alpha,\ell).
\end{equation}
The first inclusion follows from the separation condition, and the second follows from maximality.
The number of balls $B(p_{\alpha},2\ell)$ containing $x \in X$ is bounded by a constant independent of $x$ and $\ell$.
Indeed, the balls $B(p_{\alpha},\ell/2)$ are pairwise disjoint, and the volume doubling property gives a uniform bound on the number of centres $p_{\alpha}$ in any ball of radius $2\ell$.

\begin{figure}[htbp]
  \centering
  \begin{tikzpicture}[
    x=1cm,y=1cm,
    font=\fontsize{11}{13}\selectfont,text=ink,
    line cap=round,line join=round
  ]
  \begin{scope}[
      xshift=3.55cm,yshift=2.90cm,
      x=1.8cm,y=1.8cm
    ]
      \clip (-1.65,-1.35) rectangle (1.65,1.35);

      \fill[bluewash] (0,0) circle[radius=1];

      \fill[shellblue!65]
        (0.00314,-0.66981)
        -- (0.59892,-0.30050)
        -- (0.58533,0.30332)
        -- (-0.00314,0.66981)
        -- (-0.59892,0.30050)
        -- (-0.58533,-0.30332)
        -- cycle;

      \fill[bluewash] (0,0) circle[radius=0.5];

      \draw[ballline!55,line width=0.45pt]
        (0.58173,1.64943) -- (1.16723,1.31322)
        (1.73296,1.63053) -- (1.16723,1.31322)
        (1.16723,1.31322) -- (1.17685,0.67267)
        (1.17685,0.67267) -- (1.73568,0.34333)
        (0.58173,1.64943) -- (0.00125,1.31537)
        (0.58533,0.30332) -- (1.17685,0.67267)
        (0.58533,0.30332) -- (-0.00314,0.66981)
        (-0.00314,0.66981) -- (0.00125,1.31537)
        (-0.59803,1.66033) -- (-1.16335,1.34877)
        (-1.16335,1.34877) -- (-1.73238,1.65851)
        (-0.59803,1.66033) -- (0.00125,1.31537)
        (-0.00314,0.66981) -- (-0.59892,0.30050)
        (-1.16335,1.34877) -- (-1.18033,0.63988)
        (-1.18033,0.63988) -- (-0.59892,0.30050)
        (-1.73965,0.33482) -- (-1.18033,0.63988)
        (0.58533,0.30332) -- (0.59892,-0.30050)
        (1.73965,-0.33482) -- (1.18033,-0.63988)
        (0.59892,-0.30050) -- (1.18033,-0.63988)
        (-0.58173,-1.64943) -- (-0.00125,-1.31537)
        (-0.58173,-1.64943) -- (-1.16723,-1.31322)
        (-0.58533,-0.30332) -- (-1.17685,-0.67267)
        (-0.58533,-0.30332) -- (0.00314,-0.66981)
        (0.00314,-0.66981) -- (-0.00125,-1.31537)
        (-1.16723,-1.31322) -- (-1.17685,-0.67267)
        (-1.73296,-1.63053) -- (-1.16723,-1.31322)
        (-1.73568,-0.34333) -- (-1.17685,-0.67267)
        (-0.59892,0.30050) -- (-0.58533,-0.30332)
        (0.59892,-0.30050) -- (0.00314,-0.66981)
        (1.18033,-0.63988) -- (1.16335,-1.34877)
        (-0.00125,-1.31537) -- (0.59803,-1.66033)
        (0.59803,-1.66033) -- (1.16335,-1.34877)
        (1.73238,-1.65851) -- (1.16335,-1.34877);

      \draw[blueedge,line width=0.95pt]
        (0.00314,-0.66981)
        -- (0.59892,-0.30050)
        -- (0.58533,0.30332)
        -- (-0.00314,0.66981)
        -- (-0.59892,0.30050)
        -- (-0.58533,-0.30332)
        -- cycle;

      \draw[blueedge,densely dashed,line width=0.65pt]
        (0,0) circle[radius=1];
      \draw[blueedge!75,densely dotted,line width=0.65pt]
        (0,0) circle[radius=0.5];

      \fill[ink!80] (-0.59939,-0.96242)
        circle[radius=1.35pt];
      \fill[ink!80] (0.60163,-0.97058)
        circle[radius=1.35pt];
      \fill[ink!80] (-1.18372,-0.02663)
        circle[radius=1.35pt];
      \fill[ink!80] (0,0)
        circle[radius=1.35pt];
      \fill[ink!80] (1.18372,0.02663)
        circle[radius=1.35pt];
      \fill[ink!80] (-0.60163,0.97058)
        circle[radius=1.35pt];
      \fill[ink!80] (0.59939,0.96242)
        circle[radius=1.35pt];

      \draw[blueedge!75,line width=0.5pt]
        (0,0) -- (-0.5,0)
        node[midway,above=3pt,text=blueedge] {$\ell/2$};

      \draw[blueedge!75,line width=0.5pt]
        (0,0) -- (0.70711,-0.70711)
        node[midway,below right=1pt,text=blueedge] {$\ell$};

      \fill[ink] (0,0) circle[radius=1.7pt];
      \node[below left=3pt] at (0,0) {$p_\alpha$};
      \node[text=blueedge] at (0.03,0.37) {$\voronoicell_\alpha$};
    \end{scope}
  \end{tikzpicture}
\end{figure}

For $\ell < \injectivityradius(X,g)$, we define the orthogonal projection $\Pi_{\ell} \colon L^2(X,E) \to L^2(X,E)$ by
\begin{equation}
  (\Pi_{\ell} \Phi)(x) := \paralleltransport{p_{\alpha}}{x} \left( \frac{1}{\abs{\voronoicell_{\alpha}}} \int_{\voronoicell_{\alpha}} \paralleltransport{y}{p_{\alpha}} \bigl(\Phi(y)\bigr) \,\dvol_g(y) \right)
\end{equation}
for $x \in \voronoicell_{\alpha}$.
The image of $\Pi_{\ell}$ is the finite-dimensional space of sections of the form
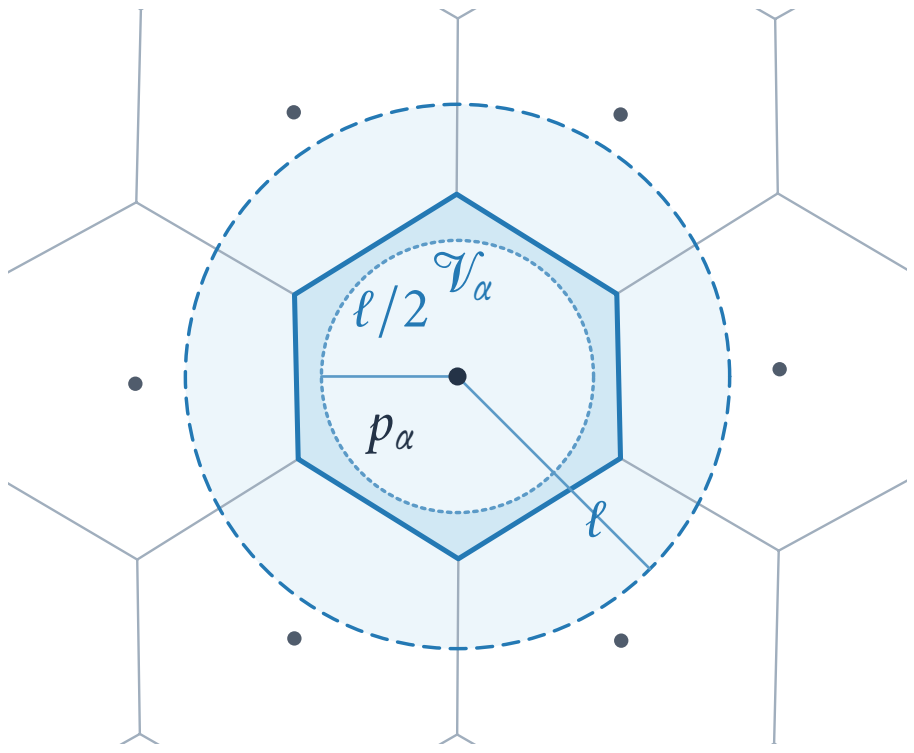
\begin{equation}
  x \mapsto \paralleltransport{p_{\alpha}}{x}(v_{\alpha})
\end{equation}
for $x \in \voronoicell_{\alpha}$ and $v_{\alpha} \in E_{p_{\alpha}}$.

\subsubsection{Completion of the proof}

\begin{proposition}\label{proposition: voronoi approximation error estimate}
  There exist $\ell_0 > 0$ and $C > 0$ such that
  \begin{equation}
    \Ltwonorm{\Phi - \Pi_{\ell}\Phi} \le C \ell \bigl( \mesoscopicdifferenceseminorm{\Phi} + \Ltwonorm{\Phi} \bigr)
  \end{equation}
  for $\ell \in (0,\ell_0)$, $a \in (0,\ell)$, and $\Phi \in L^2(X,E)$.
\end{proposition}
\begin{proof}
  Set $\ell_0 := \min \{\radiusPoincare, r_0\}$, where $r_0$ is the constant fixed at the end of Section~\ref{subsubsection: the flat model}.
  Let $\ell \in (0,\ell_0)$, and choose a finite maximal $\ell$-separated set $\{p_{\alpha}\}$.
  Take $a \in (0,\ell)$ and $\Phi \in L^2(X,E)$.
  On $B(p_{\alpha},2\ell)$, we set
  \begin{equation}
    \Phi_{\alpha}(y) := \paralleltransport{y}{p_{\alpha}} \bigl(\Phi(y)\bigr) \in E_{p_{\alpha}}.
  \end{equation}
  Then, $(\Pi_{\ell} \Phi)(x) = \paralleltransport{p_{\alpha}}{x} \bigl((\Phi_{\alpha})_{\voronoicell_{\alpha}}\bigr)$ for $x \in \voronoicell_{\alpha}$.
  By Proposition~\ref{proposition: shell poincare inequality}, we have
  \begin{equation}
    \begin{aligned}
    \int_{\voronoicell_{\alpha}} \norm{\Phi - \Pi_{\ell}\Phi}^2 \,\dvol_g ={}& \int_{\voronoicell_{\alpha}} \norm{\Phi_{\alpha} - (\Phi_{\alpha})_{\voronoicell_{\alpha}}}^2 \,\dvol_g \\
    \le{}& \int_{\voronoicell_{\alpha}} \norm{\Phi_{\alpha} - (\Phi_{\alpha})_{B(p_{\alpha},\ell)}}^2 \,\dvol_g \\
    \le{}& \constantPoincare \ell^2 \mesoscopicdifferenceseminormLocal{B(p_{\alpha},2\ell)}{\Phi_{\alpha}}^2.
    \end{aligned}
  \end{equation}

  By the holonomy estimate \eqref{inequality: holonomy bound} of Lemma~\ref{lemma: normal nbd comparison}, we have
  \begin{equation}
    \norm{\paralleltransport{x}{p_{\alpha}} \circ \paralleltransport{y}{x} - \paralleltransport{y}{p_{\alpha}}} \le \constantgeometric \ell d_g(x,y)
  \end{equation}
  for $x,y \in B(p_{\alpha},2\ell)$.
  Hence, we have
  \begin{equation}
    \norm{\Phi_{\alpha}(x) - \Phi_{\alpha}(y)} \le \norm{\Phi(x) - \paralleltransport{y}{x}(\Phi(y))} + \constantgeometric \ell d_g(x,y) \norm{\Phi(y)}.
  \end{equation}
  Thus, by the definition of $\mesoscopicdifferenceseminormLocal{B(p_{\alpha},2\ell)}{\Phi_{\alpha}}$, we have
  \begin{multline}
    \mesoscopicdifferenceseminormLocal{B(p_{\alpha},2\ell)}{\Phi_{\alpha}}^2 \\
    \le 2\mesoscopicdifferenceseminormLocal{B(p_{\alpha},2\ell)}{\Phi}^2+ 2\constantgeometric^2 \ell^2 \int_{B(p_{\alpha},2\ell)} \left(\int_{B(p_{\alpha},2\ell)} \cutoffnabla(d_g(x,y)/a) \,\frac{\dvol_g(x)}{a^n}\right) \norm{\Phi(y)}^2 \,\dvol_g(y).
  \end{multline}
  Since $d_g(x,y)<a$ on the support of $\cutoffnabla$, we have
  \begin{equation}
    \int_{B(p_{\alpha},2\ell)} \cutoffnabla(d_g(x,y)/a) \,\frac{\dvol_g(x)}{a^n} \le \norm{\cutoffnabla}_{L^{\infty}} \frac{\abs{B(y,a)}}{a^n} \le C_1,
  \end{equation}
  where $C_1$ is independent of $\alpha$, $\ell$, $a$, and $y$.
  Thus, we have
  \begin{equation}
    \mesoscopicdifferenceseminormLocal{B(p_{\alpha},2\ell)}{\Phi_{\alpha}}^2 \le2 \mesoscopicdifferenceseminormLocal{B(p_{\alpha},2\ell)}{\Phi}^2 +2 \constantgeometric^2 C_1 \ell^2 \Ltwonorm[L^2(B(p_{\alpha},2\ell))]{\Phi}^2.
  \end{equation}

  The balls $\{B(p_{\alpha},2\ell)\}$ have uniformly bounded overlap.
  In particular, each pair $(x,y) \in X \times X$ belongs to $B(p_{\alpha},2\ell) \times B(p_{\alpha},2\ell)$ for at most a uniformly bounded number $N_0$ of indices $\alpha$.
  Thus, we have
  \begin{equation}
    \begin{aligned}
    \Ltwonorm{\Phi - \Pi_{\ell}\Phi}^2 \le{}& \constantPoincare \ell^2 \sum_{\alpha} \mesoscopicdifferenceseminormLocal{B(p_{\alpha},2\ell)}{\Phi_{\alpha}}^2 \\
    \le{}& 2\constantPoincare \ell^2 N_0 \left( \mesoscopicdifferenceseminorm{\Phi}^2 + \constantgeometric^2 C_1 \ell^2 \Ltwonorm{\Phi}^2 \right) \\
    \le{}& 2\constantPoincare \ell^2 N_0 \left( \mesoscopicdifferenceseminorm{\Phi}^2 + \constantgeometric^2 C_1 \ell_0^2 \Ltwonorm{\Phi}^2 \right).
    \end{aligned}
  \end{equation}
  This completes the proof.
\end{proof}

\begin{proof}[First proof of Theorem~\ref{theorem: mesoscopic Rellich lemma}]
  Set
  \begin{equation}
    M := \sup_{j \ge 1} \left( \mesoscopicdifferenceseminorm[a_j]{\Phi_j} + \Ltwonorm{\Phi_j} \right) < \infty.
  \end{equation}
  Let $\epsilon > 0$.
  Choose $\ell \in (0,\ell_0)$ such that $C \ell M < \epsilon / 2$, where $C$ is the constant in Proposition~\ref{proposition: voronoi approximation error estimate}.
  Since $a_j \to 0$, there exists $j_0$ such that $a_j < \ell$ for every $j \ge j_0$.
  Then, by Proposition~\ref{proposition: voronoi approximation error estimate}, we have
  \begin{equation}
    \Ltwonorm{\Phi_j - \Pi_{\ell} \Phi_j} \le C \ell \bigl( \mesoscopicdifferenceseminorm[a_j]{\Phi_j} + \Ltwonorm{\Phi_j} \bigr) < \frac{\epsilon}{2}
  \end{equation}
  for every $j \ge j_0$.
  Hence, for $j \ge j_0$, each $\Phi_j$ lies within $\epsilon/2$ of the bounded set $\{\Pi_\ell \Phi_j\}_{j \ge j_0}$ in the finite-dimensional range of $\Pi_\ell$.
  Since $\epsilon > 0$ is arbitrary, $\{\Phi_j\}$ is totally bounded in $L^2(X,E)$.
  Hence, it has a strongly convergent subsequence.
  This completes the proof.
\end{proof}

\subsection{The second proof}

\subsubsection{The flat torus}

We first consider the flat model $\flattorus = V/\Lattice$ with $\flatmodel{E} = \flattorus \times S$.
For $R > 0$, we define the orthogonal projection $\iota_R$ on $L^2(\flattorus,\flatmodel{E})$ by the high-frequency truncation
\begin{equation}
  \iota_R(\Psi) := \sum_{\substack{\xi \in \duallattice{\Lattice} \\ \abs{\xi} \le R}} e_{\xi} \Psi_{\xi}
\end{equation}
for $\Psi = \sum e_{\xi} \Psi_{\xi} \in L^2(\flattorus,\flatmodel{E})$.
The range of $\iota_R$ is finite dimensional.
When $0 < a < 1/R$, we have
\begin{equation}
  \begin{aligned}
  \Ltwonorm{\Psi - \iota_R(\Psi)}^2 ={}& \sum_{\substack{\xi \in \duallattice{\Lattice} \\ \abs{\xi} > R}} \norm{\Psi_{\xi}}_S^2 \\
  \le{}& \sum_{\substack{\xi \in \duallattice{\Lattice} \\ \abs{\xi} > R}} \frac{\min\{\abs{\xi}^2, 1/a^2\}}{R^2} \norm{\Psi_{\xi}}_S^2 \\
  \le{}& \frac{1}{R^2} \sum_{\xi \in \duallattice{\Lattice}} \min \left\{ \abs{\xi}^2, 1/a^2 \right\} \norm{\Psi_{\xi}}_S^2 = \frac{1}{R^2} \mesoscopicseminorm{\Psi}^2.
  \end{aligned}
\end{equation}

Suppose that $a_j \to 0$ and that $\{\Psi_j\} \subset L^2(\flattorus,\flatmodel{E})$ satisfies
\begin{equation}
  M := \sup_{j \ge 1} \bigl( \mesoscopicseminorm[a_j]{\Psi_j} + \Ltwonorm{\Psi_j} \bigr) < \infty.
\end{equation}
Given $\epsilon > 0$, we choose $R > 0$ such that $M/R < \epsilon/2$.
Since $a_j \to 0$, there exists $j_0$ such that $a_j < 1/R$ for every $j \ge j_0$.
Hence, we have
\begin{equation}
  \Ltwonorm{\Psi_j - \iota_R(\Psi_j)} \le \frac{M}{R} < \frac{\epsilon}{2}.
\end{equation}
Hence, for $j \ge j_0$, each $\Psi_j$ lies within $\epsilon/2$ of the bounded set $\{\iota_R \Psi_j\}_{j \ge j_0}$ in the finite-dimensional range of $\iota_R$.
Since $\epsilon > 0$ is arbitrary, $\{\Psi_j\}$ is totally bounded in $L^2(\flattorus,\flatmodel{E})$.
Hence, it has a strongly convergent subsequence.
This completes the proof in the flat case.

\subsubsection{Localisation on the manifold}

To generalise the preceding argument to a general Riemannian manifold, we first need the following lemma.
Compare Lemma~\ref{lemma: comparison of difference seminorms}.
\begin{lemma}\label{lemma: another comparison of difference seminorms}
  There exist $r \in (0,r_0)$ and $\constantcomparison > 0$ such that
  \begin{equation}
    \mesoscopicseminorm{\flatmodel{\Phi}} \le \constantcomparison \bigl(\mesoscopicdifferenceseminorm{\Phi} + \Ltwonorm{\Phi} \bigr)
  \end{equation}
  for every $p \in X$, $a \in (0,r)$, and $\Phi \in L^2(X,E)$ supported in $B(p,r)$.
\end{lemma}
\begin{proof}
  By \eqref{inequality: uniform distance ratio comparison}, we can choose $r \in (0,r_0)$ sufficiently small that
  \begin{equation}\label{inequality: cutofflower and cutoffnabla comparison}
    \cutofflower(\abs{z-w}/a) \le C_1 \cutoffnabla \bigl( d_g(\exp_p(z),\exp_p(w))/a \bigr)
  \end{equation}
  for $z,w \in B(0,3r)$.

  Fix $p \in X$ and $a \in (0,r)$.
  Let $\Phi \in L^2(X,E)$ be supported in $B(p,r)$.
  By \eqref{inequality: holonomy bound}, we have
  \begin{equation}
    \norm{\flatmodel{\Phi}(z) - \flatmodel{\Phi}(w)} \le \norm{\Phi(x)-\paralleltransport{y}{x}(\Phi(y))} + \constantgeometric r d_g(x,y) \norm{\Phi(y)}
  \end{equation}
  for $z,w \in B(0,3r)$, $x := \exp_p(z)$, and $y := \exp_p(w)$.
  Hence, by \eqref{inequality: uniform distance comparison} and $\norm{\Phi(y)} = \norm{\flatmodel{\Phi}(w)}$, we have
  \begin{equation}
    \frac{\norm{\flatmodel{\Phi}(z) - \flatmodel{\Phi}(w)}^2}{\abs{z-w}^2} \le 2 \constantgeometric^2 \frac{\norm{\Phi(x)-\paralleltransport{y}{x}(\Phi(y))}^2}{d_g(x,y)^2} + 2 \constantgeometric^4 r^2 \norm{\flatmodel{\Phi}(w)}^2.
  \end{equation}
  Since $\flatmodel{\Phi}$ is supported in $B(0,r)$ and $a < r$, every non-zero contribution to its flat difference seminorm $\flatmesoscopicdifferenceseminorm{\cutofflower}{\flatmodel{\Phi}}$ comes from two points in $B(0,2r)$.
  Thus, we have
  \begin{equation}
    \begin{aligned}
      \flatmesoscopicdifferenceseminorm{\cutofflower}{\flatmodel{\Phi}}^2 ={}& \int_{B(0,2r)} \int_{B(0,2r)}  \frac{\norm{\flatmodel{\Phi}(z)-\flatmodel{\Phi}(w)}^2}{\abs{z-w}^2} \cutofflower\bigl(\abs{z-w}/a\bigr) \,\frac{dz}{a^n} \,dw \\
      \le{}& \begin{multlined}[t]
        2 \constantgeometric^2  \int_{B(0,2r)} \int_{B(0,2r)} \frac{\norm{\Phi(x)-\paralleltransport{y}{x}(\Phi(y))}^2}{d_g(x,y)^2} \cutofflower\bigl(\abs{z-w}/a\bigr) \,\frac{dz}{a^n} \,dw \\
        + 2 \constantgeometric^4 r^2 \int_{B(0,2r)} \int_{B(0,2r)} \norm{\flatmodel{\Phi}(w)}^2 \cutofflower\bigl(\abs{z-w}/a\bigr) \,\frac{dz}{a^n} \,dw.
      \end{multlined}
    \end{aligned}
  \end{equation}
   By \eqref{inequality: Jpw bound} and \eqref{inequality: cutofflower and cutoffnabla comparison}, we have
  \begin{equation}
    \int_{B(0,2r)} \int_{B(0,2r)} \frac{\norm{\Phi(x)-\paralleltransport{y}{x}(\Phi(y))}^2}{d_g(x,y)^2} \cutofflower\bigl(\abs{z-w}/a\bigr) \,\frac{dz}{a^n} \,dw \le C_1 \constantgeometric^2 \mesoscopicdifferenceseminorm{\Phi}^2.
  \end{equation}
  By \eqref{inequality: comparison of L2 norms}, we also have
  \begin{equation}
    \begin{aligned}
    \int_{B(0,2r)} \int_{B(0,2r)} \norm{\flatmodel{\Phi}(w)}^2 \cutofflower\bigl(\abs{z-w}/a\bigr) \,\frac{dz}{a^n} \,dw \le{}& \left(\int_{\R^n} \cutofflower(\abs{v}) \,dv \right) \Ltwonorm[L^2(\flattorus)]{\flatmodel{\Phi}}^2 \\
    \le{}& \left(\int_{\R^n} \cutofflower(\abs{v}) \,dv \right) \cdot \constantLtwonorm^2 \Ltwonorm[L^2(X)]{\Phi}^2 \\
    =:{}& M_{\cutofflower} \constantLtwonorm^2 \Ltwonorm[L^2(X)]{\Phi}^2.
    \end{aligned}
  \end{equation}
  Thus, we have
  \begin{equation}
    \flatmesoscopicdifferenceseminorm{\cutofflower}{\flatmodel{\Phi}}^2 \le 2\constantgeometric^4 C_1 \mesoscopicdifferenceseminorm{\Phi}^2 + 2 \constantgeometric^4 \constantLtwonorm^2 M_{\cutofflower} r^2 \Ltwonorm{\Phi}^2.
  \end{equation}

  By Proposition~\ref{proposition: equivalence of seminorms} applied with
  $\cutofflower$ in place of $\cutoffnabla$, there exists $\nu_{\cutofflower} > 0$ such that
  \begin{equation}
    \nu_{\cutofflower} \mesoscopicseminorm{\flatmodel{\Phi}} \le \flatmesoscopicdifferenceseminorm{\cutofflower}{\flatmodel{\Phi}}.
  \end{equation}
  Thus, we have
  \begin{equation}
    \mesoscopicseminorm{\flatmodel{\Phi}}^2 \le \frac{2\constantgeometric^4 C_1}{\nu_{\cutofflower}^2} \mesoscopicdifferenceseminorm{\Phi}^2 + \frac{2 \constantgeometric^4 \constantLtwonorm^2 M_{\cutofflower} r^2}{\nu_{\cutofflower}^2} \Ltwonorm{\Phi}^2.
  \end{equation}
  This completes the proof.
\end{proof}

\begin{lemma}\label{lemma: multiplication and difference seminorm}
  There exists a constant $\constantmult > 0$ such that
  \begin{equation}
    \mesoscopicdifferenceseminorm{f\Phi} \le \norm{f}_{L^{\infty}} \mesoscopicdifferenceseminorm{\Phi} + \constantmult \norm{df}_{L^{\infty}} \Ltwonorm{\Phi}
  \end{equation}
  for $a \in (0,r_0)$, $f \in C^{\infty}(X,\R)$, and $\Phi \in L^2(X,E)$.
\end{lemma}
\begin{proof}
  Since
  \begin{equation}
    f(x)\Phi(x) - \paralleltransport{y}{x}(f(y)\Phi(y)) = f(x)\bigl(\Phi(x) - \paralleltransport{y}{x}(\Phi(y))\bigr) + (f(x)-f(y)) \paralleltransport{y}{x}(\Phi(y)),
  \end{equation}
  we have
  \begin{equation}
    \begin{aligned}
    \mesoscopicdifferenceseminorm{f\Phi} ={}& \left(\int_X \int_X \norm*{\frac{f(x)\Phi(x) - \paralleltransport{y}{x}(f(y)\Phi(y))}{d_g(x,y)}}^2 \cutoffnabla(d_g(x,y)/a) \,\frac{\dvol_g(y)}{a^n} \,\dvol_g(x)\right)^{1/2} \\
    \le{}& \norm{f}_{L^{\infty}} \mesoscopicdifferenceseminorm{\Phi} + \norm{df}_{L^{\infty}} \left( \int_X \norm{\Phi(y)}^2 \left( \int_X \cutoffnabla(d_g(x,y)/a) \,\frac{\dvol_g(x)}{a^n} \right) \,\dvol_g(y) \right)^{1/2}
    \end{aligned}
  \end{equation}
  for $a \in (0,r_0)$.
  We note that
  \begin{equation}
    \int_X \cutoffnabla(d_g(x,y)/a) \,\frac{\dvol_g(x)}{a^n} \le \constantgeometric \int_{\R^n} \cutoffnabla(\abs{v}) \,dv =: (\constantmult)^2.
  \end{equation}
  Thus, we have
  \begin{equation}
    \mesoscopicdifferenceseminorm{f\Phi} \le \norm{f}_{L^{\infty}} \mesoscopicdifferenceseminorm{\Phi} + \constantmult \norm{df}_{L^{\infty}} \Ltwonorm{\Phi}.
  \end{equation}
\end{proof}

\begin{proof}[Second proof of Theorem~\ref{theorem: mesoscopic Rellich lemma}]
  Fix $r > 0$ as in Lemma~\ref{lemma: another comparison of difference seminorms}.
  Set
  \begin{equation}
    M := \sup_{j \ge 1} \bigl( \mesoscopicdifferenceseminorm[a_j]{\Phi_j} + \Ltwonorm{\Phi_j} \bigr) < \infty.
  \end{equation}
  Choose a partition of unity $f_1, \dots, f_N \colon X \to [0,1]$ such that $\sum f_{\alpha} = 1$ and
  \begin{equation}
    \support (f_{\alpha}) \subset B(p_{\alpha},r)
  \end{equation}
  for some points $p_1, \dots, p_N \in X$.
  Since $a_j \to 0$, there exists $j_0$ such that $a_j < r$ for every $j \ge j_0$.

  Fix $\alpha \in \{1,\dots,N\}$.
  We set
  \begin{equation}
    \Psi_{\alpha,j} := \flatmodel{(f_{\alpha}\Phi_j)} \in L^2(\flattorus,\flatmodel{E})
  \end{equation}
  for $j \ge j_0$.
  Each $\Psi_{\alpha,j}$ is supported in $B(0,r)$.
  By Lemmata~\ref{lemma: another comparison of difference seminorms} and~\ref{lemma: multiplication and difference seminorm}, we have
  \begin{equation}
    \begin{aligned}
    \mesoscopicseminorm[a_j]{\Psi_{\alpha,j}} \le{}& \constantcomparison \bigl( \mesoscopicdifferenceseminorm[a_j]{f_{\alpha} \Phi_j} + \Ltwonorm{f_{\alpha} \Phi_j} \bigr) \\
    \le{}& \constantcomparison \bigl( \norm{f_{\alpha}}_{L^{\infty}} \mesoscopicdifferenceseminorm[a_j]{\Phi_j} + \constantmult \norm{df_{\alpha}}_{L^{\infty}} \Ltwonorm{\Phi_j} + \norm{f_{\alpha}}_{L^{\infty}} \Ltwonorm{\Phi_j} \bigr) \\
    \le{}& \constantcomparison \bigl( \norm{f_{\alpha}}_{L^{\infty}} + \constantmult \norm{df_{\alpha}}_{L^{\infty}} + \norm{f_{\alpha}}_{L^{\infty}} \bigr) M.
    \end{aligned}
  \end{equation}
  By \eqref{inequality: comparison of L2 norms}, we have
  \begin{equation}
    \Ltwonorm[L^2(\flattorus)]{\Psi_{\alpha,j}} \le \constantLtwonorm \Ltwonorm[L^2(X)]{f_{\alpha} \Phi_j} \le \constantLtwonorm \norm{f_{\alpha}}_{L^{\infty}} M.
  \end{equation}
  Hence, we have
  \begin{equation}
    \sup_{j \ge j_0} \bigl( \mesoscopicseminorm[a_j]{\Psi_{\alpha,j}} + \Ltwonorm{\Psi_{\alpha,j}} \bigr) < \infty.
  \end{equation}
  Thus, by the flat case proved above, $\{\Psi_{\alpha,j}\}_j$ has a strongly convergent subsequence.
  Since there are only finitely many indices $\alpha$, we can choose a common subsequence, denoted by the same letters, along which all these sequences converge strongly in $L^2$.

  Finally, we have
  \begin{equation}
    \begin{aligned}
    \Ltwonorm[L^2(X)]{\Phi_j - \Phi_k} \le{}& \sum_{\alpha=1}^N \Ltwonorm[L^2(X)]{(f_{\alpha}\Phi_j) - (f_{\alpha}\Phi_k)} \\
    \le{}& \constantLtwonorm \sum_{\alpha=1}^N \Ltwonorm[L^2(\flattorus)]{\Psi_{\alpha,j} - \Psi_{\alpha,k}} \\
    \le{}& \constantLtwonorm N \max_{\alpha=1,\dots,N}\Ltwonorm[L^2(\flattorus)]{\Psi_{\alpha,j} - \Psi_{\alpha,k}}
    \end{aligned}
  \end{equation}
  for each $j,k \ge j_0$.
  Thus, $\{\Phi_j\}$ is Cauchy in $L^2(X,E)$.
  This completes the proof.
\end{proof}

\section{The coupled operators}\label{section: the coupled operators}

In this section, we couple the continuum Dirac operator with the mesoscopic Wilson--Dirac operator.
The main result is a vanishing theorem for the kernels of these coupled operators.
This theorem is the technical core of the paper.
In Section~\ref{section: index of mesoscopic Wilson Dirac operators}, we combine it with the relative $K^1$-theory and spectral flow recalled in Section~\ref{section: Spectral flow and the index} to prove the main theorem of the paper.
We keep the geometric data $(X,g,E,\gamma,\cliffordmultiplication,\nabla)$, the cut-off functions $\cutoffDirac$, $\cutoffWilson$, and $\cutoffnabla$, and the normalisation constants $C_D$ and $C_W$ fixed.

For $m\in[-1,1]$ and $t\in[0,1]$ with $0<a<\injectivityradius(X,g)/10$ as before, we consider the coupled operator
\begin{equation}
  \coupledoperator{m}{t} := \begin{pmatrix}
    D+m\gamma & t \id \\
    t \id & -\bigl(\mesoscopicWilsonDiracoperator + m\gamma\bigr)
  \end{pmatrix}
\end{equation}
on $L^2(X,E) \oplus L^2(X,E)$ with domain $L^2_1(X,E) \oplus L^2(X,E)$.
The coupled operator $\coupledoperator{m}{t}$ is self-adjoint by Proposition~\ref{proposition: mesoscopic WD is bounded self-adjoint}.

\begin{theorem}\label{theorem: vanishing theorem for coupled operators}
  There exists $\constantAvan \in (0,\injectivityradius(X,g)/10)$ such that, for all $a \in (0,\constantAvan)$, we have
  \begin{equation}
    \ker\coupledoperator{m}{t}=\{0\}
  \end{equation}
  for every $(m,t) \in [-1,1] \times [0,1]$ with either $m = \pm 1$ or $t = 1$.
\end{theorem}

 
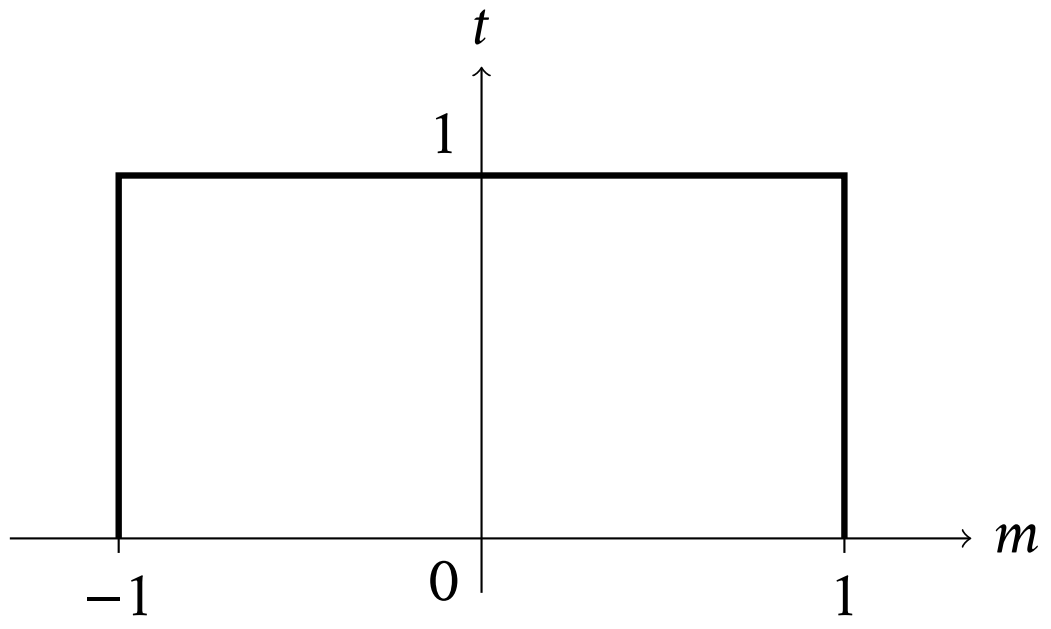
\begin{figure}[htbp]
  \centering
  \begin{tikzpicture}[x=2.4cm,y=2.4cm]
    \draw[->,thin] (-1.3,0) -- (1.35,0) node[right] {$m$};
    \draw[->,thin] (0,-0.15) -- (0,1.3) node[above] {$t$};

    \draw[very thick]
      (-1,0) -- (-1,1) -- (1,1) -- (1,0);

    \foreach \s in {-1,1}
      \draw[thin] (\s,0) -- (\s,-0.04) node[below] {$\s$};
    \node[below left] at (0,0) {$0$};
    \node[above left] at (0,1) {$1$};
  \end{tikzpicture}
\end{figure}

\begin{proof}
  We argue by contradiction.
  Suppose that there exist $a_j \to 0$, $(m_j,t_j) \in [-1,1] \times [0,1]$ with either $m_j = \pm 1$ or $t_j = 1$, and $(\Phi_j,\Psi_j) \in L^2_1(X,E) \oplus L^2(X,E)$ such that
  \begin{equation}
    \begin{pmatrix}
      D+m_j\gamma & t_j \id \\
      t_j \id & -\bigl(\mesoscopicWilsonDiracoperator[a_j] + m_j\gamma\bigr)
    \end{pmatrix}
    \begin{pmatrix}
      \Phi_j \\
      \Psi_j
    \end{pmatrix}
    =
    \begin{pmatrix}
      0 \\
      0
    \end{pmatrix}
  \end{equation}
  and
  \begin{equation}\label{equation: contradiction normalisation}
    \Ltwonorm{\Phi_j}^2 + \Ltwonorm{\Psi_j}^2 = 1
  \end{equation}
  for $j = 1,2,3,\dots$.
  Equivalently, $\Phi_j$ and $\Psi_j$ satisfy
  \begin{equation}\label{equation: contradiction DPhi equation}
    D \Phi_j = - m_j \gamma \Phi_j - t_j \Psi_j
  \end{equation}
  and
  \begin{equation}\label{equation: contradiction WDPsi equation}
    \mesoscopicWilsonDiracoperator[a_j] \Psi_j = -m_j \gamma \Psi_j + t_j\Phi_j.
  \end{equation}
  By compactness, after passing to a subsequence if necessary, we may assume that $(m_j,t_j) \to (m,t) \in [-1,1] \times [0,1]$ with either $m = \pm 1$ or $t = 1$.
  We note that $m^2 + t^2 \ge 1$.

  By \eqref{equation: contradiction DPhi equation} and \eqref{equation: contradiction normalisation}, we have
  \begin{equation}
    \Ltwonorm{D\Phi_j} = \Ltwonorm{m_j \gamma \Phi_j + t_j \Psi_j} \le \abs{m_j}\Ltwonorm{\Phi_j} + \abs{t_j}\Ltwonorm{\Psi_j} \le 2.
  \end{equation}
  Hence, the elliptic estimate for $D$ implies $\{\Phi_j\}$ is bounded in $L^2_1(X,E)$.
  Thus, by the ordinary Rellich lemma, we may pass to a further subsequence for which
  \begin{equation}\label{equation: strong convergence Phij}
    \Phi_j \to \Phi \text{ strongly in } L^2(X,E)
  \end{equation}
  and
  \begin{equation}\label{equation: weak convergence Phij}
    \Phi_j \rightharpoonup \Phi \text{ weakly in } L^2_1(X,E)
  \end{equation}
  for some $\Phi \in L^2_1(X,E)$.
  In the same way, by \eqref{equation: contradiction WDPsi equation} and \eqref{equation: contradiction normalisation}, we have
  \begin{equation}
    \Ltwonorm{\mesoscopicWilsonDiracoperator[a_j] \Psi_j} \le \Ltwonorm{\Psi_j} + \Ltwonorm{\Phi_j} \le 2,
  \end{equation}
  and, by Corollary~\ref{corollary: mesoscopic elliptic and Rellich}, we may pass to a further subsequence for which
  \begin{equation}\label{equation: strong convergence Psij}
    \Psi_j \to \Psi \text{ strongly in } L^2(X,E)
  \end{equation}
  and
  \begin{equation}\label{equation: weak convergence Psij}
    \mesoscopicWilsonDiracoperator[a_j]\Psi_j \rightharpoonup D\Psi \text{ weakly in } L^2(X,E)
  \end{equation}
  for some $\Psi \in L^2_1(X,E)$.

  By the strong convergence \eqref{equation: strong convergence Phij} and \eqref{equation: strong convergence Psij} and the normalisation \eqref{equation: contradiction normalisation}, we have
  \begin{equation}\label{equation: contradiction limit normalisation}
    \Ltwonorm{\Phi}^2 + \Ltwonorm{\Psi}^2 = 1.
  \end{equation}
  Taking weak limits in \eqref{equation: contradiction DPhi equation} and \eqref{equation: contradiction WDPsi equation}, we have
  \begin{equation}\label{equation: contradiction limit DPhi equation}
    (D + m\gamma) \Phi + t \Psi = 0
  \end{equation}
  and
  \begin{equation}\label{equation: contradiction limit WDPsi equation}
    (D + m\gamma) \Psi - t\Phi = 0.
  \end{equation}
  Since $D \gamma + \gamma D = 0$ and $D$ is self-adjoint, we have
  \begin{equation}
    \begin{aligned}
    0 ={}& \Ltwonorm{(D+m\gamma)\Phi + t\Psi}^2 + \Ltwonorm{(D+m\gamma)\Psi - t\Phi}^2 \\
    ={}& \Ltwonorm{D\Phi}^2 + \Ltwonorm{D\Psi}^2 + (m^2+t^2)(\Ltwonorm{\Phi}^2+\Ltwonorm{\Psi}^2) \\
    \ge{}&  (m^2+t^2)(\Ltwonorm{\Phi}^2+\Ltwonorm{\Psi}^2).
    \end{aligned}
  \end{equation}
  But $m^2+t^2 \ge 1$, contradicting \eqref{equation: contradiction limit normalisation}.
  This completes the proof.
\end{proof}

We next verify the Fredholm property.
Recall that
\begin{equation}
  m_W := \int_{\R^n} \frac{1}{\abs{u}} \frac{1}{C_W} \cutoffWilson(\abs{u}) \,du,
\end{equation}
which was introduced in Section~\ref{subsection: the mesoscopic Wilson term}.

\begin{proposition}\label{proposition: mesoscopicWD plus mgamma is Fredholm}
  If $a < m_W$, then $\mesoscopicWilsonDiracoperator + m\gamma$ is a Fredholm operator for every $m \in [-1,1]$.
\end{proposition}
\begin{proof}
  The decomposition $\mesoscopicWilsonterm = (m_W/a)\id + \mesoscopicWilsontermcompactpert$ gives
  \begin{equation}
    \mesoscopicWilsonDiracoperator + m\gamma = \left(\frac{m_W}{a} + m \right)\gamma + \mesoscopicDiracoperator + \mesoscopicWilsontermcompactpert \gamma.
  \end{equation}
  The last two terms are compact operators by Proposition~\ref{proposition: mesoscopic Dirac is hilbert--schmidt and selfadjoint} and Proposition~\ref{proposition: mesoscopic Wilson term is Fredholm and self-adjoint}.
  Since $m_W/a > 1$, the first term is invertible for every $m \in [-1,1]$.
  Thus, $\mesoscopicWilsonDiracoperator + m\gamma$ is Fredholm.
\end{proof}

\begin{proposition}\label{proposition: coupled operators are fredholm}
  If $a < m_W$, then $\coupledoperator[a]{m}{t}$ is a Fredholm operator for every $(m,t) \in [-1,1] \times [0,1]$.
\end{proposition}
\begin{proof}
  The operator $D+m\gamma$ is invertible for $m \ne 0$.
  Set
  \begin{equation}
    \Pi_m := \begin{cases}
      0 & \text{ if } m \ne 0 \\
      \text{the orthogonal projection onto $\ker D$} & \text{ if } m = 0.
    \end{cases}
  \end{equation}
  Then, $T_m := D + m\gamma + \Pi_m$ is invertible for every $m \in [-1,1]$.
  By the elliptic estimate, $T_m^{-1}$ is bounded from $L^2(X,E)$ to $L^2_1(X,E)$.
  Hence, by the ordinary Rellich lemma, it is compact as an operator on $L^2(X,E)$.

  For each $(m,t) \in [-1,1] \times [0,1]$, we have a factorisation
  \begin{equation}
    \coupledoperator{m}{t}
    +
    \begin{pmatrix}
      \Pi_m & 0 \\
      0 & 0
    \end{pmatrix}
    =
    \begin{pmatrix}
      \id & 0 \\
      tT_m^{-1} & \id
    \end{pmatrix}
    \begin{pmatrix}
      T_m & 0 \\
      0 & -\bigl( \mesoscopicWilsonDiracoperator + m\gamma \bigr) - t^2 T_m^{-1}
    \end{pmatrix}
    \begin{pmatrix}
      \id & tT_m^{-1} \\
      0 & \id
    \end{pmatrix}.
  \end{equation}
The rightmost and the leftmost factors are bounded and invertible.
The middle factor is Fredholm because $\mesoscopicWilsonDiracoperator + m\gamma$ is Fredholm by Proposition~\ref{proposition: mesoscopicWD plus mgamma is Fredholm} and $T_m^{-1}$ is compact.
Thus, the left-hand side is Fredholm.
Since $\Pi_m$ has finite rank, the coupled operator $\coupledoperator{m}{t}$ is Fredholm.
This completes the proof.
\end{proof}

 Recall that $\constantAvan$ is defined in Theorem~\ref{theorem: vanishing theorem for coupled operators}.
\begin{corollary}\label{corollary: invertibility at endpoints}
  If $a < \min\{\constantAvan,m_W\}$, the coupled operator $\coupledoperator{m}{t}$ is invertible for every $(m,t) \in [-1,1] \times [0,1]$ with either $m = \pm 1$ or $t=1$.
  In particular, the operators $\mesoscopicWilsonDiracoperator \pm \gamma$ are invertible.
\end{corollary}
\begin{proof}
  The coupled operator $\coupledoperator{m}{t}$ is self-adjoint.
  By Theorem~\ref{theorem: vanishing theorem for coupled operators}, its kernel vanishes for $a < \constantAvan$.
  By Proposition~\ref{proposition: coupled operators are fredholm}, it is Fredholm for $a < m_W$.
  Thus, it is invertible for $a < \min\{\constantAvan,m_W\}$.
\end{proof}

\section{Spectral flow and the index}\label{section: Spectral flow and the index}

In this section, we recall a model for the relative $K^1$-group and its relation to spectral flow.

\subsection{A model for the relative \texorpdfstring{$K^1$}{K1}-group}

We use a Hilbert bundle model for the relative $K^1$-group.
See \cite{MR5076866}*{Section~2.2} and \cite{afffmoy}*{Section~5}.

Let $Z$ be a compact Hausdorff space and $\mathcal{H} \to Z$ be a complex Hilbert bundle.
We adopt the following conventions for a Hilbert bundle.
The structure group is equipped with the operator norm topology.
The fibres are assumed to be separable and are allowed to be finite-dimensional.

A family $h = \{ h_z \}_{z \in Z}$ of densely defined closed possibly unbounded operators on the fibres of $\mathcal{H}$ is called \emph{Riesz continuous} if its bounded transforms
\begin{equation}
  h_z (1 + h_z^* h_z)^{-1/2}
\end{equation}
form a norm-continuous family in local unitary trivialisations.
The domains of the operators are allowed to vary with $z$.

Let $A \subset Z$ be a closed subspace.
A $K^1$-cycle $(\mathcal{H},h)$ over $(Z,A)$ consists of a complex Hilbert bundle $\mathcal{H} \to Z$ and a Riesz-continuous family $h = \{h_z\}_{z \in Z}$ of densely defined self-adjoint Fredholm operators on its fibres such that $h_z$ is invertible for $z \in A$.
Isomorphic cycles are identified.

\begin{definition}\label{definition: equivalence of K1 cycles}
  Two $K^1$-cycles $(\mathcal{H}_0,h_0)$ and $(\mathcal{H}_1,h_1)$ over $(Z,A)$ are equivalent if there exist a $K^1$-cycle $(\mathcal{H}_2,h_2)$ over $(Z,Z)$ and a family
  $\{\widetilde{h}_s\}_{s \in [0,1]}$ of densely defined self-adjoint Fredholm operators on $\mathcal{H}_0 \oplus \mathcal{H}_1 \oplus \mathcal{H}_2$ satisfying the following conditions:
  \begin{enumerate}
    \item The family is Riesz continuous in $(z,s)\in Z \times [0,1]$.
    \item We have $\widetilde{h}_0 = -h_0 \oplus h_1 \oplus h_2$.
    \item $(\widetilde{h}_s)_z$ is invertible whenever $z \in A$ or $s=1$.
  \end{enumerate}
  The set of equivalence classes is denoted by $K^1(Z,A)$.
\end{definition}

Direct sum gives the addition in $K^1(Z,A)$, and the inverse of $[(\mathcal{H},h)]$ is represented by $(\mathcal{H},-h)$.
A cycle that is invertible over all of $Z$ represents $0$.

We write $[h] \in K^1(Z,A)$ for the equivalence class of a cycle $(\mathcal{H},h)$ when the underlying Hilbert bundle is clear from the context.

The following direct consequence of the definition will be used in the next section.

\begin{proposition}\label{proposition: element equality criterion}
  Let $h_0$ and $h_1$ be Riesz-continuous families of self-adjoint Fredholm operators over $[-1,1]$ that are invertible over $\{\pm 1\}$.
  Suppose that $\{C(m,t)\}_{(m,t) \in [-1,1] \times [0,1]}$ is a Riesz-continuous family of densely defined self-adjoint Fredholm operators such that
  \begin{equation}
    C(m,0) = h_0(m) \oplus (-h_1(m))
  \end{equation}
  and $C(m,t)$ is invertible if either $m = \pm 1$ or $t = 1$.
  Then, we have
  \begin{equation}
    [h_0] = [h_1]
  \end{equation}
  in $K^1([-1,1],\{\pm 1\})$.
\end{proposition}

\begin{remark}
  The main theorem of this paper can also be proved directly using the additivity and homotopy invariance of spectral flow, without introducing $K$-theory.
  We formulate the argument in terms of relative $K^1$-groups with a view to future extensions, for example, to families of operators.
\end{remark}

We next recall the bounded perturbation criterion for Riesz continuity.
In~\cite{afffmoy}*{Proposition~21}, this is proved using the theory of operator-Lipschitz functions and Besov spaces.
For completeness, we give an elementary alternative proof based on the resolvent identity.
\begin{proposition}\label{proposition: bounded perturbations and Riesz continuity}
  Let $h_0$ and $h_1$ be self-adjoint operators on a complex Hilbert space with a common dense domain.
  Suppose that $h_1-h_0$ extends to a bounded operator.
  Then, we have
  \begin{equation}
    \operatornorm*{r(h_1) - r(h_0)} \le \operatornorm{h_1-h_0},
  \end{equation}
  where $r(\lambda) := \lambda/\sqrt{1+\lambda^2}$.
\end{proposition}
\begin{proof}
  For $R>0$, set
  \begin{equation}
    r_R(\lambda) := \frac{2}{\pi} \int_0^R \frac{\lambda}{1+\lambda^2+s^2} \,ds = \frac{2}{\pi} r(\lambda) \arctan \bigl(R/\sqrt{1+\lambda^2}\bigr).
  \end{equation}
  Then, $\abs{r_R(\lambda)} \le 1$ and $r_R(\lambda) \to r(\lambda)$ as $R \to \infty$.
  By the spectral theorem, $r_R(h_j)$ converges strongly to $r(h_j)$.

  We next prove an estimate for $r_R(h_1)-r_R(h_0)$ that is uniform in $R$.
  For a self-adjoint operator $h$ and $c>0$, set
  \begin{equation}
    Q_c(h) := h(h^2+c^2)^{-1} = \frac{(h-ic)^{-1}+(h+ic)^{-1}}{2}.
  \end{equation}
  Then, we have
  \begin{equation}
    r_R(h_j) = \frac{2}{\pi} \int_0^R Q_{\sqrt{1+s^2}}(h_j) \,ds.
  \end{equation}
  The resolvent identity gives
  \begin{equation}
    (h_1 \pm ic)^{-1}-(h_0 \pm ic)^{-1} = -(h_1 \pm ic)^{-1} (h_1-h_0) (h_0 \pm ic)^{-1}.
  \end{equation}
  Thus, we have
  \begin{equation}
    Q_c(h_1)-Q_c(h_0) = -\frac{1}{2} \bigl( (h_1-ic)^{-1} (h_1-h_0) (h_0-ic)^{-1} + (h_1+ic)^{-1} (h_1-h_0) (h_0+ic)^{-1} \bigr).
  \end{equation}
  Since $\operatornorm{(h_j \pm ic)^{-1}} \le c^{-1}$, we have
  \begin{equation}
    \operatornorm{Q_c(h_1)-Q_c(h_0)} \le \frac{\operatornorm{h_1-h_0}}{c^2}.
  \end{equation}
  Integrating this estimate with $c=\sqrt{1+s^2}$, we have
  \begin{equation}
    \begin{aligned}
    \operatornorm{r_R(h_1)-r_R(h_0)} \le{}& \frac{2}{\pi} \int_0^R \operatornorm{Q_{\sqrt{1+s^2}}(h_1) -Q_{\sqrt{1+s^2}}(h_0)}\,ds \\
    \le{}& \frac{2}{\pi} \operatornorm{h_1-h_0} \int_0^R \frac{ds}{1+s^2} \\
    \le{}& \operatornorm{h_1-h_0}.
    \end{aligned}
  \end{equation}

  Finally, for every vector $u$, strong convergence gives
  \begin{equation}
    \norm*{(r(h_1)-r(h_0))u} = \lim_{R \to \infty} \norm*{(r_R(h_1)-r_R(h_0))u} \le \operatornorm{h_1-h_0}\norm{u}.
  \end{equation}
  This completes the proof.
\end{proof}

\begin{corollary}\label{corollary: bounded perturbation of a fixed operator}
  Let $h$ be a densely defined self-adjoint operator on a complex Hilbert space.
  Let $\{K_z\}_{z \in Z}$ be a norm-continuous family of bounded self-adjoint operators on the same Hilbert space parametrised by a topological space $Z$.
  Then, the perturbed family $h_z := h + K_z$ is Riesz continuous.
\end{corollary}
\begin{proof}
   By Proposition~\ref{proposition: bounded perturbations and Riesz continuity}, we have
  \begin{equation}
    \operatornorm{r(h_z)-r(h_w)} \le \operatornorm{K_z-K_w}
  \end{equation}
  for $z,w \in Z$.
  Thus, the norm continuity of $z \mapsto K_z$ implies that of $z \mapsto r(h_z)$.
\end{proof}

\begin{corollary}\label{corollary: Riesz-continuous family}
  For each $a \in (0,\injectivityradius(X,g)/10)$, the families
  \begin{gather}
    m \mapsto D + m\gamma \\
    m \mapsto \mesoscopicWilsonDiracoperator + m\gamma \\
    (m,t) \mapsto \coupledoperator{m}{t}
  \end{gather}
  are Riesz continuous.
\end{corollary}

\subsection{Spectral flow}

For completeness, we recall the definition of the spectral flow for a Riesz-continuous family $T = \{T_m\}_{m \in [-1,1]}$ of self-adjoint Fredholm operators.
Set $F(m) := T_m (1+T_m^2)^{-1/2}$.
Choose a subdivision $-1 = m_0 < \dots < m_N = 1$ and $\epsilon_j > 0$ such that, for $m \in [m_{j-1},m_j]$, the points $\pm \epsilon_j$ lie outside the spectrum of $F(m)$ and the spectral projection onto $[-\epsilon_j,\epsilon_j]$ has finite rank.
Then, we define
\begin{equation}
  \spectralflow (T) := \sum_{j=1}^N \left( \dim \image 1_{[0,\epsilon_j)}(F(m_j)) - \dim \image 1_{[0,\epsilon_j)}(F(m_{j-1})) \right),
\end{equation}
where $1_{[0,\epsilon_j)}(F(m))$ denotes the spectral projection.
Such choices exist, and the resulting integer $\spectralflow (T)$ is independent of them.
See~\cite{MR2114489}*{Definition 1.1}.

For families that are invertible at the endpoints, spectral flow is additive under direct sum and invariant under homotopies that remain invertible at the endpoints.
It induces an isomorphism
\begin{equation}\label{equation: K1 and Z isomorphism}
  \spectralflow \colon K^1([-1,1],\{\pm 1\}) \to \Z.
\end{equation}
See~\cite{MR5076866}*{Theorem 2.36}.
The class of multiplication by $m \in [-1,1]$ on the trivial complex line bundle is mapped to $1$ by $\spectralflow$.

\subsection{Index and spectral flow}

The index of a continuum Dirac operator $D$ is defined by
\begin{equation}
  \indexAtiyahSinger D := \trace \bigl( \gamma|_{\kernel D} \bigr).
\end{equation}
Since $D$ anti-commutes with $\gamma$, it has the form
\begin{equation}
  D = \begin{pmatrix}
    0 & D^- \\
    D^+ & 0
  \end{pmatrix}
\end{equation}
with respect to the eigenbundle decomposition of $E = E^+ \oplus E^-$ by $\gamma$.
We have $(D^+)^* = D^-$ and
\begin{equation}
  \indexAtiyahSinger D = \dim \kernel D^+ - \dim \kernel D^- = \dim \kernel D^+ - \dim \cokernel D^+.
\end{equation}

\begin{proposition}\label{proposition: sf is equal to index}
  The family $\{D + m\gamma\}_{m \in [-1,1]}$ defines a class in $K^1([-1,1],\{\pm 1\})$, and
  \begin{equation}
    \spectralflow \bigl(\{D + m\gamma\}_{m \in [-1,1]}\bigr) = \indexAtiyahSinger D.
  \end{equation}
\end{proposition}
\begin{figure}[htbp]
  \centering
  \begin{tikzpicture}[
    x=3cm,y=1.65cm,
    font=\small,text=ink,
    line cap=round,line join=round
  ]
    \draw[ballline,dashed,line width=0.65pt]
      plot[domain=-1:1,samples=101]
      (\x,{sqrt(0.81+\x*\x)});
    \draw[ballline,dashed,line width=0.65pt]
      plot[domain=-1:1,samples=101]
      (\x,{-sqrt(0.81+\x*\x)});

    \node[above,text=ballline] at (0.68,1.40)
      {$\sqrt{\mu^2+m^2}$};
    \node[below,text=ballline] at (0.68,-1.40)
      {$-\sqrt{\mu^2+m^2}$};

    \draw[->,thin,ink] (-1.18,0) -- (1.23,0)
      node[right] {$m$};
    \draw[->,thin,ink] (0,-1.78) -- (0,1.78)
      node[above] {$\lambda$};

    \draw[thin,ink] (-1,0) -- (-1,-0.05)
      node[below] {$-1$};
    \draw[thin,ink] (1,0) -- (1,-0.05)
      node[below] {$1$};
    \node[left,inner sep=3pt] at (0,-0.24) {$0$};

    \draw[blueedge,line width=1pt]
      (-1,-1) -- (1,1);
    \draw[goldedge,line width=1pt]
      (-1,1) -- (1,-1);

    \draw[->,blueedge,line width=1pt]
      (-0.72,-0.72) -- (-0.48,-0.48);
    \draw[->,goldedge,line width=1pt]
      (-0.72,0.72) -- (-0.48,0.48);

    \node[below right,text=blueedge,inner sep=3pt]
      at (0.73,0.73) {$\lambda=m$};
    \node[above right,text=goldedge,inner sep=3pt]
      at (0.73,-0.73) {$\lambda=-m$};
  \end{tikzpicture}
\end{figure}
\begin{proof}
  By Corollary~\ref{corollary: Riesz-continuous family}, the family $\{D + m\gamma\}_{m \in [-1,1]}$ is Riesz continuous.
  For $\Phi \in L^2_1(X,E)$, the anti-commutation relation gives
  \begin{equation}
    \Ltwonorm{(D+m\gamma)\Phi}^2 = \Ltwonorm{D \Phi}^2 + m^2 \Ltwonorm{\Phi}^2.
  \end{equation}
  Hence, $D+m\gamma$ is invertible for $m \ne 0$.
  Thus, $\{D + m\gamma\}_{m \in [-1,1]}$ defines a class in $K^1([-1,1],\{\pm 1\})$.

  The orthogonal decomposition
  \begin{equation}
    L^2(X,E) = \kernel D \oplus (\kernel D)^{\perp}
  \end{equation}
  is preserved by $D$ and $\gamma$.
  There exists $C > 0$ such that
  \begin{equation}
    \Ltwonorm{D \Phi} \ge C \Ltwonorm{\Phi}
  \end{equation}
  for $\Phi \in L^2_1(X,E) \cap (\kernel D)^{\perp}$.
  Hence, the restriction of $D+m\gamma$ to $(\kernel D)^{\perp}$ is invertible for every $m$.
  Thus, this summand contributes $0$ to the spectral flow.

  On $\kernel D$, the family is $m\gamma$.
  The family $m\gamma$ acts by multiplication by $m$ and $-m$ on $\kernel D^+$ and $\kernel D^-$ respectively.
  Thus, the total contribution to the spectral flow is equal to $\dim \kernel D^+ - \dim \kernel D^- = \indexAtiyahSinger D$.
  This completes the proof.
\end{proof}

\begin{remark}
  The suspension isomorphism
  \begin{equation}
    K^0(\mathrm{pt}) \cong K^1([-1,1],\{\pm 1\})
  \end{equation}
  sends the class of the cycle defined by $D$ to the class defined by $\{D+m\gamma\}_{m\in[-1,1]}$.
  See \cite{MR5076866}*{Sections~2.2.2--2.2.3}.
\end{remark}

\begin{remark}
  The relation between the index and spectral flow in Proposition~\ref{proposition: sf is equal to index} also admits an interpretation in terms of domain-wall fermion Dirac operators.
  See~\cite{MR4179728}.
\end{remark}

\section{Index of mesoscopic Wilson--Dirac operators}\label{section: index of mesoscopic Wilson Dirac operators}

In this section, we prove our main theorem.
Since the mesoscopic Wilson--Dirac operator is self-adjoint, its ordinary Fredholm index is zero.
Since it does not anti-commute with $\gamma$, we cannot define its index as the trace of $\gamma$ on its kernel either.
See Theorem~\ref{theorem: no-go theorem}.
We thus define its index by the spectral flow of the family $\{\mesoscopicWilsonDiracoperator + m\gamma\}_{m\in[-1,1]}$ for sufficiently small $a$.
We then prove that this index agrees with $\indexAtiyahSinger D$.

By Proposition~\ref{proposition: mesoscopicWD plus mgamma is Fredholm} and Corollary~\ref{corollary: Riesz-continuous family}, if $a < m_W$, then $\{\mesoscopicWilsonDiracoperator + m\gamma\}_{m \in [-1,1]}$ is a Riesz-continuous family of self-adjoint Fredholm operators.

\begin{definition}\label{definition: mesoscopic index}
  We define the \emph{mesoscopic index} by
  \begin{equation}
    \indexAtiyahSinger \bigl( \mesoscopicWilsonDiracoperator \bigr) := \spectralflow \left( \bigl\{ \mesoscopicWilsonDiracoperator + m\gamma \bigr\}_{m \in [-1,1]} \right)
  \end{equation}
  for $a \in (0,\min\{m_W, \injectivityradius(X,g)/10\})$.
\end{definition}

\begin{theorem}\label{theorem: main theorem}
  There exists $\constantAindex > 0$ such that
  \begin{equation}
    \indexAtiyahSinger D = \indexAtiyahSinger \bigl( \mesoscopicWilsonDiracoperator \bigr)
  \end{equation}
  for $a \in (0,\constantAindex)$.
\end{theorem}
\begin{proof}
  Set $\constantAindex := \min \{\constantAvan, m_W\}$, where $\constantAvan$ is defined in Theorem~\ref{theorem: vanishing theorem for coupled operators}.
  Fix $a \in (0,\constantAindex)$.
  By Proposition~\ref{proposition: coupled operators are fredholm} and Corollary~\ref{corollary: Riesz-continuous family}, the coupled operators $\coupledoperator{m}{t}$ for $m \in [-1,1]$ and $t \in [0,1]$ form a Riesz-continuous family of self-adjoint Fredholm operators.
  By Corollary~\ref{corollary: invertibility at endpoints}, this family is invertible at either $m = \pm 1$ or $t = 1$.
  At $t = 0$, we have
  \begin{equation}
    \coupledoperator{m}{0} = (D+m\gamma) \oplus \bigl(- \bigl(\mesoscopicWilsonDiracoperator + m\gamma \bigr) \bigr).
  \end{equation}
  Hence, by Proposition~\ref{proposition: element equality criterion}, we have
  \begin{equation}
    [D+m\gamma] = \bigl[ \mesoscopicWilsonDiracoperator + m\gamma \bigr].
  \end{equation}
  Then, by \eqref{equation: K1 and Z isomorphism}, Proposition~\ref{proposition: sf is equal to index}, and Definition~\ref{definition: mesoscopic index}, we have
  \begin{equation}
    \indexAtiyahSinger D = \indexAtiyahSinger \bigl( \mesoscopicWilsonDiracoperator \bigr)
  \end{equation}
  for $a \in (0,\constantAindex)$.
  This completes the proof.
\end{proof}

\begin{remark}
  We can replace the interval $[-1,1]$ by any fixed interval $[-M,M]$ with $M > 0$.
  Indeed, in the proof of Theorem~\ref{theorem: vanishing theorem for coupled operators}, the bounds $2$ are replaced by $M+1$, and the final estimate becomes $m^2 + t^2 \ge \min \{M^2, 1\}$.
  The Fredholm arguments apply for $m_W/a > M$.
  Thus, the threshold for $a$ depends on $M$.
\end{remark}

\section*{Acknowledgements}
This paper is part of a joint project with Mikio Furuta, who passed away before it was completed.
The second author prepared the final manuscript and is responsible for any errors.
We are grateful to Hidenori Fukaya, Tsuyoshi Kato, Tetsuya Onogi, and Satoshi Yamaguchi for helpful discussions.
This work was supported by JSPS KAKENHI Grant Number JP26K06778.

\section*{Use of generative AI}
The authors used ChatGPT and Claude as aids in producing the TikZ code for the figures, in checking the proofs, and in revising the exposition and the English of the manuscript.
Theorem~\ref{theorem: no-go theorem}) and Corollary~\ref{corollary: norm-resolvent convergence} were suggested by Claude in discussions on the manuscript.
The second author verified these drafts and rewrote them.
The other results and their proofs are due to the authors.
All suggestions were examined by the second author, who made the final decisions and takes full responsibility for the content.

\end{document}